\documentclass[12pt,reqno]{NumPDEsArticle}

\usepackage[T1]{fontenc}
\usepackage[english]{babel}
\usepackage{csquotes}
\usepackage{enumitem}
\usepackage{amsmath,amssymb}
\usepackage{booktabs}
\usepackage{biblatex}
\usepackage{stmaryrd}
\usepackage{diagbox}
\usepackage[caption=false]{subfig}
\definecolor{best}{HTML}{2ca02c}
\definecolor{row}{HTML}{1f77b4}
\definecolor{col}{HTML}{bcbd22}
\usepackage{color, colortbl}
\usepackage{moreverb}
\usepackage[mode=buildnew]{standalone}

\usepackage{NumPDEsMacros}

\usepackage{hyperref}
\usepackage{tikz}
\usetikzlibrary{calc,matrix}
\usepackage{pgfplots}
\usepackage{pgfplotstable}
\pgfplotsset{%
    compat=1.18,%
    every axis/.style={scale only axis},%
    grid style={densely dotted, semithick},%
}
\pgfplotstableset{%
    include outfiles=true,%
    force remake=true,%
    trim cells=true,%
    col sep=comma,%
    row sep=newline%
}
\usepackage{booktabs}
\usepackage{colortbl}
\usepackage{diagbox}

\definecolor{myyellow}{HTML}{bcbd22}
\colorlet{colorcolmin}{myyellow!40!white}
\definecolor{myblue}{HTML}{1f77b4}
\colorlet{colorrowmin}{myblue!40!white}
\definecolor{mygreen}{HTML}{2ca02c}
\colorlet{colorbothmin}{mygreen!40!white}

\pgfplotsset{
    colormap={parula}{
        rgb=(0.2081,0.1663,0.5292)
        rgb=(0.2116,0.1898,0.5777)
        rgb=(0.2123,0.2138,0.627)
        rgb=(0.2081,0.2386,0.6771)
        rgb=(0.1959,0.2645,0.7279)
        rgb=(0.1707,0.2919,0.7792)
        rgb=(0.1253,0.3242,0.8303)
        rgb=(0.0591,0.3598,0.8683)
        rgb=(0.0117,0.3875,0.882)
        rgb=(0.006,0.4086,0.8828)
        rgb=(0.0165,0.4266,0.8786)
        rgb=(0.0329,0.443,0.872)
        rgb=(0.0498,0.4586,0.8641)
        rgb=(0.0629,0.4737,0.8554)
        rgb=(0.0723,0.4887,0.8467)
        rgb=(0.0779,0.504,0.8384)
        rgb=(0.0793,0.52,0.8312)
        rgb=(0.0749,0.5375,0.8263)
        rgb=(0.0641,0.557,0.824)
        rgb=(0.0488,0.5772,0.8228)
        rgb=(0.0343,0.5966,0.8199)
        rgb=(0.0265,0.6137,0.8135)
        rgb=(0.0239,0.6287,0.8038)
        rgb=(0.0231,0.6418,0.7913)
        rgb=(0.0228,0.6535,0.7768)
        rgb=(0.0267,0.6642,0.7607)
        rgb=(0.0384,0.6743,0.7436)
        rgb=(0.059,0.6838,0.7254)
        rgb=(0.0843,0.6928,0.7062)
        rgb=(0.1133,0.7015,0.6859)
        rgb=(0.1453,0.7098,0.6646)
        rgb=(0.1801,0.7177,0.6424)
        rgb=(0.2178,0.725,0.6193)
        rgb=(0.2586,0.7317,0.5954)
        rgb=(0.3022,0.7376,0.5712)
        rgb=(0.3482,0.7424,0.5473)
        rgb=(0.3953,0.7459,0.5244)
        rgb=(0.442,0.7481,0.5033)
        rgb=(0.4871,0.7491,0.484)
        rgb=(0.53,0.7491,0.4661)
        rgb=(0.5709,0.7485,0.4494)
        rgb=(0.6099,0.7473,0.4337)
        rgb=(0.6473,0.7456,0.4188)
        rgb=(0.6834,0.7435,0.4044)
        rgb=(0.7184,0.7411,0.3905)
        rgb=(0.7525,0.7384,0.3768)
        rgb=(0.7858,0.7356,0.3633)
        rgb=(0.8185,0.7327,0.3498)
        rgb=(0.8507,0.7299,0.336)
        rgb=(0.8824,0.7274,0.3217)
        rgb=(0.9139,0.7258,0.3063)
        rgb=(0.945,0.7261,0.2886)
        rgb=(0.9739,0.7314,0.2666)
        rgb=(0.9938,0.7455,0.2403)
        rgb=(0.999,0.7653,0.2164)
        rgb=(0.9955,0.7861,0.1967)
        rgb=(0.988,0.8066,0.1794)
        rgb=(0.9789,0.8271,0.1633)
        rgb=(0.9697,0.8481,0.1475)
        rgb=(0.9626,0.8705,0.1309)
        rgb=(0.9589,0.8949,0.1132)
        rgb=(0.9598,0.9218,0.0948)
        rgb=(0.9661,0.9514,0.0755)
        rgb=(0.9763,0.9831,0.0538)
    }
}

\definecolor{TUblue}{rgb}{0,0.4,0.6}
\definecolor{TUgray}{rgb}{0.3922,0.3882,0.3882}
\definecolor{TUgreen}{rgb}{0,0.4941,0.4431}
\definecolor{TUmagenta}{rgb}{0.7294,0.2745,0.5098}
\definecolor{TUyellow}{rgb}{0.8824,0.5373,0.1333}

\usepackage{pdftexcmds}
\makeatletter
\newcommand{\strequal}[2]{\pdf@strcmp{#1}{#2}==0}
\makeatother

\newcommand\drawslopetriangle[4][ST]{
  \pgfplotsextra
  {
    \pgfkeys{/pgf/fpu=true}
    \pgfmathsetmacro\leftcoord{#3};
    \pgfmathsetmacro\rightcoord{10*#3};
    \pgfmathsetmacro\bottomcoord{#4};
    \pgfmathsetmacro\topcoord{10^(#2)*#4};
    \pgfkeys{/pgf/fpu=false}

    \coordinate (#1-BL) at (axis cs:\leftcoord,\bottomcoord);
    \coordinate (#1-BR) at (axis cs:\rightcoord,\bottomcoord);
    \coordinate (#1-TL) at (axis cs:\leftcoord,\topcoord);

    \shadedraw[%
      bottom color = black!20,%
      middle color = black!5,%
      top color    = white,%
      draw         = black%
    ]
      (#1-TL) -- (#1-BL) node[midway, left=-2pt] {\scriptsize\(#2\)}
      -- (#1-BR) node[midway, below=-2pt] {\scriptsize\(1\)} -- (#1-TL);
  }
}
\newcommand\drawswappedslopetriangle[4][SST]{
  \pgfplotsextra
  {
    \pgfkeys{/pgf/fpu=true}
    \pgfmathsetmacro\leftcoord{#3/10};
    \pgfmathsetmacro\rightcoord{#3};
    \pgfmathsetmacro\topcoord{#4};
    \pgfmathsetmacro\bottomcoord{10^(-#2)*#4};
    \pgfkeys{/pgf/fpu=false}

    \coordinate (#1-TR) at (axis cs:\rightcoord,\topcoord);
    \coordinate (#1-BR) at (axis cs:\rightcoord,\bottomcoord);
    \coordinate (#1-TL) at (axis cs:\leftcoord,\topcoord);

    \shadedraw[%
      bottom color = black!20,%
      middle color = black!5,%
      top color    = white,%
      draw         = black%
    ]
      (#1-BR) -- (#1-TR) node[midway, right=-2pt] {\scriptsize\(#2\)}
      -- (#1-TL) node[midway, above=-2pt] {\scriptsize\(1\)} -- (#1-BR);
  }
}
\newcommand\drawslopetriangleup[4][STU]{
  \pgfplotsextra
  {
    \pgfkeys{/pgf/fpu=true}
    \pgfmathsetmacro\leftcoord{#3};
    \pgfmathsetmacro\rightcoord{10*#3};
    \pgfmathsetmacro\bottomcoord{#4};
    \pgfmathsetmacro\topcoord{10^(#2)*#4};
    \pgfkeys{/pgf/fpu=false}

    \coordinate (#1-BL) at (axis cs:\leftcoord,\bottomcoord);
    \coordinate (#1-BR) at (axis cs:\rightcoord,\bottomcoord);
    \coordinate (#1-TR) at (axis cs:\rightcoord,\topcoord);

    \shadedraw[%
      bottom color    = black!20,%
      middle color = black!5,%
      top color = white,%
      draw         = black%
    ]
      (#1-BL) -- (#1-BR) node[midway, below=-2pt] {\scriptsize\(1\)}
      -- (#1-TR) node[midway, right=-2pt] {\scriptsize\(#2\)} -- (#1-BL);
  }
}

\newcommand\vvvert{|\mkern-1.5mu|\mkern-1.5mu|}
\newcommand\lbracket{[\mkern-4.5mu[}
\newcommand\rbracket{]\mkern-4.5mu]}

\newcommand\Eta{\textup{H}}
\newcommand\Zeta{\textup{Z}}

\newcommand\A{\mathbb{A}}
\newcommand\N{\mathbb{N}}
\newcommand\R{\mathbb{R}}
\newcommand\T{\mathbb{T}}
\newcommand\dist{\textbf{\textup{d}}}

\newcommand\QQ{\mathcal{Q}}
\newcommand\RR{\mathcal{R}}
\newcommand\TT{\mathcal{T}}
\newcommand\MM{\mathcal{M}}

\newcommand\UU{\mathcal{U}}
\newcommand\XX{\mathcal{X}}

\newcommand\elll{{\underline{\ell}}}
\newcommand\kk{{\underline{k}}}
\newcommand\jj{{\underline{j}}}
\newcommand\mm{{\underline{m}}}
\newcommand\mmu{{\underline{\mu}}}

\newcommand\lambdaalg{\lambda_{\textup{alg}}}
\newcommand\lambdasym{\lambda_{\textup{sym}}}
\newcommand\lambdalin{\lambda_{\textup{lin}}}
\newcommand\qalg{q_{\textup{alg}}}
\newcommand\qsym{q_{\textup{sym}}}
\newcommand\qlin{q_{\textup{lin}}}
\newcommand\qred{q_{\textup{red}}}
\newcommand\qctr{q_{\textup{ctr}}}

\newcommand\Clin{C_{\textup{lin}}}
\newcommand\Crel{C_{\textup{rel}}}
\newcommand\Cdrel{C_{\textup{drel}}}
\newcommand\Cstab{C_{\textup{stab}}}

\let\div\relax
\DeclareMathOperator\div{div}
\DeclareMathOperator*\essinf{ess\,inf\vphantom{p}}
\DeclareMathOperator*\esssup{ess\,sup}
\DeclareMathOperator\cost{cost}
\DeclareMathOperator\refine{\tt refine}

\renewcommand{\d}[1]{\mathop{}\!\mathrm{d}#1}

\makeatletter
\newcommand{\labeltext}[2]{%
    \@bsphack
    \csname phantomsection\endcsname 
    \def\@currentlabel{#1}{\label{#2}}%
    \@esphack
}
\makeatother

\title[Iterative solvers in adaptive FEM]{\Large Iterative solvers
in adaptive FEM \\\bigskip \normalsize\mdseries
Adaptivity yields quasi-optimal computational runtime}



\author{Philipp Bringmann}
\author{Ani Miraçi}
\author{Dirk Praetorius}

\email{philipp.bringmann@asc.tuwien.ac.at (corresponding author)}
\email{ani.miraci@asc.tuwien.ac.at}
\email{dirk.praetorius@asc.tuwien.ac.at}

\thanks{This research was funded by the Austrian Science Fund (FWF) projects 
\href{https://www.fwf.ac.at/en/research-radar/10.55776/F65}{10.55776/F65} 
(SFB F65 ``Taming complexity in PDE systems''), 
\href{https://www.fwf.ac.at/en/research-radar/10.55776/I6802}{10.55776/I6802} 
(international project I6802 ``Functional error estimates for PDEs on 
unbounded domains''), and 
\href{https://www.fwf.ac.at/en/research-radar/10.55776/P33216}{10.55776/P3\-3216} 
(standalone project P33216 ``Computational nonlinear PDEs'').}

\begin{document}
\maketitle

\begin{abstract}
	This chapter provides an overview of state-of-the-art 
  adaptive finite element methods (AFEMs) for the numerical solution of 
  second-order elliptic partial differential equations (PDEs), where the 
  primary focus is on the optimal interplay of local mesh refinement and 
  iterative solution of the arising discrete systems. Particular emphasis 
  is placed on the thorough description of the essential ingredients 
  necessary to design adaptive algorithms of optimal complexity, i.e., 
  algorithms that mathematically guarantee the optimal rate of convergence 
  with respect to the overall computational cost and, hence, time. Crucially, 
  adaptivity induces reliability of the computed numerical approximations by 
  means of a-posteriori error control. This ensures that the error committed 
  by the numerical scheme is bounded from above by computable quantities. 
  The analysis of the adaptive algorithms is based on the study of appropriate 
  quasi-error quantities that include and balance different components of the 
  overall error. Importantly, the quasi-errors stemming from an adaptive 
  algorithm with contractive iterative solver satisfy a centerpiece concept, 
  namely, full R-linear convergence. This guarantees that the adaptive 
  algorithm is essentially contracting this quasi-error at each step and it 
  turns out to be the cornerstone for the optimal complexity of AFEM. 
  The unified analysis of the adaptive algorithms is presented in the context 
  of symmetric linear PDEs. Extensions to goal-oriented, non-symmetric, 
  as well as non-linear PDEs are presented with suitable nested iterative 
  solvers fitting into the general analytical framework of the linear 
  symmetric case. Numerical experiments highlight the theoretical results and 
  emphasize the practical relevance and gain of adaptivity with iterative 
  solvers for numerical simulations with optimal complexity.
\end{abstract}

%
%

\section{Motivation}
\label{BMP:sec:intro}

\subsection{Concept of adaptivity and historical overview}
\label{BMP:sec:concept}

The design of modern numerical schemes
for approximating the solution
to partial differential equations (PDEs)
proceeds in certain steps.
First, the equation is transformed into a suitable
\emph{variational formulation}.
This provides the mathematical framework allowing to verify
existence and uniqueness of the solution.
Second, the \emph{discretization} of the formulation
is typically based on dissecting the underlying
domain into a computational mesh, choosing an appropriate
finite element space, e.g.,
consisting of globally continuous piecewise polynomials of a
given degree, and solving the discretized variational formulation
for this discrete space.
Such a scheme is called finite element method (FEM)
and results in a system of equations with a finite
number of unknowns to be solved on a computer.
However, generating and refining the mesh is
crucial for the quality of the discrete solution in
approximating the unavailable -- and possibly singular --
solution to the PDE.
Indeed, the accuracy of the approximation will be heavily spoiled by
singularities of the solution, if they are not sufficiently
resolved by the mesh.
Third, \emph{iterative solvers}
are of utmost importance
to compute a numerical approximation.
In particular, once the discrete problem is translated to a
system of linear equations, the choice of the algebraic
solver is essential for memory- and time-efficient simulations.
Moreover, non-linear discrete problems must be
linearized and, clearly, the interaction of linearization and inexact solution
of the resulting systems of linear equations
pose additional challenges to the structure of the numerical scheme and its 
analysis.

From a practical point of view,
numerical simulations must not only
ensure \emph{reliability}
(i.e., a guaranteed user-specified accuracy)
but should moreover provide
\emph{efficiency} (i.e., considerably short and, potentially, 
minimal computing time).
For this purpose,
computable \emph{a-posteriori error estimators}
play a major role.
They are evaluated after computing a numerical approximation
to the solution with the aim of assessing its quality and,
thereby, deciding whether a desired accuracy has been achieved or not.
Moreover, they enable the efficient usage
of computational resources 
(like memory consumption and computing time) by steering the
adaptive mesh-refinement process
and balancing different parts of the error like
discretization error, linearization error, and algebraic
solver error.
The newly refined mesh allows to compute a more accurate
approximation which in turn, through the a-posteriori
estimator, allows to further refine the mesh where needed.
The continuation of this \emph{nested iterative process},
while yielding ever better numerical results, is inherently
cumulative in nature.
On the one hand,
these algorithms allow for the best possible
accuracy with respect to the
dimension of the discrete spaces of approximation,
i.e.,~\emph{rate-optimality}.
On the other hand,
the additional computational effort for the mesh refinement
should pay off to achieve the desired error threshold
within the shortest possible runtime
and, thus, keeping the computational cost minimal,
i.e.,~\emph{cost-optimality}.
In both cases,
the analysis is not built upon the study of the error itself,
but on (essentially equivalent) \emph{quasi-error} quantities
consisting of 
the a-posteriori estimator and other error components.

In order to achieve such optimality results,
the need for suitable iterative solvers and appropriate stopping 
criteria balancing
reasonable computational cost and good accuracy
becomes centerpiece.
In practice, the approximation should
be constructed in a number of operations proportional to the
dimension of the discrete space
(resp. number of simplices in the triangulation),
i.e., at~\textit{linear cost}, while ensuring
an \textit{error contraction} property with respect to the
previous available approximation.
Finally, the overall number of iterative solver steps
is controlled thanks to the \textit{a-posteriori} estimator
through an adaptive \textit{stopping criterion}.

Before introducing the model problem
in Section~\ref{BMP:sec:model_problem},
this final paragraph recalls key contributions from
the academic history of the field, where we refer, in particular, to the 
exhaustive review~\cite{bbco2024}.
Some first pioneering works on the mathematics of adaptive mesh refinement
include~\cite{br1978,zr1979,gkzb1983}.
Since adaptive algorithms do not drive the mesh size to zero
in the whole domain, their convergence is not guaranteed by
the classical a-priori error analysis, see, e.g., 
\cite{BrennerScott_08, ErnGuermondFEM2}.
While convergence was already addressed in~\cite{bv1984} for 
one-dimensional problems, the first convergence results of adaptive FEMs 
for two or higher space dimensions
were provided by
\cite{doerfler1996,mns2000}.
Eventually, the notion of nonlinear approximation classes
enabled the proof of optimal convergence rates
with respect to the dimension of the discrete space, 
see~\cite{bdd2004,stevenson2007}.
The succeeding generalizations to a myriad of applications
led to a summarizing axiomatic framework in~\cite{cfpp2014} and 
the role of adaptive stopping criteria being emphasized 
in~\cite{ErnVohralik2013}.
While even the early contributions like~\cite{stevenson2007}
address the iterative algebraic solution,
all convergence results remain subject to the assumption of
the iterative solution being sufficiently close to the exact
algebraic solution; see, e.g.,~\cite{stevenson2007,cfpp2014}.
This severe algorithmic restriction was removed
for the proof of full R-linear convergence
in~\cite{banach2} for energy-contractive iterative solvers
and in~\cite{bfmps2023} for general contractive solvers.

\subsection{Symmetric second-order model problem}
\label{BMP:sec:model_problem}

The focus of the present chapter
is initially drawn to symmetric second-order linear elliptic PDEs
as a model problem with further extensions,
discussed at the end of the chapter.
Let \(\Omega \subset \R^{d}\), \(d \in \N\), denote the
computational domain, which is mathematically assumed to be
a bounded Lipschitz domain with polytopal boundary
\(\partial\Omega\).
This chapter employs the notation of Lebesgue and Sobolev
spaces as used in standard textbooks on PDEs,
e.g.,~\cite[Sect.~5]{evans1998}.
Given a symmetric coefficient tensor
\(\boldsymbol{A} \in [L^\infty(\Omega)]^{d \times d}_{\textup{sym}}\)
and load functions $f \in L^2(\Omega)$
and \(\boldsymbol{f} \in [L^2(\Omega)]^d\),
we consider the elliptic boundary value problem of
seeking \(u^\star \in \XX \coloneqq H^1_0(\Omega)\)
such that
\begin{equation}
    \label{BMP:eq:strong}
    -\div (\boldsymbol{A}\nabla u^\star)
    =
    f - \div \boldsymbol{f}
    \text{ in } \Omega
    \quad\text{and}\quad
    u^\star
    =
    0
    \text{ on } \Gamma \coloneqq \partial\Omega.
\end{equation}
We assume that \(\boldsymbol{A}\) is uniformly elliptic in the sense
that the minimal and maximal eigenvalues are uniformly
bounded almost everywhere, i.e.,
\begin{equation}
    \label{BMP:eq:ellipticity}
    0
    <
    \essinf_{x \in \Omega}
    \lambda_{\textup{min}}(\boldsymbol{A}(x))
    \leq
    \esssup_{x \in \Omega}
    \lambda_{\textup{max}}(\boldsymbol{A}(x))
    <
    \infty.
\end{equation}
Then, we define the bilinear form \(a : \XX \times \XX \to \R\) 
using the dot notation henceforth designating the vector-product in the 
sense of
\[
    a(u, v)
    \coloneqq
    \int_\Omega (\boldsymbol{A} \nabla u) \cdot \nabla v \, \d{x}
    \ 
    =
    \sum_{i,j=1}^d \int_\Omega \boldsymbol{A}_{ij} 
    \partial_i u \partial_j v \, \d{x}   
    \text{ for all } u, v \in \XX = H^1_0(\Omega).
\]
The uniform bounds~\eqref{BMP:eq:ellipticity} ensure that
\(a(\cdot,\cdot)\) is a scalar product on \(\XX\) and thus induces an 
energy norm denoted as \(\vvvert \cdot \vvvert\).
The Riesz--Fischer theorem~\cite[Thm.~2 in \S{} D.2]{evans1998}
guarantees existence and uniqueness of the solution $u^\star \in \XX$
to~\eqref{BMP:eq:strong} and the equivalent variational formulation
\begin{equation}
    \label{BMP:eq:weak}
    a(u^\star, v)
    =
    \!\int_\Omega (f \, v + \boldsymbol{f} \cdot \nabla v) \, \d{x} \
    = 
    \! \int_\Omega f \, v  \, \d{x}
    + \sum_{i=1}^d
    \!\int_\Omega \boldsymbol{f}_i \partial_i v \, \d{x} 
    \eqqcolon
    F(v)
    \text{ for all } v \in \XX.
\end{equation}
The discretization of problem~\eqref{BMP:eq:weak}
employs a conforming simplicial mesh $\TT_H$ of $\Omega$,
i.e., triangles for \(d = 2\) and
tetrahedra for \(d = 3\).
Conformity means that
$\bigcup_{T \in \TT_H} T = \overline \Omega$ and
that the intersection of two different simplices
is either empty or a common lower-dimensional sub-simplex,
e.g., a common vertex or edge for \(d = 2\) or a common
vertex, edge, or face for \(d = 3\).
The usual conforming finite element space
of globally continuous and piecewise polynomial functions
of degree at most \(p \in \N\) reads
\begin{align}
    \label{BMP:eq:fem-space}
    \! \! \XX_H
    \coloneqq
    \{
        v_H \in\XX
        \;:\;
        v_H\vert_T 
        \text{ is a polynomial of degree $\le p$, for all $T \in \TT_H$} 
    \}.
\end{align}
This leads to the discrete formulation seeking
the unique solution $u_H^\star \in \XX_H$ to
\begin{align}
    \label{BMP:eq:discrete}
    a(u_H^\star, v_H)
    =
    F(v_H)
    \quad \text{for all }
    v_H \in \XX_H.
\end{align}
Although the index \(H\) historically originates from meshes
with a quasi-uniform mesh size \(H\),
this notation indicates all discrete objects throughout this chapter
also with respect to highly graded meshes.
Related discrete objects have the same index, e.g., 
\(v_H\) is a discrete function in the space \(\XX_H\) corresponding 
to the triangulation \(\TT_H\). 
Different indices are used to distinguish between different triangulations 
and the corresponding discrete objects, e.g., 
\(\TT_H\), \(\TT_h\), \(\TT_\ell\), etc.

\subsection{Outline of the chapter}
\label{BMP:sec:outline}

In the following, the algorithm and its constituting modules are presented in
Section~\ref{BMP:sec:algo}. 
The main results on optimality
are showcased in Section~\ref{BMP:sec:convergence}. 
Section~\ref{BMP:sec:numerics} highlights the numerical
performance of the presented adaptive algorithm for a
variety of test cases and confirms its optimality results in
accordance with the theory.
The proofs of the main results are 
presented in detail in Section~\ref{BMP:sec:proofs}.
The Sections~\ref{BMP:sec:goafem}--\ref{BMP:sec:ailfem} are
devoted to extending the model problem~\eqref{BMP:eq:strong} to a
goal-oriented framework, or a setting exhibiting 
non-symmetry or non-linearity.
Throughout, particular emphasis is put on the algorithmic decision 
for the termination of the iterative solver and its cost-optimal interplay 
with adaptive mesh refinement. An additional focus is on parameter-robust 
convergence, where mathematics guarantees convergence without any constraints 
on the choice of fine-tuned parameters.

%
%

\section{Adaptive mesh-refining algorithm}
\label{BMP:sec:algo}

This section discusses the assumptions on the
modules \texttt{SOLVE}, \texttt{ESTIMATE}, \texttt{MARK},
and \texttt{REFINE} as employed in the adaptive algorithm.
The module \texttt{REFINE} is presented first
since the successively improved sequence of meshes is fundamental
for the discretization and computation in the
remaining modules.
The formulation of the algorithm follows in
Subsection~\ref{BMP:sec:AFEM}
as well as the introduction
of the notions of quasi-errors and
computational cost in Subsection~\ref{BMP:sec:quasi-errors}
below.
Subsection~\ref{BMP:sec:convergence}
presents the convergence and complexity results.

\subsection{Module \texttt{REFINE}}
\label{BMP:sec:refine}

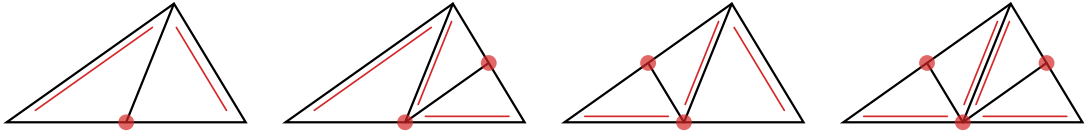
\begin{figure}
    \centering
        \resizebox{0.9\textwidth}{!}{\begin{tikzpicture}[line cap=round]

\definecolor{pyRed}{HTML}{d62728}

\tikzstyle{edge}=[thick]
\tikzstyle{marked}=[circle,fill=pyRed,inner sep=2pt,opacity=0.7]

\newcommand{\referenceEdge}[3]{%
	\coordinate (M) at ($0.333*#1+0.333*#2+0.333*#3$);%
	\coordinate (Am) at ($#1!0.3!(M)$);%
	\coordinate (Bm) at ($#2!0.3!(M)$);%
	\draw[semithick,pyRed] (Am) -- (Bm);}

\newcommand{\drawParentTriangle}{%
	\coordinate (P1) at (0,0);
	\coordinate (P2) at (3,0);
	\coordinate (P3) at (2.1,1.5);
	\draw[edge] (P1) -- (P2) -- (P3) -- cycle;
}

\pgfmathsetmacro{\xdist}{3.5}

\begin{scope}
	\drawParentTriangle
	\coordinate (N1) at ($0.5*(P1)+0.5*(P2)$);
	\draw[edge] (P3) -- (N1);
	\referenceEdge{(P3)}{(P1)}{(N1)}
	\referenceEdge{(P2)}{(P3)}{(N1)}
	\node[marked] at (N1) {};
\end{scope}

\begin{scope}[shift={(1*\xdist,0)}]
	\drawParentTriangle
	\coordinate (N1) at ($0.5*(P1)+0.5*(P2)$);
	\coordinate (N2) at ($0.5*(P2)+0.5*(P3)$);
	\draw[edge] (P3) -- (N1);
	\draw[edge] (N1) -- (N2);
	\referenceEdge{(P3)}{(P1)}{(N1)}
	\referenceEdge{(N1)}{(P2)}{(N2)}
	\referenceEdge{(P3)}{(N1)}{(N2)}
	\node[marked] at (N1) {};
	\node[marked] at (N2) {};
\end{scope}

\begin{scope}[shift={(2*\xdist,0)}]
	\drawParentTriangle
	\coordinate (N1) at ($0.5*(P1)+0.5*(P2)$);
	\coordinate (N3) at ($0.5*(P1)+0.5*(P3)$);
	\draw[edge] (P3) -- (N1);
	\draw[edge] (N1) -- (N3);
	\referenceEdge{(P2)}{(P3)}{(N1)}
	\referenceEdge{(P1)}{(N1)}{(N3)}
	\referenceEdge{(N1)}{(P3)}{(N3)}
	\node[marked] at (N1) {};
	\node[marked] at (N3) {};
\end{scope}

\begin{scope}[shift={(3*\xdist,0)}]
	\drawParentTriangle
	\coordinate (N1) at ($0.5*(P1)+0.5*(P2)$);
	\coordinate (N2) at ($0.5*(P2)+0.5*(P3)$);
	\coordinate (N3) at ($0.5*(P1)+0.5*(P3)$);
	\draw[edge] (P3) -- (N1);
	\draw[edge] (N3) -- (N1);
	\draw[edge] (N1) -- (N2);
	\referenceEdge{(P1)}{(N1)}{(N3)}
	\referenceEdge{(N1)}{(P2)}{(N2)}
	\referenceEdge{(P3)}{(N1)}{(N3)}
	\referenceEdge{(N1)}{(P3)}{(N2)}
	\node[marked] at (N1) {};
	\node[marked] at (N2) {};
	\node[marked] at (N3) {};
\end{scope}

\end{tikzpicture}}
    \caption{%
        Visualization in 2D of the finite number of patterns 
        of new triangles obtained by NVB. 
        Each triangle has an associated refinement edge to be bisected in 
        case the triangle is marked for refinement. NVB marks edges, here 
        indicated by red dots, in order to avoid hanging nodes. The red
        line indicates the new refinement 
        edge -- opposite to the newest vertex --
        on the new triangles.
    }
    \label{BMP:fig:NVB:visualization}
\end{figure}

The convergence analysis of adaptive FEM crucially depends
on a mesh-refine\-ment strategy which avoids over-refinement 
and the degeneration of the shapes of simplices.
From an analytical perspective,
the standard choice is newest-vertex bisection (NVB). 
This strategy is denoted by \(\refine(\cdot)\) 
throughout this chapter. A visualization in the two-dimensional case 
is given in Figure~\ref{BMP:fig:NVB:visualization}.
Let $\TT_h = \refine(\TT_H, \MM_H)$ denote the coarsest conforming
triangulation of \(\Omega\)
obtained from $\TT_H$ by bisecting at least all marked elements of
$\MM_H \subseteq \TT_H$.
A suitable closure procedure ensures that no hanging
nodes (or hanging edges in 3D) exist in the refined mesh
\(\TT_h\).
This requirement is subject to an appropriate enumeration
and connectivity of the vertices of the initial mesh $\TT_0$;
see~\cite{stevenson2008} for $d \ge 2$.
The recent work~\cite{dgs2023} introduced an initialization
strategy for arbitrary initial meshes $\TT_0$ for all $d \ge 2$;
see also the seminal work~\cite{Mitchell91} for $d \ge 2$ for admissible 
$\TT_0$ as well as~\cite{kpp2013} 
for the case $d=2$ and any initial mesh.
The existence of hanging nodes cannot be avoided
for hexahedral meshes, imposing additional
difficulties in the analysis; see, e.g.,~\cite[Sect.~6.3]{bn2010} for details.

Without a second argument, $\TT_h \in \refine(\TT_H)$
represents all triangulations $\TT_h$ obtainable
by finitely many steps (including zero) of NVB refinements
(with arbitrary marked elements in each step).
We define $\T \coloneqq \refine(\TT_0)$
as the set of all meshes which can be
generated from the initial simplicial mesh $\TT_0$ of
$\Omega$ by use of $\refine(\cdot)$.
NVB-generated meshes are indeed locally nested, i.e., 
for any two simplices \(T \in \TT \in \T\)
and \(T' \in \TT' \in \T\) with \(\vert T \cap T' \vert > 0\),
in the sense of the $d$-dimensional measure,
it holds either \(T = T'\) or \(T \subsetneqq T'\) or
\(T' \subsetneqq T\).
Moreover, the NVB refinement ensures that all meshes $\TT_H\in\T$
are uniformly shape regular. This means that for
\(\rho_T\) being the radius of the largest ball inscribed
in \(T\) and \(\operatorname{diam}(T)\) being the diameter of \(T\), there
exists a shape-regularity constant \(\gamma > 0\) such that
\begin{align}
    \label{BMP:eq:shape-regularity}
    \gamma
    \coloneqq
    \sup_{\TT \in \T}
    \max_{T \in \TT}
    \frac{\operatorname{diam}(T)}{\rho_T}
    < \infty.
\end{align}

Every $\TT_H\in\T$ corresponds to an associated
finite element space $\XX_H \subsetneqq \XX$ defined
via~\eqref{BMP:eq:fem-space}.
Note that, for $\TT_H,\TT_h\in\T$ with $\TT_h\in\refine(\TT_H)$,
there holds nestedness of the spaces $\XX_H\subseteq\XX_h$.

Throughout this chapter, the tilde notation \( A \lesssim B \)
abbreviates \( A \leq C \, B\) for a constant \(C > 0\) that
is independent of the mesh size, 
but might depend on \(\TT_0\) (and, thus, the shape-regularity
parameter \(\gamma\)),
the input data \(\boldsymbol{A}\), \(f\), and \(\boldsymbol{f}\),
as well as on the polynomial degree \(p\) of the discretization.
Moreover, \(A \eqsim B\) represents the equivalence
\(A \lesssim B\) and \(B \lesssim A\).

\begin{figure}
    \centering
        \includegraphics[width=0.22\textwidth]{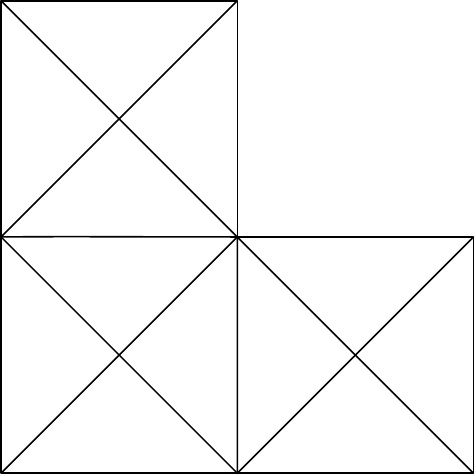} \hfil
        \includegraphics[width=0.22\textwidth]{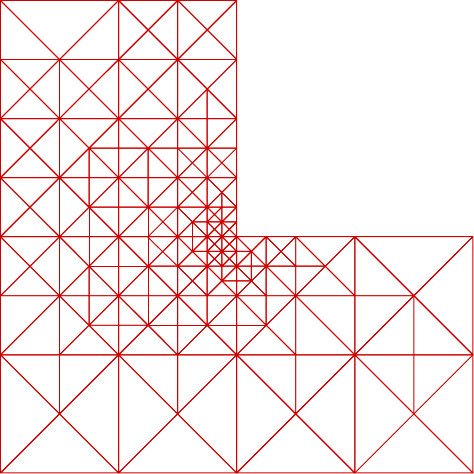} \hfil
        \includegraphics[width=0.22\textwidth]{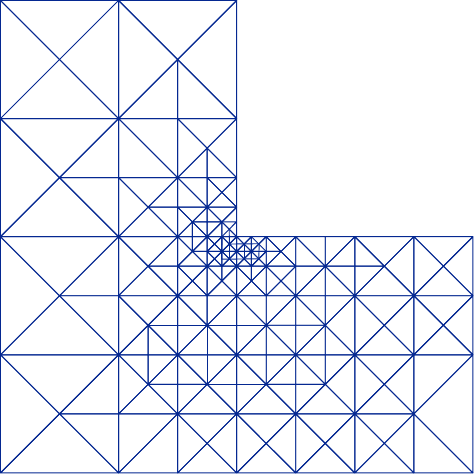} \hfil
        \includegraphics[width=0.22\textwidth]{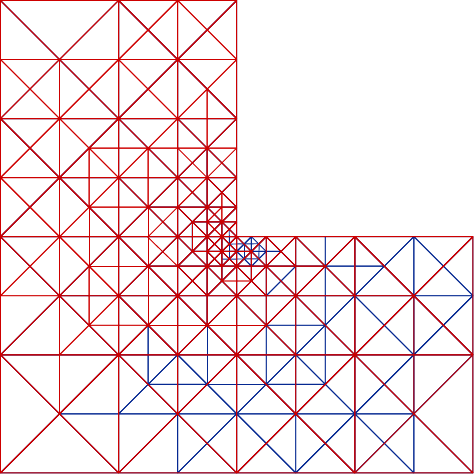}
    \caption{%
        Example in 2D of the overlay of 
        two distinct triangulations \(\TT_h \) and \(\TT_{h'} \), 
        represented in red and 
        blue. The meshes stem from the refinement of the same coarse mesh 
        \(\TT_H \),  
        represented in black, and their overlay 
        \( \TT_h \oplus \TT_{h'} \), i.e., the coarsest common 
        refinement, is illustrated in the rightmost figure.
    }
    \label{BMP:fig:overlay:visualization}
\end{figure}

\begin{figure}
    \centering
    \includegraphics[width=0.22\textwidth]{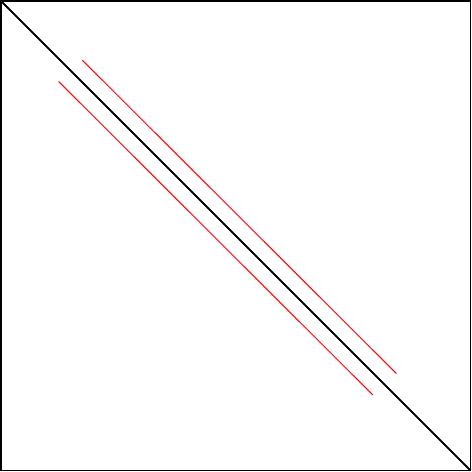} \hfil
    \includegraphics[width=0.22\textwidth]{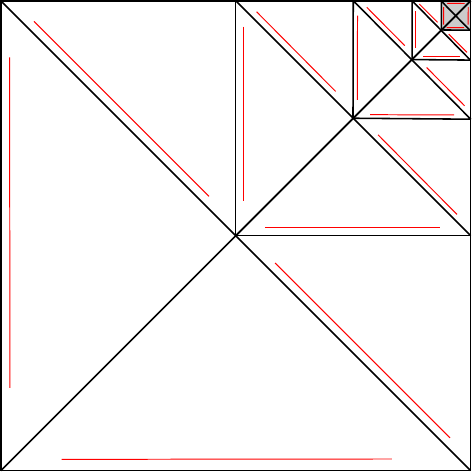} \hfil
    \includegraphics[width=0.22\textwidth]{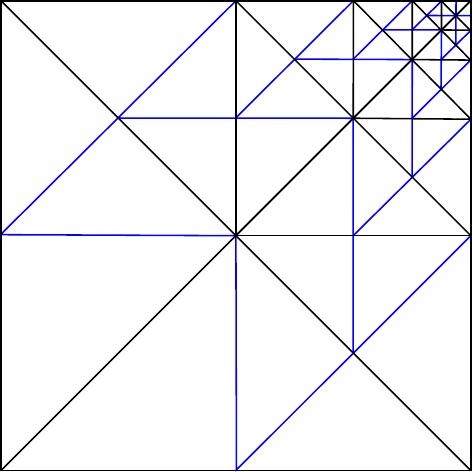}
    \caption{%
        Example in 2D of non-locality of the mesh closure. 
        From the initial triangulation \( \TT_0\) in the left, 
        iterative NVB refinement leads to the middle triangulation 
        \(\TT_\ell\). Only four elements 
        (highlighted in gray, each refinement edge represented by a 
        red line) of \(\TT_\ell\) are marked 
        for refinement. Nonetheless, each triangle of \(\TT_\ell\) 
        is bisected at least 
        once for the new triangulation \(\TT_{\ell+1}\) 
        in the right to be conforming 
        (and avoid hanging nodes) 
        and thus the refinement propagates through the entire mesh.
    }
    \label{BMP:fig:mesh:closure:visualization}
\end{figure}

Our analysis relies on the following 
properties of the mesh refinement:
\begin{enumerate}[label={\textbf{(R\arabic*)}}]
    \item
        \labeltext{R1}{BMP:axiom:children}
        \textbf{splitting property:} Each refined element is
        split into finitely many children, i.e.,
        for all $\TT_H \in \T$ and all $\MM_H \subseteq \TT_H$,
        the triangulation $\TT_h = \refine(\TT_H, \MM_H)$ satisfies that
        \begin{align*}
            \# (\TT_H \setminus \TT_h) + \# \TT_H \leq \# \TT_h
            \lesssim
            \# (\TT_H \setminus \TT_h) + \# (\TT_H \cap \TT_h);
        \end{align*}
    \item
        \labeltext{R2}{BMP:axiom:overlay}
        \textbf{overlay estimate:} For all meshes $\TT_{H} \in
        \T$ and $\TT_h, \TT_{h'} \in \refine(\TT_{H})$,
        there exists a common refinement $\TT_h \oplus
        \TT_{h'} \in \refine(\TT_h) \cap
        \refine(\TT_{h'})
        \subseteq \refine(\TT_{H})$ such
        that
        \begin{align*}
            \# (\TT_h \oplus \TT_{h'}) \leq \# \TT_h + \# \TT_{h'} - \# \TT_{H};
        \end{align*}
    \item
        \labeltext{R3}{BMP:axiom:closure}
        \textbf{mesh-closure estimate:}\;
        For each sequence $(\TT_\ell)_{\ell \in \N_0}$ of
        successively refined
        meshes, i.e.,
        $\TT_{\ell+1} \coloneqq \refine(\TT_\ell,\MM_\ell)$
        with $\MM_\ell \subseteq \TT_\ell$
        for all $\ell \in \N_0$,
        it holds  that
        \begin{align*}
            \# \TT_\ell - \# \TT_0 \
            \lesssim
            \sum_{j=0}^{\ell -1} \# \MM_j.
        \end{align*}
\end{enumerate}
We refer to~\cite{bdd2004,stevenson2007,stevenson2008,ckns2008,kpp2013,gss2014}
for the validity of~\eqref{BMP:axiom:children}--\eqref{BMP:axiom:closure}
for NVB-refinement algorithms.
In particular, it follows from the splitting 
property~\eqref{BMP:axiom:children} that the computation of 
\(\TT_h = \refine(\TT_H,\MM_H)\) can be realized at linear cost
\( \mathcal{O} (\# \TT_H ) \), i.e., the number of operations is bounded
by a fixed multiple of the number of elements in \(\TT_H\).
Figure~\ref{BMP:fig:overlay:visualization} illustrates the overlay
\(\TT_h \oplus \TT_{h'}\) of two 
meshes. The estimate~\eqref{BMP:axiom:overlay} can be seen by noting that 
each coarse simplex \(T \in \TT_H\) has children in \(\TT_h\), and 
\(\TT_{h'}\), respectively, which can be represented by two binary trees 
thanks to the structure of NVB refinement. 
Moreover, note that mesh closure is in general 
\emph{non-local}, i.e., the number of new simplices arising 
in one refinement step cannot be bound by the number of marked elements alone; 
see Figure~\ref{BMP:fig:mesh:closure:visualization} for an example. 
Nonetheless, the estimate~\eqref{BMP:axiom:closure} ensures control on the 
number of generated elements over a mesh hierarchy by the sum of marked 
elements.

\subsection{Module \texttt{SOLVE}}
\label{BMP:sec:solve}

Recall that even when the linear system arising from the
discretization~\eqref{BMP:eq:discrete}
is symmetric and positive definite, its exact solution
is in practice prohibitively expensive. Thus, an iterative solver will
be used to compute an approximation.
One step of the iterative solver for the numerical solution
of the discrete problem~\eqref{BMP:eq:discrete} is formally
denoted as an operator
\(\Psi_H: \XX^* \times \XX_H \to \XX_H\), i.e.,
\(u_H^k \coloneqq \Psi_H(F; u_H^{k-1})\)
is the new approximation constructed by one solver step
from the available approximation $u_H^{k-1}$ and
the right-hand side \(F\) in the so-called dual space
\(\XX^* \equiv H^{-1}(\Omega)\) consisting of all linear and continuous
forms \(F:\XX \to \R \).
A crucial assumption for the later convergence analysis is
the strict contraction of the algebraic solver in the energy norm,
defined by
\[
    \vvvert u \vvvert
    \coloneqq
    \langle A \nabla u, \nabla u \rangle_{L^2(\Omega)}^{1/2}
    \coloneqq
    \bigg( \int_\Omega \boldsymbol{A} \nabla u \cdot \nabla u \, \d{x} \bigg)^{1/2}
    \quad \text{for all } u \in \XX.
\]
This means that there exists
some uniform constant \(0 < \qalg < 1\) at each iteration,
i.e.,
\begin{align}
    \label{BMP:eq:algebra}
    \tag{C}
    \vvvert u_H^\star - u_H^k \vvvert
    \le
    \qalg \,
    \vvvert u_H^\star - u_H^{k-1} \vvvert.
\end{align}
Examples for such solvers include
conjugate gradient methods with optimal additive Schwarz
or BPX-type preconditioners \cite{cnx2012} or
optimal geometric multigrid methods
which respect the multilevel structure of locally graded meshes
of the later adaptive algorithm. For the latter, we refer to \cite{wz2017} 
for discretizations with \(p = 1\) as well as to~\cite{imps2024} 
for higher-order robust solvers and also~\cite{SMPZ08}.
Note that, here, the exact discrete solution $u_H^{\star}$
is never computed but only approximated by $u_H^{k}$.
However, the contraction \eqref{BMP:eq:algebra} and the triangle inequality
imply that
\begin{align}
    \label{BMP:eq:algebra:apost}
    \frac{1 - \qalg}{\qalg} \,\vvvert u_H^{\star} - u_H^{k}\vvvert
    \le
    \vvvert  u_H^{k} - u_H^{k-1}\vvvert
    \le
    (1 + \qalg) \, \vvvert u_H^{\star} - u_H^{k-1}\vvvert .
\end{align}
This means that $\vvvert u_H^{k} - u_H^{k-1}\vvvert $
provides a (computable) a-posteriori estimator
for the \emph{algebraic error} $\vvvert u_H^{\star} - u_H^{k}\vvvert $.
Hence, this estimator allows to stop the algebraic solver when the
algebraic error
is dominated by a fixed factor \( \lambdaalg >0 \) 
times the discretization error measured by 
some computable quantity \(\eta_H(u_H^k)\), i.e.,
\begin{align}
\label{BMP:eq:algebra:stopping}
    \vvvert u_H^k - u_H^{k-1} \vvvert
    \le
    \lambdaalg \,
    \eta_H(u_H^k).
\end{align}
This can be interpreted as numerically balancing the error contributions: 
When the
algebra is sufficiently resolved, the solver should be stopped and
a new mesh-refinement step should take place.
The final index of the iterative solver being the minimal \(k \in \N\), 
where~\eqref{BMP:eq:algebra:stopping} is satisfied, is denoted by \(\kk[H]\).
In this notation,
the dependency on the discretization parameter
\(H\) is omitted if it is clear from the context.
The contraction~\eqref{BMP:eq:algebra} 
immediately establishes a bound for
the distance between the
initial iterate \(u_H^0\)
and the final iterate \(u_H^{\kk}\)
\begin{equation}
    \label{BMP:eq:flc:stability}
    \vvvert u_H^{\kk} - u_H^0 \vvvert
    \le
    \vvvert u_H^\star - u_H^{\kk} \vvvert
    +
    \vvvert u_H^\star - u_H^0 \vvvert
    \le
    (1 + \qalg^{\kk[H]})\,
    \vvvert u_H^\star - u_H^0 \vvvert
    \le
    2 \,
    \vvvert u_H^\star - u_H^0 \vvvert.
\end{equation}
Note that one solver step of the optimally preconditioned
conjugate gradient method \cite{cnx2012} or the geometric multigrid 
method~\cite{wz2017,imps2024} can indeed be implemented at linear cost
\(\mathcal{O} (\# \TT_H) \).

\subsection{Module \texttt{ESTIMATE}}
\label{BMP:sec:estimate}

To account for the discretization error, i.e.,
the inaccuracy of the numerical scheme in approximating
the continuous solution \(u^\star\) to~\eqref{BMP:eq:weak},
we use a residual-based a-posteriori estimator.
To formulate the latter, we require additional regularity of the data.
More precisely, we assume that \(\boldsymbol{A}\vert_T\)
and \(\boldsymbol{f}\vert_T\) are Lipschitz continuous
for all \(T \in \TT_H\).
Since \(v_H \in \XX_H\) is a globally continuous,
\(\TT_H\)-piecewise polynomial of degree \(p \in \N\),
its gradient \(\nabla v_H\) is in general
a discontinuous \(\TT_H\)-piecewise polynomial of degree
\(p - 1 \in \N_0\).
In particular, \(\boldsymbol{A} \nabla v_H\) can jump
across the \((d-1)\)-dimensional faces
$T \cap T'$ of neighboring simplices
\(T, T' \in \TT_H\).
In contrast, note that the exact solution \(u^\star \in \XX\)
of~\eqref{BMP:eq:weak} has vanishing jump
of the normal component
\((\boldsymbol{A}\nabla u^\star - \boldsymbol{f}) \cdot n\)
across all interior faces,
where \(n\) is a unit normal vector on $T \cap T'$.
To measure inconsistency of any
finite element function \(v_H \in \XX_H\) in this respect,
let
\[
    \lbracket (\boldsymbol{A} \nabla v_H - \boldsymbol{f}) \cdot n \rbracket\,
    \vert_{T \cap T'}
    \coloneqq
    \big(
        (\boldsymbol{A} \nabla v_H - \boldsymbol{f})\vert_{T}
        -
        (\boldsymbol{A} \nabla v_H - \boldsymbol{f})\vert_{T'}
    \big) \cdot n
\]
denote the normal jump as scalar-valued function
on the interior face \(T \cap T'\)
(with, e.g., \((\boldsymbol{A} \nabla v_H -
\boldsymbol{f})\vert_{T}\) understood as the trace of a
continuous function from \(T\) to $T \cap T'$).
For any discrete function \(v_H \in \XX_H\),
the residual-based quantity
\begin{subequations}
\label{BMP:eq:estimator}
\begin{equation}
    \label{BMP:eq:estimator:a}
    \eta_H(T,v_H)^2
    \coloneqq
    \vert T \vert^{2/d}\,
    \Vert f + \div(\boldsymbol{A} \nabla v_H  
    - \boldsymbol{f}) \Vert_{L^2(T)}^2
    +
    \vert T \vert^{1/d}\,
    \Vert \lbracket (\boldsymbol{A} \nabla v_H - \boldsymbol{f})
    \cdot n \rbracket \Vert_{L^2(\partial T \cap \Omega)}^2
\end{equation}
is computable for all \(T \in \TT_H\) (up to a quadrature error).
Abbreviate any sum of the local contributions by
\begin{equation}
    \label{BMP:eq:estimator:b}
    \eta_H(\MM_H, v_H)
    \coloneqq
    \bigg(\sum_{T \in \MM_H} \eta_H(T, v_H)^2 \bigg)^{1/2}
    \quad \text{for all } \MM_H \subseteq \TT_H
\end{equation}
\end{subequations}
and set $\eta_H(v_H) \coloneqq \eta_H(\TT_H, v_H)$.
The first summand in~\eqref{BMP:eq:estimator:a} measures the volume
residual of the finite element function \(v_H\) with respect to the strong
form of the PDE~\eqref{BMP:eq:strong}.
Note that any \(v_H \in \XX_H\) satisfies the boundary conditions
in~\eqref{BMP:eq:strong} exactly so that no boundary residual is needed here.
As noted above, the second summand in~\eqref{BMP:eq:estimator:a} can be
understood as a consistency error related to 
\( (\boldsymbol{A} \nabla v_H - \boldsymbol{f}) \notin H (\div; \Omega)\).
Thus, \eqref{BMP:eq:estimator} indicates, at least heuristically,
which simplices have a large contribution to
the discretization error \(\vvvert u^\star - u_H^\star \vvvert\), 
allowing to steer the adaptive
mesh refinement.

For further details on the construction and analysis
of residual-based error estimators,
the reader is referred, e.g., to the 
monographs~\cite{AinsworthOden2000,verfurth2013}.
We emphasize that the a-posteriori error
indicator~\eqref{BMP:eq:estimator} can be computed at overall
linear cost \(\mathcal{O} (\# \TT_H ) \).

The framework of the \emph{axioms of adaptivity} from~\cite{cfpp2014}
summarizes the many contributions to the convergence analysis
with rates.
The following abstract conditions, that indeed hold for
the residual a-posteriori error estimator~\eqref{BMP:eq:estimator},
guarantee plain convergence
(without rates) of the resulting adaptive algorithm.
They involve generic constants
$\Cstab, \Crel, \Cdrel > 0$ and $0 < \qred < 1$, where
\( \Crel \) depends only on \(\gamma\)-shape regularity, while
\( \Cstab, \Cdrel \) depend additionally on the polynomial degree
\(p\), and \( \qred = 2^{-1/(2d)}\) for NVB refinement.
For any \(\TT_h \in \texttt{refine}(\TT_H)\), let $u_H^\star$ and 
$u_h^\star$ denote the
corresponding exact discrete solutions
that are never computed in practice,
but are only used in the analysis.
\begin{enumerate}[label={(A\arabic*)}]
\bf
\item stability:\labeltext{A1}{BMP:axiom:stability}\rm\,
\(
    | \eta_h(\UU_H, v_h) - \eta_H(\UU_H,v_H) |
    \le
    \Cstab \vvvert v_h - v_H \vvvert
\)
for all $v_h\in\XX_h$, $v_H \in \XX_H$, and
any subset $\UU_H \subseteq \TT_H \cap \TT_h$ of unrefined
simplices;
\bf
\item reduction:\labeltext{A2}{BMP:axiom:reduction}\rm\,
$\eta_h(\TT_h \backslash \TT_H, v_H) \le \qred \,
\eta_H(\TT_H \backslash \TT_h, v_H)$ for all 
$v_H \in \XX_H\subseteq \XX_h$;
\bf
\item reliability:\labeltext{A3}{BMP:axiom:reliability}\rm\,
$\vvvert u^\star - u_H^\star \vvvert \le \Crel \, \eta_H(u_H^\star)$;
\bf
\item
discrete reliability:\labeltext{A4}{BMP:axiom:discrete_reliability}\rm\,
$ \vvvert u_h^\star - u_H^\star\vvvert \le \Cdrel \, \eta_H(\TT_H
\backslash \TT_h, u_H^\star)$.
\end{enumerate}
Furthermore, the exact discrete solutions satisfy the following
Galerkin orthogonality
\begin{equation}
    \label{BMP:eq:orthogonality}
    a(u^\star - u_H^\star, v_H) = 0
    \quad\text{for all } v_H \in \XX_H,
\end{equation}
leading, in particular, to the following Pythagorean identity
\begin{equation}
    \label{BMP:eq:Pythagoras}
    \vvvert u^\star - u_H^\star \vvvert^2
    +
    \vvvert u_H^\star - v_H \vvvert^2
    =
    \vvvert u^\star - v_H \vvvert^2
    \quad\text{for all } v_H \in \XX_H.
\end{equation}
This immediately results in the
C\'ea-type quasi-best approximation result
\begin{equation}\label{BMP:eq:cea}
	\vvvert  u^\star - u_H^\star \vvvert \,
	\le \,
	\vvvert  u^\star - v_H \vvvert 
	\quad \text{for all } v_H \in \XX_H.
\end{equation}

Note that the residual a-posteriori error estimator in addition to
reliability~\eqref{BMP:axiom:reliability} satisfies a so-called
\emph{efficiency} property
\begin{equation}
    \label{BMP:eq:efficiency}
    \eta_H(u_H^\star)
    \lesssim
    \vvvert u^\star - u_H^\star \vvvert \,
    + \,
    \operatorname{osc}_H (u_H^\star),
\end{equation}
with the latter term representing the \emph{data oscillation}. It is 
defined as 
\begin{align}
    \begin{split}
    \label{BMP:eq:data:oscillation}
    \operatorname{osc}_H (v_H)^2 
    \coloneqq 
        &\sum_{T \in \TT_H}  \Bigl[
            \vert T \vert^{2/d} \Vert (1 - \Pi^{q}_T) (f - \div \boldsymbol{f} 
        + \div \boldsymbol{A} \nabla v_H) \Vert_{L^2(T)}^2 
        \\
        &
        + \sum_{E \subset \partial T \setminus \Gamma}
        \vert T \vert^{1/d} 
        \Vert (1 - \Pi^{q}_E) 
        \lbracket (\boldsymbol{A} \nabla v_H - \boldsymbol{f})
        \cdot n \rbracket \Vert_{L^2(E)}^2
        \Bigr],
    \end{split}
\end{align}
where \(\Pi^{q}_{T}\) is the  
\(L^2(T)\)-orthogonal projection
onto the space of polynomials of degree 
\(q = 2 \,(p-1)\) on the simplex \(T\) and, 
similarly, \(\Pi^{q}_E\) is the 
\(L^2(E)\)-orthogonal projection
onto the polynomials of degree \(q\) 
on the interior face \(E\).
In particular, one can prove the generalized Céa-type estimate    
\begin{equation}\label{BMP:eq:generalized:cea}
	\vvvert  u^\star - u_H^\star \vvvert + \text{osc}_H(u_H^\star)
	\lesssim
    \vvvert  u^\star - v_H \vvvert + \text{osc}_H(v_H)
	\quad \text{for all } v_H \in \XX_H
\end{equation}
which, together with reliability~\eqref{BMP:axiom:reliability} 
and efficiency~\eqref{BMP:eq:efficiency}, results in
\begin{equation}
    \label{BMP:eq:equivalence:err+osc:est}
    \eta_H(u_H^\star)
    \eqsim
    \inf_{v_H \in \XX_H} 
    \left( \vvvert u^\star - v_H \vvvert + \text{osc}_H(v_H) \right).
\end{equation}

However, since the convergence results in Subsection~\ref{BMP:sec:convergence}
focus on the estimator $\eta_H(u_H^\star)$ only,
efficiency is not explicitly required in the analysis
\cite[Sect.~3.2]{cfpp2014}.
For a-posteriori error analysis of data oscillation terms, we refer to, e.g.,
\cite{KreuzerVeeser2021}, see also references therein.

\subsection{Module \texttt{MARK}}
\label{BMP:sec:mark}

The a-posteriori estimator~\eqref{BMP:eq:estimator}
for the discretization error 
is used within a marking strategy to
indicate the error distribution among simplices
and determine which ones should be refined.
Preferably, a small set of simplices whose estimator
contributions exceed a given percentage $0 < \theta \le 1$
of $\eta_H(v_H)$ are marked.
More precisely, the up-to-date numerical analysis employs
the D\"orfler criterion
introduced in \cite{doerfler1996} for reasons that
are made precise below:
Determine a set $\MM_H \subseteq \TT_H$
of minimal cardinality
such that
\begin{equation}
    \label{BMP:eq:marking_criterion}
    \theta\, \eta_H(u_H^{\kk})^2
    \le
    \eta_H(\MM_H, u_H^{\kk})^2.
\end{equation}
Note that small \(\theta \approx 0\) leads to few marked elements (and, hence,
highly adapted meshes), whereas large \(\theta \approx 1\) essentially
enforces uniform mesh refinement.
A straight-forward implementation of the 
strategy~\eqref{BMP:eq:marking_criterion}
would rely on sorting the estimator contributions
which would result in suboptimal
\(\mathcal{O}(\#\TT_H \log\#\TT_H) \) complexity.
However, a modification of the QuickSelect algorithm
in~\cite{pp2020} allows for linear complexity
on average.
Moreover, \cite{stevenson2007} proposes a
realization of~\eqref{BMP:eq:marking_criterion} based on binning, which determines
$\MM_H \subseteq \TT_H$ satisfying~\eqref{BMP:eq:marking_criterion}
and containing, up to a factor 2, the minimal number of elements.

Having marked elements \(\MM_H \subseteq \TT_H\) for refinement, the 
refinement procedure generates a finer triangulation
\(\TT_h\) such that at least all marked simplices
\(\MM_H \subseteq \TT_H \setminus \TT_h\)
with large estimator contributions are refined.
The use of the D\"orfler marking 
criterion~\eqref{BMP:eq:marking_criterion} with arbitrary 
adaptivity parameter
\(0< \theta\le 1\)
guarantees the following estimator reduction property.
\begin{lemma}[perturbed reduction of the estimator for nested approximations]
    Let \(\TT_H\) be a given triangulation with corresponding 
    finite element space \(\XX_H\).
    Given \(u_H^{\kk} \in \XX_H \) and \( 0<\theta <1 \), 
    let \(\MM_H \subseteq \TT_H\) be the elements marked for 
    refinement by~\eqref{BMP:eq:marking_criterion}. 
    Let \(\TT_h = \textup{\texttt{refine}}(\TT_H, \MM_H)\) with
    corresponding finite element space \(\XX_h\)
    and set \(u_h^0 \coloneqq u_H^\kk \in \XX_H \subseteq \XX_h\).
    Then,
    the solver-generated iterates \( u_h^k \in \XX_h\) 
    with \( 0 \le k \le \kk\),
    satisfy
    \begin{equation}
        \label{BMP:eq:flc:reduction_iterates}
        \eta_h(u_h^{\kk})
        \le
        q_{\theta} \,
        \eta_H(u_H^{\kk})
        +
        2 \, \Cstab \,
        \vvvert
            u_h^\star - u_H^{\kk}
        \vvvert,
        \text{ where }
        q_{\theta}
        \coloneqq
        \big[1 - (1 - \qred^2) \theta\big]^{1/2} <1.
    \end{equation}
\end{lemma}

\begin{proof}
    Since all marked simplices are refined,
    the inclusion $\MM_H \subseteq \TT_H \setminus \TT_h$ holds
    and the D\"orfler marking
    criterion~\eqref{BMP:eq:marking_criterion} yields
    \[
        - \eta_H(
            \TT_H \setminus \TT_h;
            u_H^{\kk}
        )^2
        \leq
        - \eta_H(\MM_H; u_h^{\kk})^2
        \stackrel{\eqref{BMP:eq:marking_criterion}}{\leq}
        - \theta \,
        \eta_H(u_H^{\kk})^2.
    \]
    This estimate, together with 
    reduction~\eqref{BMP:axiom:reduction}, leads to
    \begin{align}
        \begin{split}
        \label{BMP:eq:reduction_est}
        \eta_h(u_H^{\kk})^2
        &=
        \eta_h(\TT_h \cap \TT_H; u_H^{\kk})^2
        +
        \eta_h(\TT_h \setminus \TT_H; u_H^{\kk})^2
        \\
        &\stackrel{\mathclap{\eqref{BMP:axiom:reduction}}}{\leq}
        \eta_H(\TT_h \cap \TT_H; u_H^{\kk})^2
        +
        \qred^2 \,
        \eta_H(\TT_H \setminus \TT_h; u_H^{\kk})^2
        \\
        &=
        \eta_H(u_H^{\kk})^2
        -
        (1 - \qred^2) \,
        \eta_H(\TT_H \setminus \TT_h; u_H^{\kk})^2
        \\
        \;
        &\leq
        \bigl[
            1 - (1 - \qred^2) \theta
        \bigr]
        \,
        \eta_H(u_H^{\kk})^2 \,
        =  q_\theta^2 \, \eta_H(u_H^{\kk})^2 .
        \end{split}
    \end{align}
    Combining this with stability~\eqref{BMP:axiom:stability}, the
    estimate \eqref{BMP:eq:flc:stability},
    and nested iteration \(u_h^0 = u_H^{\kk}\),
    results in
    \begin{align*}
        \eta_h(u_h^{\kk})
        &\leq
        \eta_h(u_H^{\kk})
        +
        \Cstab \,
        \vvvert
            u_h^{\kk} - u_H^{\kk}
        \vvvert
        \leq
        q_\theta \,
        \eta_H(u_H^{\kk})
        +
        \Cstab \,
        \vvvert
            u_h^{\kk} - u_H^{\kk}
        \vvvert
        \\
        &\le
        q_\theta \,
        \eta_H(u_H^{\kk})
        +
        2 \, \Cstab \,
        \vvvert
            u_h^\star - u_H^{\kk}
        \vvvert.
        \qedhere
    \end{align*}
\end{proof}

\subsection{Algorithm 1: Adaptive FEM}
\label{BMP:sec:AFEM}
\labeltext{1}{BMP:algorithm:afem}

\noindent
{\bfseries Input:}
Initial mesh \(\TT_0\),
polynomial degree \(p \in \N\),
initial iterate \(u_0^0 \coloneqq u_0^{\kk} \coloneqq 0\),
marking parameter \(0 < \theta \le 1\),
solver parameter \(\lambdaalg > 0\).

\medskip
\noindent
{\bfseries Adaptive loop:}
For all \(\ell = 0, 1, 2, \dots\),
repeat the following steps {\rmfamily(I)--(III)}:
\begin{enumerate}[leftmargin=2.5em]
    \item[\rmfamily(I)]
        {\ttfamily SOLVE \& ESTIMATE.}
        For all \(k = 1, 2, 3, \dots\),
        repeat the steps {\rmfamily(a)--(b)}:
        \begin{itemize}[leftmargin=2em]
            \item[\rmfamily(a)]
                Compute
                \(
                    u_\ell^k
                    \coloneqq
                    \Psi_\ell(F; u_\ell^{k-1})
                \)
                and the corresponding refinement indicators
                \\
                \(\eta_\ell(T; u_\ell^k)\)
                for all \(T \in \TT_\ell\).
            \item[\rmfamily(b)]
                Terminate \(k\)-loop and
                define \(\kk[\ell] \coloneqq k\)
                provided that
                \begin{equation}
                    \label{BMP:eq:j_stopping_criterion}
                    \vvvert u_\ell^k - u_\ell^{k-1} \vvvert
                    \le
                    \lambdaalg \,
                    \eta_\ell(u_\ell^k).
                \end{equation}
        \end{itemize}
    \item[\rmfamily(II)]
        {\ttfamily MARK.}
        Determine a set \(\MM_\ell \subseteq \TT_\ell\)
        that satisfies the Dörfler 
        criterion~\eqref{BMP:eq:marking_criterion}.
        \item[\rmfamily(III)]
            {\ttfamily REFINE.}
            Generate the new mesh
            \(\TT_{\ell+1} \coloneqq \refine(\TT_\ell, \MM_\ell)\)
            by NVB and employ nested iteration
            \(u_{\ell+1}^0 \coloneqq u_\ell^{\kk}\).
\end{enumerate}

\noindent
{\bfseries Output:}
Sequences of successively
refined triangulations \(\TT_\ell\),
discrete approximations \(u_\ell^k\),
and corresponding error estimators \(\eta_\ell(u_\ell^k)\).
\medskip

\begin{remark}
    Let us comment on the termination 
    criterion~\eqref{BMP:eq:j_stopping_criterion} for the iterative solver used 
    in Algorithm~\ref{BMP:algorithm:afem}.
    \begin{enumerate}[label={\normalfont (\roman*)}]
        \item  According to Lemma~\ref{BMP:prop:overall-error-estimator} below, the sum 
        \(\vvvert u_\ell^k - u_\ell^{k-1} \vvvert
        +
        \eta_\ell(u_\ell^k)\)
        provides a computable upper bound for the error 
        \(\vvvert u^\star - u_\ell^k \vvvert\).
        Stopping the iterative solver with the minimal index \(k \in \N\) such that 
        \eqref{BMP:eq:j_stopping_criterion} holds, can thus be understood as 
        numerical equilibration of both summands.
        \item Moreover, as seen in~\eqref{BMP:eq:algebra:apost} above, the 
        term \(\vvvert u_\ell^k - u_\ell^{k-1} \vvvert\) provides control 
        on the algebraic error \(\vvvert u_\ell^\star - u_\ell^{k} \vvvert\), 
        while 
        \(\vvvert u^\star - u_\ell^\star \vvvert 
        \lesssim  \eta_\ell(u_\ell^\star) \lesssim  
        \eta_\ell(u_\ell^k) +\vvvert u_\ell^\star - u_\ell^{k} \vvvert \)
        controls the discretization error. E.g., for 
        \(\lambdaalg = 1/10\) and omitting other constants, the stopping 
        criterion~\eqref{BMP:eq:j_stopping_criterion} may also be interpreted that
        the algebraic error should be below \(10 \% \) 
        of the discretization error.
        \item An additional consequence of~\eqref{BMP:eq:j_stopping_criterion} 
        is that \(\eta_\ell(u_\ell^\kk) \) provides a computable upper bound 
        for the error \(\vvvert u^\star - u_\ell^\kk \vvvert \) of each 
        respective final iterate; 
        see~\eqref{BMP:eq:overall-error-estimator-last-iter} below.
        \item Finally, it can happen that the iterative solver does not 
        terminate. In this case, it follows from 
        Theorem~\ref{BMP:thm:full-linear-convergence} below that 
        \(u^\star = u_\ell^\star\) and 
        \(\eta_\ell(u_\ell^k)  \rightarrow\eta_\ell(u_\ell^\star)=0 \)
        as \(k \rightarrow \infty\).
    \end{enumerate}
\end{remark}

Since the numerical analysis below studies asymptotic behavior, this algorithm
generates an infinite sequence \( (u_\ell^k)_{\ell \in \N_0} \).
However, a practical implementation may terminate when 
the estimated error is below a user-prescribed threshold \(\tau > 0\), 
i.e., 
\begin{align}
    \label{BMP:eq:termination_criterion}
    \vvvert u_\ell^k - u_\ell^{k-1} \vvvert
    +
    \eta_\ell(u_\ell^k)
    \leq \tau,
\end{align}
or after a maximal dimension \(\dim(\XX_\ell)\) or runtime 
is reached.
The following lemma shows that the left-hand side
of~\eqref{BMP:eq:termination_criterion} provides indeed
a-posteriori error control on
\( \vvvert u^\star - u_\ell^{k} \vvvert \)
and that the stopping criterion~\eqref{BMP:eq:j_stopping_criterion}
is designed in a way that it ensures reliability of
\( \eta_\ell(u_\ell^\kk) \) for the a-posteriori
error control of
\( \vvvert u^\star - u_\ell^{\kk} \vvvert \).

\begin{lemma}[a-posteriori control of the overall error]
    \label{BMP:prop:overall-error-estimator}
    For \(\ell \in \N_0\) and \(k \ge 1\), let \(u_\ell^k \in \XX_\ell \)
    be computed by the adaptive loop of Algorithm~\ref{BMP:algorithm:afem}.
    Then, there holds
    \begin{equation}
        \label{BMP:eq:overall-error-estimator}
        \vvvert u^\star - u_\ell^{k} \vvvert
        \lesssim
        \vvvert u_\ell^k - u_\ell^{k-1} \vvvert
        +\eta_\ell(u_\ell^k)
    \end{equation}
    and for \( k = \kk[\ell]\) 
        \begin{equation}
            \label{BMP:eq:overall-error-estimator-last-iter}
            \vvvert u^\star - u_\ell^{\kk} \vvvert
            \lesssim
            \eta_\ell(u_\ell^\kk).
        \end{equation}
\end{lemma}

\begin{proof}
    The reliability~\eqref{BMP:axiom:reliability}, the
    stability~\eqref{BMP:axiom:stability}, and a-posteriori control
    for the algebraic error~\eqref{BMP:eq:algebra:apost} yield the result
    \begin{align*}
        \vvvert u^\star - u_\ell^{k} \vvvert
        &\le
        \vvvert u^\star - u_\ell^\star \vvvert
            +  \vvvert u_\ell^\star - u_\ell^{k} \vvvert
        \stackrel{\eqref{BMP:axiom:reliability}}{\le}
        \Crel \, \eta_\ell(u_\ell^\star)
            +  \vvvert u_\ell^\star - u_\ell^{k} \vvvert \\
        &
        \stackrel{\mathclap{\eqref{BMP:axiom:stability}}}{\le}
            \Crel \, \eta_\ell(u_\ell^k)
        + (\Crel \Cstab +1)\, \vvvert u_\ell^\star - u_\ell^{k} \vvvert.
        \\
        &
        \stackrel{\mathclap{\eqref{BMP:eq:algebra:apost}}}{\le}
        \Crel \, \eta_\ell(u_\ell^k)
        + (\Crel \Cstab +1)\,
        \frac{\qalg}{1 - \qalg} \,
        \vvvert u_\ell^k - u_\ell^{k-1} \vvvert.
    \end{align*}
    In particular, thanks to the stopping
    criterion~\eqref{BMP:eq:j_stopping_criterion}, there holds
    \[
        \vvvert u^\star - u_\ell^{\kk} \vvvert
        \le
        \big[
        \Crel \, + (\Crel \Cstab +1)\,
        \frac{\qalg}{1 - \qalg} \, \lambdaalg \big]
        \,
        \eta_\ell(u_\ell^\kk).
        \qedhere
    \]
\end{proof}

\subsection{Notions of quasi-errors and computational cost}
\label{BMP:sec:quasi-errors}

The sequential nature of Algorithm~\ref{BMP:algorithm:afem}
motivates the countably infinite index set
\[
    \QQ
    \coloneqq
    \{
        (\ell, k) \in \N_0^2
        \;:\;
        u_\ell^k \in \XX
        \text{ is defined in Algorithm~\ref{BMP:algorithm:afem}}
    \}.
\]
The set \(\QQ\) is equipped with
the linear order
\[
    (\ell', k') \leq (\ell, k)
    :\Longleftrightarrow
    u_{\ell'}^{k'} \text{ is used no later than }
    u_\ell^k
    \text{ in Algorithm~\ref{BMP:algorithm:afem}}
\]
and the total step counter
\[
    \vert \ell, k \vert
    \coloneqq
    \#\{
        (\ell', k') \in \QQ
        \;:\;
        (\ell', k') \leq (\ell, k)
    \}
    \in
    \N_0
    \quad
    \text{for all }
    (\ell, k) \in \QQ.
\]
Consistently with the stopping indices in
Algorithm~\ref{BMP:algorithm:afem}, we define
\begin{align*}
    \elll
    &\coloneqq
    \sup \{
        \ell \in \N_0
        \;:\; (\ell, 0) \in \QQ
    \}
    \in
    \N_0
    \cup
    \{\infty\},
    \\
    \kk[\ell]
    &\coloneqq
    \sup \{
        k \in \N_0
        \;:\; (\ell, k) \in \QQ
    \}
    \in
    \N
    \cup
    \{\infty\},
    \text{ whenever } (\ell, 0) \in \QQ.
\end{align*}
To simplify notation,
the dependency of indices
is omitted whenever the relation is clear
from the context, e.g.,
\(\kk \coloneqq \kk[\ell]\) in
\(u_\ell^{\kk} = u_{\ell}^{\kk[\ell]}\)
or \((\ell, \kk) = (\ell, \kk[\ell])\).

It is important to note that each of the modules 
in Subsections~\ref{BMP:sec:refine}--\ref{BMP:sec:mark} above
can be realized in linear complexity, i.e.,
for a given mesh \(\TT_\ell\),
the number of arithmetic and logical operators
is proportional to the number of elements \(\#\TT_\ell\), which in turn, 
for any fixed polynomial degree, is proportional to the number of degrees of 
freedom, 
i.e., the dimension of the discrete finite element space \(\XX_\ell\).
Due to the nested nature of the adaptive algorithm, the total
computational cost (and, hence, empirical computational time)
to compute \(u_\ell^k\) is thus proportional to
\begin{equation}
    \label{BMP:eq:cost}
    \cost(\ell, k)
    \coloneqq
    \sum_{\substack{
        (\ell^\prime, k^\prime) \in \QQ
        \\
        \vert \ell^\prime, k^\prime \vert
        \le \vert \ell, k\vert
    }}
    \# \TT_{\ell'}
    \quad \text{for all } (\ell, k) \in \QQ.
\end{equation}
For the empirical comparison of the performance for various polynomial degrees \(p\),
a cost term is introduced using the dimension of the corresponding
discrete spaces
\begin{equation}
    \label{BMP:eq:cost_p}
    \cost_p(\ell, k)
    \coloneqq
    \sum_{\substack{
        (\ell^\prime, k^\prime) \in \QQ
        \\
        \vert \ell^\prime, k^\prime \vert
        \le \vert \ell, k\vert
    }}
    \dim(\XX_{\ell'})
    \quad \text{for all } (\ell, k) \in \QQ.
\end{equation}
Note that the two cost terms~\eqref{BMP:eq:cost} and~\eqref{BMP:eq:cost_p} are
equivalent with equivalence constants depending on the polynomial degree.

The following quasi-error is introduced to measure 
the algebraic and the discretization error
\begin{equation}
    \label{BMP:eq:quasi-error}
    \Eta_\ell^k
    \coloneqq
    \vvvert u_\ell^\star - u_\ell^k \vvvert
    +
    \eta_\ell(u_\ell^k)
    \quad \text{for all } (\ell, k) \in \QQ.
\end{equation}
Arguing as in the proof of Lemma~\ref{BMP:prop:overall-error-estimator}
above, the following
equivalences can be shown for the final iterates
\[
    \eta_\ell(u_\ell^\kk)
    \eqsim
    \Eta_\ell^\kk
    = \vvvert u_\ell^\star - u_\ell^\kk \vvvert
    +
    \eta_\ell(u_\ell^\kk)
    \eqsim
    \vvvert u^\star - u_\ell^\kk \vvvert
    +
    \eta_\ell(u_\ell^\kk),
\]
linking the present quasi-error~\eqref{BMP:eq:quasi-error} with
that of~\cite{ckns2008}.
The contraction of the algebraic error ensures
the following perturbed contraction of the weighted quasi-error.
\begin{lemma}[contraction of quasi-error of final iterates]
    \label{BMP:lem:contraction_quasi-error}
    Let \(\gamma \in (0,1)\) such that
    \(
        \gamma
        <
        (1 - \qalg)/(2 \Cstab)
    \).
    Consider the weighted quasi-error of the final iterates
    \begin{equation}
        \label{BMP:eq:quasi-error-final-iterates}
        \Eta_\ell
        \coloneqq
        \vvvert
            u_\ell^\star - u_\ell^{\kk}
        \vvvert
        +
        \gamma \,
        \eta_\ell(u_\ell^{\kk})
        \quad\text{for all }
        (\ell,\kk) \in \QQ.
    \end{equation}
    Then, there exists a constant
    \(0 < \qctr < 1\) such that
    \begin{equation}
        \label{BMP:eq:contr-quasierr}
        \Eta_{\ell+1}
        \le
        \qctr \, \Eta_\ell
        +
        \qctr \,
        \vvvert u_{\ell+1}^\star - u_\ell^\star \vvvert
        \quad \text{for all } (\ell+1, \kk) \in \QQ.
    \end{equation}
\end{lemma}

\begin{proof}
    The solver contraction~\eqref{BMP:eq:algebra}, nested iteration
    $u_{\ell+1}^{0} = u_{\ell}^{\kk}$, and
    the estimator reduction~\eqref{BMP:eq:flc:reduction_iterates}
    prove, for $\qctr \coloneqq \max\{ \qalg + 2 \Cstab
    \gamma, q_\theta \} \in (0,1)$,
    \begin{align*}
            \vvvert
                u_{\ell+1}^\star - u_{\ell+1}^{\kk}
            \vvvert
            +
            \gamma \,
            \eta_{\ell+1}(u_{\ell+1}^{\kk})
            &
            \stackrel{\mathcal{\eqref{BMP:eq:algebra}}}{\le}
            \qalg^{\kk [\ell +1]} \,
            \vvvert
                u_{\ell+1}^\star - u_\ell^{\kk}
            \vvvert
            +
            \gamma \,
            \eta_{\ell+1}(u_{\ell+1}^{\kk})
            \\
            &
            \stackrel{\mathcal{\eqref{BMP:eq:flc:reduction_iterates}}}{\leq}
            (\qalg + 2 \Cstab \gamma) \,
            \vvvert
                u_{\ell+1}^\star - u_\ell^{\kk}
            \vvvert
            +
            q_\theta \, \gamma \,
            \eta_\ell(u_\ell^{\kk})
            \\
            &\le
            \qctr \,
            \bigl[
                \vvvert
                    u_{\ell+1}^\star - u_\ell^{\kk}
                \vvvert
                +
                \gamma \,
                \eta_\ell(u_\ell^{\kk})
            \bigr].
    \end{align*}
    The triangle inequality concludes the proof
    of~\eqref{BMP:eq:contr-quasierr}.
\end{proof}

\subsection{Convergence and complexity}
\label{BMP:sec:convergence}

The first main result asserts full R-linear convergence
of the quasi-error \(\Eta_\ell^k\) introduced in~\eqref{BMP:eq:quasi-error}.
More precisely, estimate~\eqref{BMP:eq:full-linear-convergence} below proves that
independently of the algorithmic decision for either local
mesh refinement or yet another solver step, Algorithm~\ref{BMP:algorithm:afem} 
essentially causes contraction of the
quasi-error.
In particular, we see that the D\"orfler marking
criterion~\eqref{BMP:eq:marking_criterion} is sufficient to prove
R-linear convergence of
\( \Eta_\ell^\kk \eqsim \eta_\ell(u_\ell^{\kk}) \).
\begin{theorem}[parameter-robust full R-linear convergence]
    \label{BMP:thm:full-linear-convergence}
    Let \(0 < \theta \le 1\) and \(\lambdaalg > 0\)
    be arbitrary.
    Then, Algorithm~\ref{BMP:algorithm:afem} guarantees
    the existence of constants \(\Clin > 0\)
    and \(0 < \qlin < 1\) such that
    \begin{equation}
        \label{BMP:eq:full-linear-convergence}
        \Eta_\ell^k
        \le
        \Clin \,
        \qlin^{
            \vert \ell, k \vert
            - \vert \ell^\prime, k^\prime \vert
        } \,
        \Eta_{\ell^\prime}^{k^\prime}
        \quad\text{for all }
        (\ell, k), (\ell^\prime, k^\prime) \in \QQ
        \text{ with }
        \vert \ell, k \vert > \vert \ell^\prime, k^\prime \vert.
    \end{equation}
    In particular, this yields parameter-robust convergence
    \begin{align}
    \label{BMP:eq:full-linear-convergence-robust}
        \vvvert u^\star - u_\ell^k \vvvert
        \lesssim \Eta_\ell^k
        \le
        \Clin \,
        \qlin^{
            \vert \ell, k \vert
        } \,
        \Eta_0^0 \rightarrow 0 \quad \text{as }
        \vert \ell, k \vert \rightarrow \infty.  
    \end{align}
\end{theorem}
The proof of Theorem~\ref{BMP:thm:full-linear-convergence} is
given in Section~\ref{BMP:sec:full_linear_convergence} below.
Theorem~\ref{BMP:thm:full-linear-convergence} guarantees convergence
(without rates) for any
choice of the parameters \(\theta\) and \(\lambdaalg\) in
Algorithm~\ref{BMP:algorithm:afem}.

The formal statement of optimal convergence rates employs
the notion of nonlinear approximation classes \cite{bdd2004}.
The computational effort for the solution of the
discrete problem~\eqref{BMP:eq:discrete} on some mesh
\(\TT_{H} \in \T_N \coloneqq \{ \TT \in \T\,:\, \# \TT
- \# \TT_0 \leq N\}\)
compared to the solution on the initial mesh \(\TT_0\)
is bounded (up to a multiplicative constant)
by \(N \in \N_0\).
The value
\begin{align}
\label{BMP:eq:best-error}
    e(u^\star, \TT_0, N)
    \coloneqq
    \min_{
        \TT_{\textup{opt}}
        \in
        \T_N
    }
    \eta_{\textup{opt}}(u^\star_{\textup{opt}})
\end{align}
represents the best possible error measured by the estimator
\(\eta_{\textup{opt}}\)
on the optimal (but not available) mesh \(\TT_{\textup{opt}}\)
with at most \(N\) additional simplices compared to \(\TT_0\).
Given \(s > 0\),
the approximation class \(\A_s\)
consists of all \(u^\star\) such that
\begin{align}
\label{BMP:eq:approx-class}
    \Vert u^\star \Vert_{\A_s}
    \coloneqq
    \sup_{N \in \N_0}
    (N + 1)^s
    e(u^\star, \TT_0, N)
    < \infty.
\end{align}
When the above expression is finite,
the regularity parameter \(s > 0\) 
represents an achievable convergence rate:
The best possible error \(e(u^\star, \TT_0, N)\)
decreases at least at rate $-s$ with respect to the number of additional 
elements \(N\) and thus
the multiplicative terms in~\eqref{BMP:eq:approx-class} counterbalance.
A sequence of meshes \(\TT_\ell \in \T\) is called
\emph{quasi-optimal}, if there holds
\[
    \sup_{\ell \in \N_0}
    (\# \TT_\ell)^s
    \eta_\ell(u_\ell^\star)
    \eqsim
    \Vert u^\star \Vert_{\A_s}
    \quad\text{for all } s > 0,
\]
i.e., \((\eta_\ell(u_\ell^\star))_{\ell \in \N_0}\)
converges with rate \(-s\) in the sense of
\[
    \eta_\ell(u_\ell^\star)
    \lesssim
    \Vert u^\star \Vert_{\A_s}
    (\# \TT_\ell)^{-s}
    \quad \text{for all } \ell \in \N_0
    \ \text{and all } s>0.
\]
Since this estimate holds for all \( s > 0\) and since 
\(\Vert u^\star \Vert_{\A_s}\)
is finite if and only if the rate \(-s\) is possible, this
constitutes the optimal convergence rate
with respect to the number of elements and, hence, for fixed polynomial
degree, with respect to the number of degrees of freedom.
However, the rate with respect to the overall 
computational cost,
equivalent to 
\(\cost(\ell, k)\) from~\eqref{BMP:eq:cost},
is of more practical interest.
In addition to the discretization error measured via
\(\eta_\ell(u_\ell^\star)\), one has to include the algebraic error
\(\vvvert u^\star_{\ell} - u_\ell^k \vvvert\)
as contribution to the quasi-error. 
Stability~\eqref{BMP:axiom:stability} ensures
\(
    \vvvert u_\ell^\star - u_\ell^k \vvvert
    +
    \eta_\ell(u_\ell^\star)
    \eqsim
    \vvvert u_\ell^\star - u_\ell^k \vvvert
    +
    \eta_\ell(u_\ell^k)
    = \Eta_\ell^k
\), advocating for the quasi-error~\eqref{BMP:eq:quasi-error}.
This leads to the following notion of optimal complexity 
for sufficiently small adaptivity parameters \(\theta\) and 
\( \lambdaalg \).
\begin{theorem}[optimal complexity]
    \label{BMP:thm:optimal_complexity}
    There exist upper bounds \(0 < \theta^\star \le 1\)
    and \(\lambdaalg^\star > 0\) such that,
    for any \(0 < \theta < \theta^\star\) and
    \(0 < \lambdaalg < \lambdaalg^\star\),
    the following holds:
    Algorithm~\ref{BMP:algorithm:afem} guarantees that
    \[
        \Vert u^\star \Vert_{\A_s}
        \lesssim
        \sup_{\substack{(\ell,k) \in \QQ}}
        \cost(\ell, k)^s\,
        \Eta_\ell^k
        \lesssim
        \max\{
            \Vert u^\star \Vert_{\A_s},
            \Eta_0^0
        \}
        \quad\text{for all } s > 0,
    \]
    where the underlying constants depend on \(s\). 
    In particular, \( \Vert u^\star \Vert_{\A_s} < \infty \) 
    is equivalent to the quasi-error \(\Eta_\ell^k\) decaying with 
    rate \(-s\) with respect to the overall computational cost, 
    and, hence, time.
\end{theorem}
The proof of Theorem~\ref{BMP:thm:optimal_complexity} is given in Section~\ref{BMP:sec:optimal_complexity} below.

%
%

\section{Numerical experiments}
\label{BMP:sec:numerics}

This section provides numerical investigations
underlining the theoretical results from
Subsection~\ref{BMP:sec:convergence}.
The experiments employ the \textsc{Matlab} object-oriented code\footnote{The 
\textsc{Matlab} code for the reproduction of all
the experiments is openly available under
\href{https://www.tuwien.at/mg/asc/praetorius/software/mooafem}%
{https://www.tuwien.at/mg/asc/praetorius/software/mooafem}} 
from~\cite{MooAFEM} and \(hp\)-robust geometric multigrid method from~\cite{imps2024} 
for the arising linear systems.%

Following~\cite{kellogg1974}, we consider the benchmark 
problem~\eqref{BMP:eq:strong}
on the square domain \(\Omega \coloneqq (-1, 1)^2\)
with a strong jump in 
the piecewise constant scalar diffusion coefficient
\(\boldsymbol{A} (x) \coloneqq 
a(x)  \, I_{2 \times 2} \in L^\infty(\Omega)\)
for 
\(a(x) \coloneqq 161.447\,638\,797\,588\,1\) if \(x_1 x_2 > 0\)
and \(a(x) \coloneqq 1\) if \(x_1 x_2 < 0\).
This setting models an interface problem.
The exact weak solution in polar coordinates reads
\(u^\star(r, \phi) \coloneqq r^\alpha \mu(\phi)\) 
with constants \(\alpha = 0.1\),
\(\beta = -14.922\,565\,104\,551\,52\), \(\delta = \pi/4\),
and \(\mu (\phi)\) defined as 
\[
    \mu(\phi)
    \coloneqq
    \begin{cases}
        \cos((\pi/2 - \beta)\alpha) \,
        \cos((\phi - \pi/2 + \delta) \alpha)
        & \text{if } 0 \leq \phi < \pi / 2,
        \\
        \cos(\delta \alpha)\, \cos((\phi - \pi + \beta) \alpha)
        & \text{if } \pi / 2 \leq \phi < \pi,
        \\
        \cos(\beta \alpha) \, \cos((\phi - \pi - \delta) \alpha)
        & \text{if } \pi \leq \phi < 3 \pi / 2,
        \\
        \cos((\pi/2 - \delta)\alpha) \,
        \cos((\phi - 3\pi/2 - \beta) \alpha)
        & \text{if } 3 \pi / 2 \leq \phi < 2 \pi.
    \end{cases}
\]
The solution determines the inhomogeneous Dirichlet boundary conditions
\(u_{\textup{D}}(x) \coloneqq u^\star(x)\) 
for \(x \in \partial\Omega\)
and leads to \(f \equiv 0\)
and \(\boldsymbol{f} \equiv 0\).
The parameter \(\alpha\) gives the regularity
of the solution \(u \in H^{1+\alpha-\varepsilon}(\Omega)\)
for all \(\varepsilon > 0\), having a strong point singularity 
at the origin where the interfaces intersect.
The boundary data
\(u_{\textup{D}}\) is approximated by
nodal interpolation \cite{mns2003,fpp2014}, leading to
an additional boundary-data oscillation term (for sufficiently smooth 
\(u_{\textup{D}}\), see~\cite{fpp2014} in 2D and~\cite{afkpp2013} 
for arbitrary spatial dimension),
in the error estimator
\begin{align*}
    \eta_H(T,v_H)^2
    &\coloneqq
    \vert T \vert\,
    \Vert \div(a \nabla v_H) \Vert_{L^2(T)}^2
    +
    \vert T \vert^{1/2}\,
    \Vert [a \nabla v_H \cdot n] \Vert_{L^2(\partial T \setminus \Gamma)}^2
    \\
    &\hphantom{{}\coloneqq{}}
    +
    \sum_{E \subset \partial T \cap \Gamma}
    \vert T \vert^{1/2}
    \Vert 
    (1 - \Pi_{E}^{p-1}) \, \partial u_{\textup{D}} / \partial s
    \Vert_{L^2(E)},
\end{align*}
where \(\Pi_{E}^{p-1}\) is the \(L^2(E)\)-orthogonal projection
onto the space of polynomials of degree \(p-1\) on the boundary face 
\(E \subset \Gamma\) and 
\(\partial u_{\textup{D}} / \partial s\) is the arc-length derivative 
of \(u_{\textup{D}}\).

\begin{figure}
    \centering  
    \subfloat[Initial mesh \(\TT_0\).]{%
        \label{BMP:fig:Kellogg:initial_mesh}
        \begin{tikzpicture}
    \begin{axis}[%
        axis equal image,%
        width=5.0cm,%
        xmin=-1.15, xmax=1.15,%
        ymin=-1.15, ymax=1.15,%
        font=\footnotesize%
    ]
        \addplot graphics [xmin=-1, xmax=1, ymin=-1, ymax=1]
        {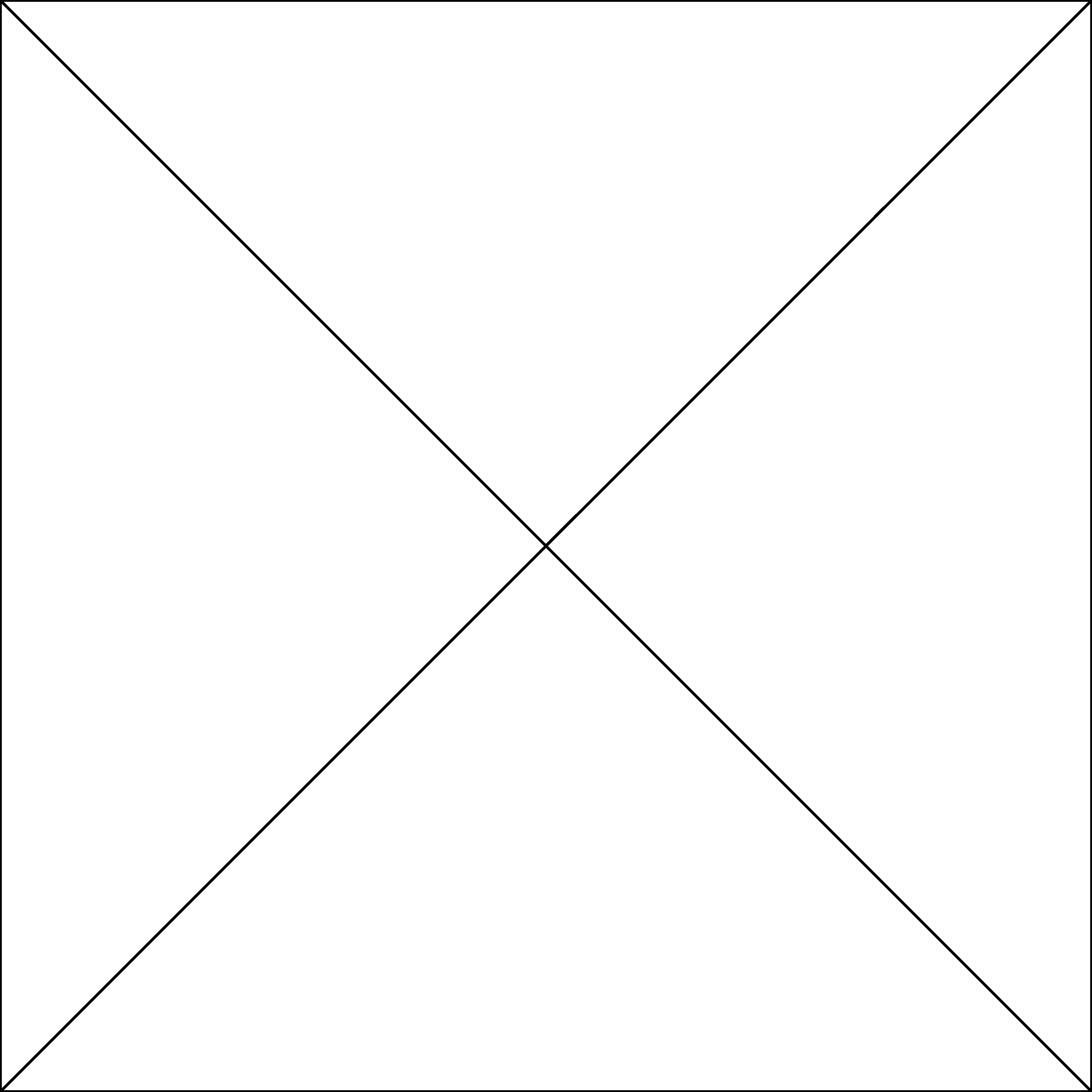};
    \end{axis}
\end{tikzpicture}
    }
    \hfil
    \subfloat[Adaptively refined mesh on level \(\ell = 15\) with \(518\) triangles.]{%
        \label{BMP:fig:Kellogg:mesh}
        \begin{tikzpicture}
    \begin{axis}[%
        axis equal image,%
        width=5.0cm,%
        xmin=-1.15, xmax=1.15,%
        ymin=-1.15, ymax=1.15,%
        font=\footnotesize%
    ]
        \addplot graphics [xmin=-1, xmax=1, ymin=-1, ymax=1]
        {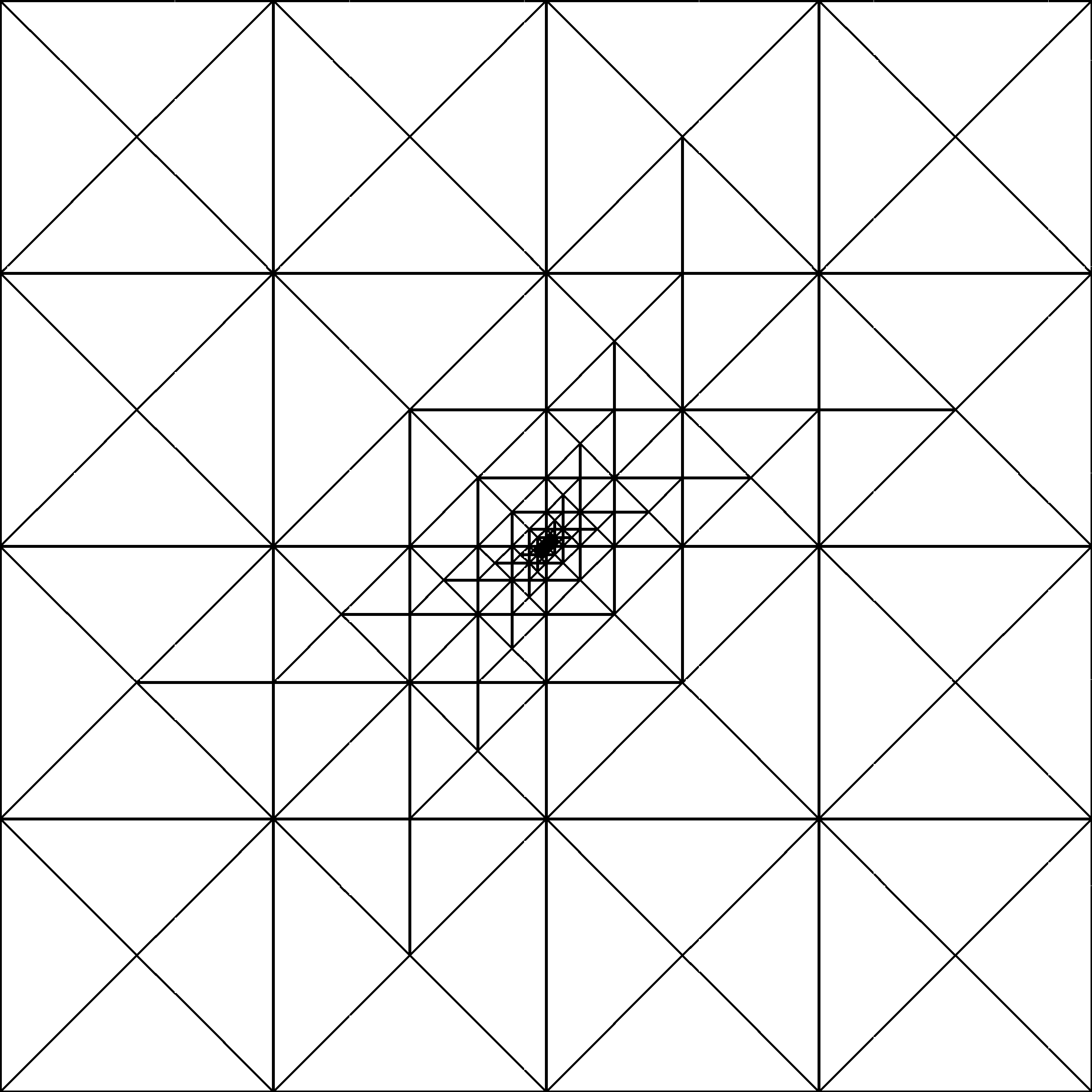};
    \end{axis}
\end{tikzpicture}
    }

    \subfloat[Exact solution \(u^\star\).]{%
        \label{BMP:fig:Kellogg:exact_solution}
        \begin{tikzpicture}
    \pgfplotsset{/pgf/number format/fixed}
    \begin{axis}[%
        width=4.2cm,%
        xmin=-1.1, xmax=1.1,%
        ymin=-1.1, ymax=1.1,%
        zmin=-0.09, zmax=0.09,%
        font=\footnotesize,%
    ]
        \addplot3 graphics [%
            points={%
                (1,-1,0) => (0.1,320.6-195.8)
                (1,0,-0.0782172) => (176.9,320.6-297.1)
                (-1,1,0) => (481.9,320.6-124.7)
                (-1,-1,0.0812259) => (128.9,320.6-0.0)
            }%
            ]
            {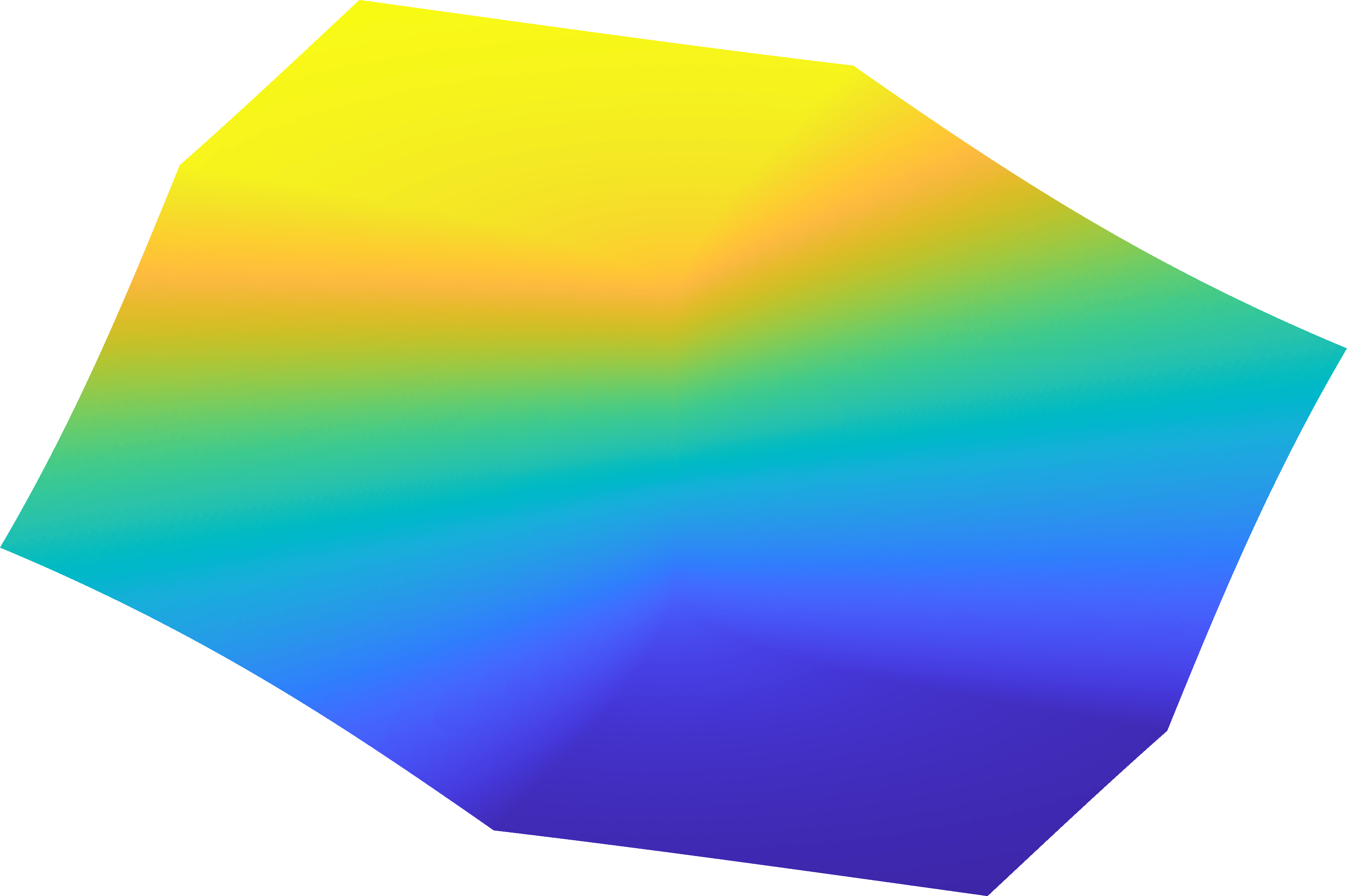};
    \end{axis}
\end{tikzpicture}
    }
    \hfil
    \subfloat[Discrete solution \(u_{15}^{\kk} \in \XX_{15}\).]{%
        \label{BMP:fig:Kellogg:solution}
        \begin{tikzpicture}
    \pgfplotsset{/pgf/number format/fixed}
    \begin{axis}[%
        width=4.2cm,%
        xmin=-1.1, xmax=1.1,%
        ymin=-1.1, ymax=1.1,%
        zmin=-0.09, zmax=0.09,%
        font=\footnotesize,%
    ]
        \addplot3 graphics [%
            points={%
                (-1,-1,0.0812259) => (128.9,321.1-0.1)
                (-1,1,0) => (481.8,321.1-125.2)
                (0,1,-0.0782172) => (417.5,321.1-262.0)
                (1,1,-0.0812259) => (353.2,321.1-320.9)
            }%
            ]
            {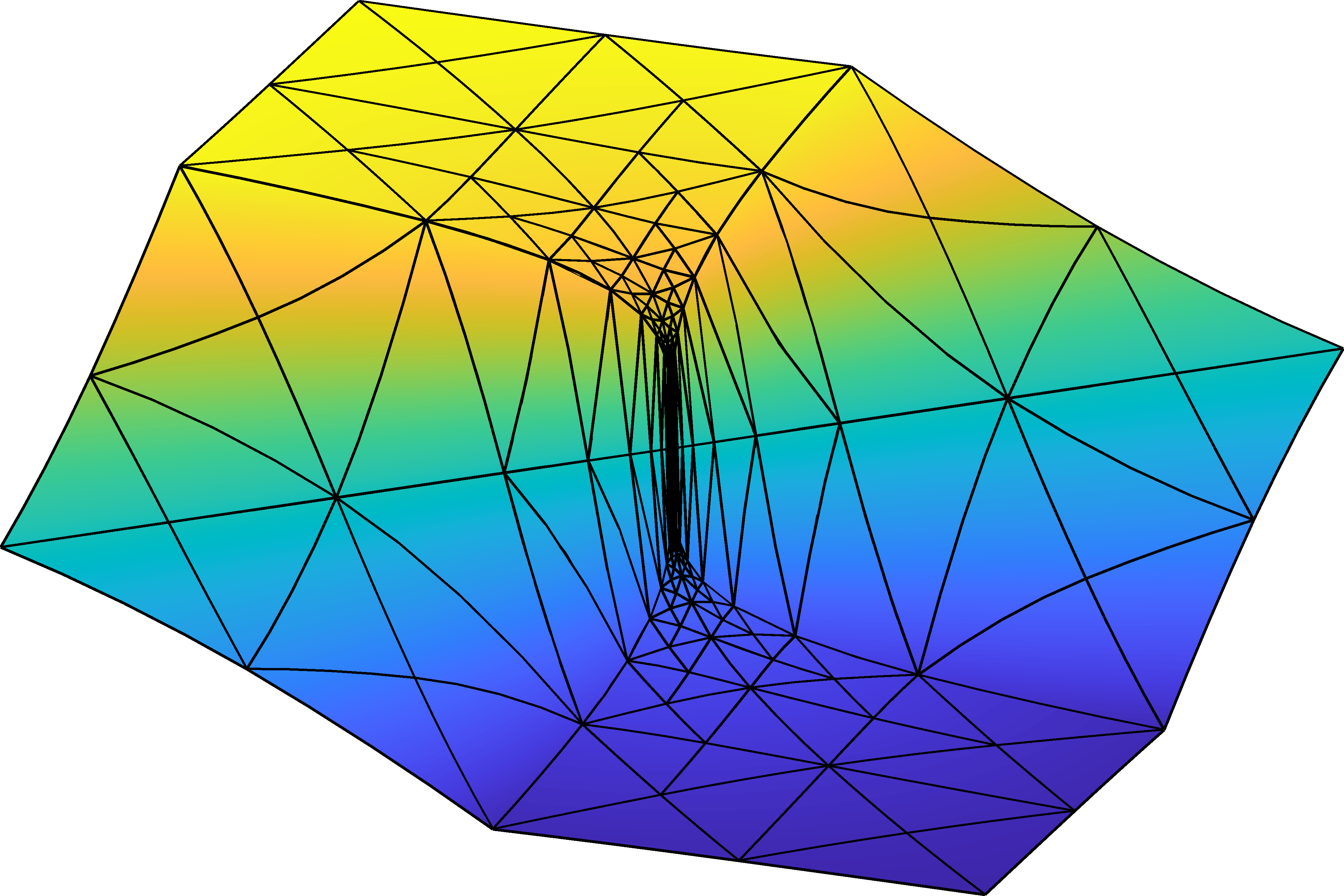};
    \end{axis}
\end{tikzpicture}
    }
    \caption{%
        Initial mesh, adaptively refined mesh, exact solution \(u^\star\),
        and discrete solution \(u_\ell^\kk\)
        for the Kellogg benchmark problem from Subsection~\ref{BMP:sec:numerics}.
        The results are generated by Algorithm~\ref{BMP:algorithm:afem} with
        polynomial degree \(p = 2\), bulk parameter \(\theta = 0.5\),
        and stopping parameter \(\lambdaalg = 0.01\).
    }
    \label{BMP:fig:Kellogg:mesh_solution}
\end{figure}

Since the exact solution \(u^\star\), 
see Figure~\ref{BMP:fig:Kellogg:exact_solution},
exhibits a strong singularity at the origin, 
Algorithm~\ref{BMP:algorithm:afem} aggressively refines 
the initial mesh from Figure~\ref{BMP:fig:Kellogg:initial_mesh} 
towards this point
as displayed in Figure~\ref{BMP:fig:Kellogg:mesh}.
Figure~\ref{BMP:fig:Kellogg:solution} shows the corresponding discrete
solution \(u_\ell^{\kk}\) on the level \(\ell = 15\).

\begin{figure}
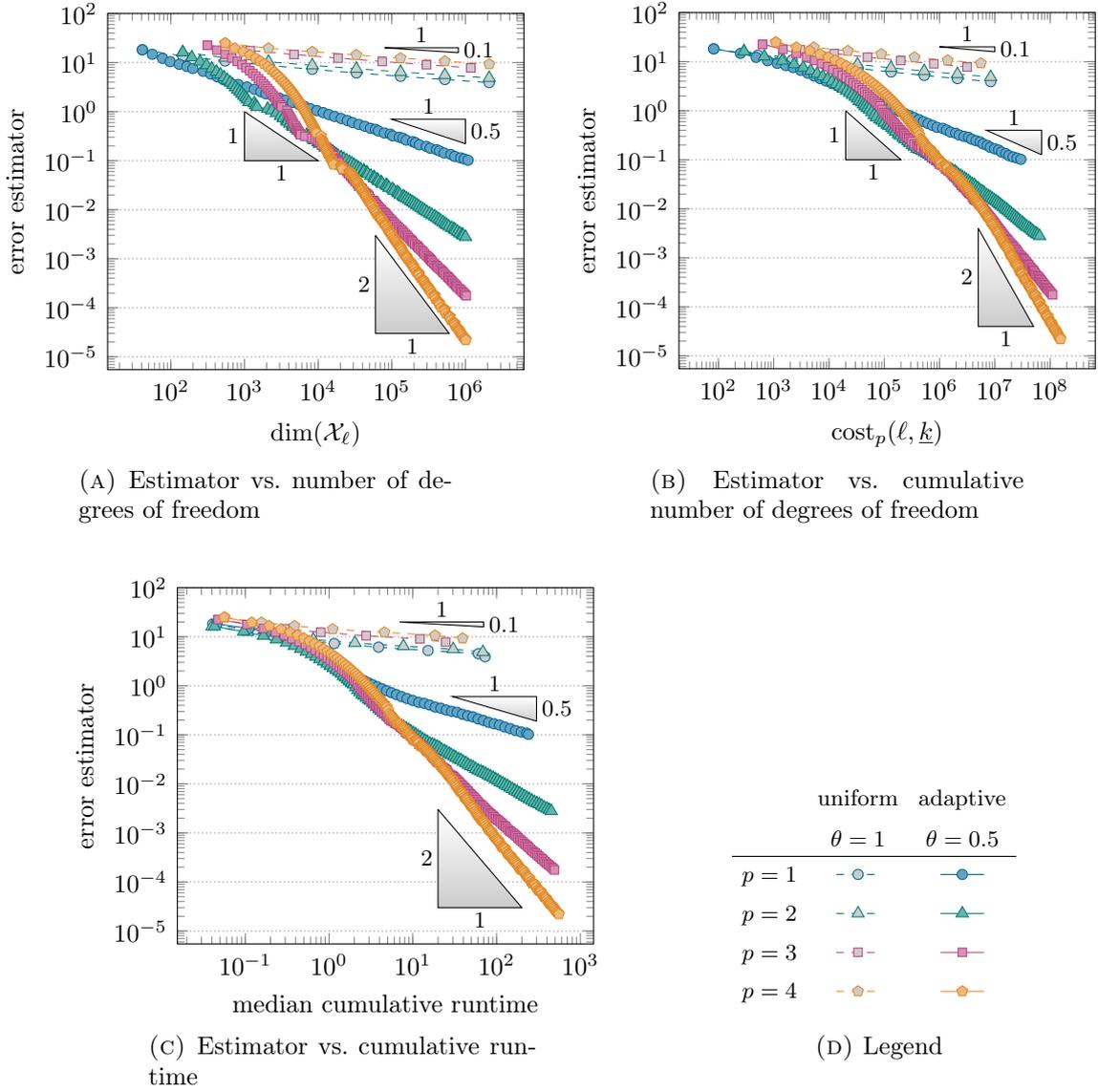

    \centering
    \subfloat[Estimator vs.\ number of degrees of freedom]{%
        \label{BMP:fig:Kellogg:convergence_nDofs}
        \input{chapter_afem_BMP/figures/Fig5a_Kellogg_problem_convergence_ndof.tex}
    }
    \hfil
    \subfloat[Estimator vs.\ cumulative number of degrees of freedom]{%
        \label{BMP:fig:Kellogg:convergence_cum_nDofs}
        \input{chapter_afem_BMP/figures/Fig5b_Kellogg_problem_convergence_cumndof.tex}
    }
    \medskip

    \subfloat[Estimator vs.\ cumulative runtime]{%
        \label{BMP:fig:Kellogg:convergence_runtime}
        \input{chapter_afem_BMP/figures/Fig5c_Kellogg_problem_convergence_cumruntime.tex}
    }
    \hfil
    \subfloat[Legend]{%
        \label{BMP:fig:Kellogg:legend}
        \begin{tikzpicture}[>=stealth]
    \colorlet{col1}{TUblue}
    \colorlet{col2}{TUgreen}
    \colorlet{col3}{TUmagenta}
    \colorlet{col4}{TUyellow}
    \colorlet{col5}{purple}
    \pgfplotsset{%
        degdefault/.style = {%
            mark = *,%
            mark size = 2pt,%
            every mark/.append style = {solid},%
            gray,%
            every mark/.append style = {fill = gray!60!white}%
        },%
        deg1/.style = {%
            degdefault,%
            col1,%
            every mark/.append style = {fill = col1!60!white}%
        },%
        deg2/.style = {%
            degdefault,%
            mark = triangle*,%
            mark size = 2.75pt,%
            col2,%
            every mark/.append style = {fill = col2!60!white}%
        },%
        deg3/.style = {%
            degdefault,%
            mark = square*,%
            mark size = 1.66pt,%
            col3,%
            every mark/.append style = {fill = col3!60!white}%
        },%
        deg4/.style = {%
            degdefault,%
            mark = pentagon*,%
            mark size = 2.2pt,%
            col4,%
            every mark/.append style = {fill = col4!60!white}%
        },%
        deg5/.style = {%
            degdefault,%
            mark = diamond*,%
            mark size = 2.75pt,%
            col5,%
            every mark/.append style = {fill = col5!60!white}%
        },%
        uniform/.style = {%
            dashed,%
            every mark/.append style = {%
                fill = black!20!white
            }%
        },%
        adaptive/.style = {%
            solid%
        },%
    }

    \matrix [
        matrix of nodes,
        anchor = south,
        font = \scriptsize,
        column 1/.style={anchor=base east},
    ] at (0,0) {
        & uniform
        & adaptive
        \\
        & \(\theta = 1\)
        & \(\theta = 0.5\)
        \\
        \hline \\
        \(p=1\)
        & \ref*{leg:est:unif:1}
        & \ref*{leg:est:adap:1}
        \\
        \(p=2\)
        & \ref*{leg:est:unif:2}
        & \ref*{leg:est:adap:2}
        \\
        \(p=3\)
        & \ref*{leg:est:unif:3}
        & \ref*{leg:est:adap:3}
        \\
        \(p=4\)
        & \ref*{leg:est:unif:4}
        & \ref*{leg:est:adap:4}
        \\
    };
\end{tikzpicture}
    }

    \caption{%
        Convergence plot of Algorithm~\ref{BMP:algorithm:afem}
        to solve the Kellogg benchmark problem from Subsection~\ref{BMP:sec:numerics}
        for various polynomial degrees.
        The adaptivity parameters read \(\theta = 0.5\) and \(\lambdaalg = 0.01\).
        All graphs display values of the residual-based error estimator 
        \(\eta_\ell (u_\ell^{\kk})\)
        from~\eqref{BMP:eq:estimator}.
    }
    \label{BMP:fig:Kellogg:convergence}
\end{figure}

The notion of optimal convergence rates relates to some
measure of computational costs.
While the motivation for studying the number of degrees of 
freedom stems from the goal of using meshes as coarse as possible,
the cumulative cost from~\eqref{BMP:eq:cost_p}
take into account all the work of computing a 
solution on an adaptive mesh.
In practice, one is more interested in the cumulative
runtime to compute a solution of a desired accuracy.
All three notions of optimality are investigated
in Figure~\ref{BMP:fig:Kellogg:convergence}.
The timing results throughout this chapter display the 
median of the total runtime obtained from 5 identical runs.
When the adaptive algorithm is implemented in linear complexity,
as in this case, the
Figures~\ref{BMP:fig:Kellogg:convergence_nDofs}--\ref{BMP:fig:Kellogg:convergence_runtime}
exhibit the same optimal rate \(-p/2\) of estimator decrease
for fixed polynomial degree~\(p\).
We underline that the rates for adaptive refinement are significantly larger
than the suboptimal rate of \(0.1\) attained with uniform refinement
for any polynomial degree \(p\).

\begin{figure}
    \centering
    \input{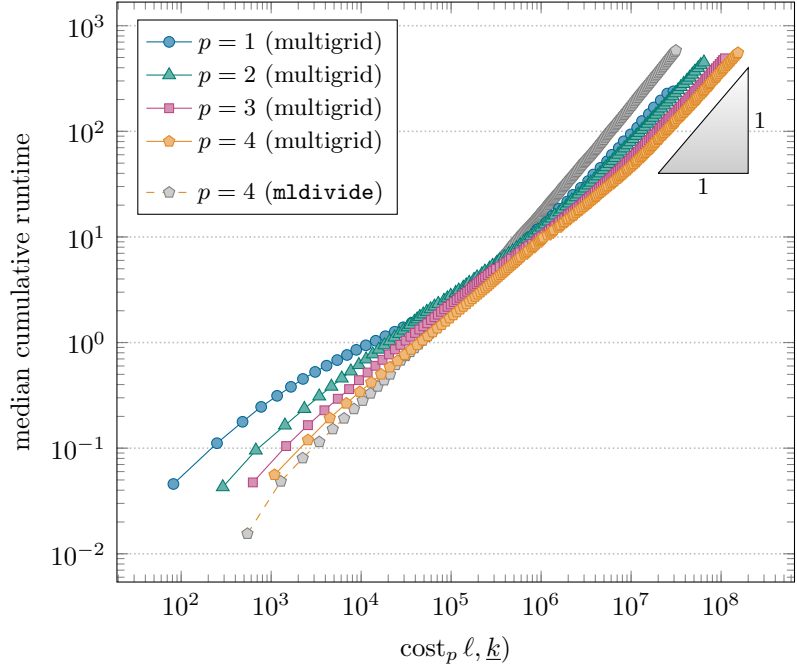}
    \caption{%
        Plot of the cumulative runtime in
        Algorithm~\ref{BMP:algorithm:afem}
        with \texttt{mldivide} and iterative algebraic solvers
        versus the cumulative number of degrees of freedom
        to solve the benchmark problem from
        Subsection~\ref{BMP:sec:numerics}
        for various polynomial degrees.
        The chosen adaptivity parameters read \(\theta = 0.5\)
        and \(\lambdaalg = 0.01\).
        The plot displays the runtime of one out of five
        identical runs chosen by the median of the total runtime.
    }
    \label{BMP:fig:Kellogg:runtime}
\end{figure}

Figure~\ref{BMP:fig:Kellogg:runtime} confirms the linear complexity because, 
independently of~\(p\), 
the slope of total runtime versus the cumulative cost attains one.
Moreover, already at \(10^4\) degrees of freedom,
the runtime of the adaptive algorithm with the iterative algebraic solver
is lower than with \textsc{Matlab}'s highly optimized 
built-in direct solver \texttt{mldivide}. 
This showcases from a practical point of view why 
optimal iterative solvers are indispensable for linear complexity of 
Algorithm~\ref{BMP:algorithm:afem}. 

\begin{figure}
    \centering
    \input{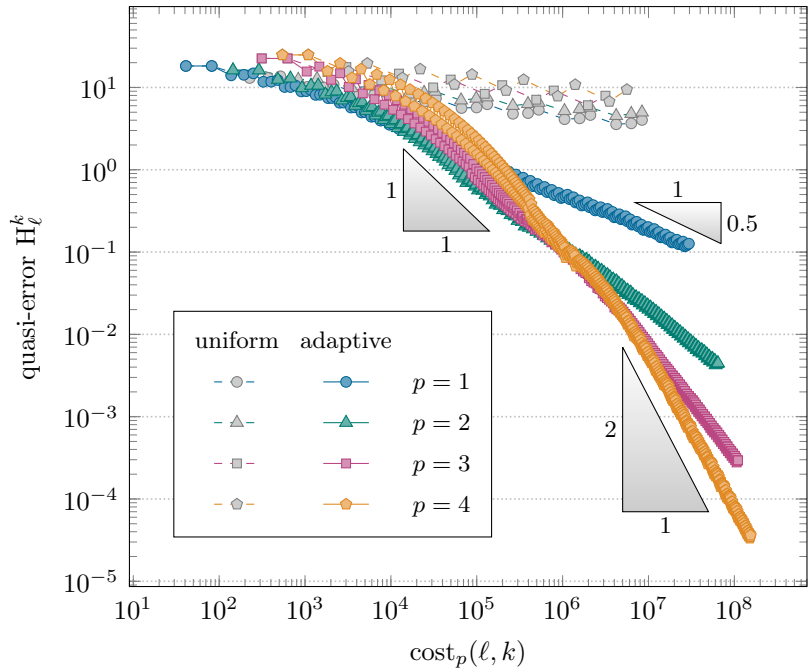}
    \caption{%
        Convergence plot of Algorithm~\ref{BMP:algorithm:afem}
        for various polynomial degrees.
        The chosen adaptivity parameters read \(\theta = 0.5\)
        and \(\lambdaalg = 0.01\).
        The graphs display the values of the quasi-error \(\Eta_\ell^k\)
        from~\eqref{BMP:eq:quasi-error} in each algebraic iteration step.
    }
    \label{BMP:fig:Kellogg:quasi_error}
\end{figure}

In Figure~\ref{BMP:fig:Kellogg:quasi_error}, we investigate 
full R-linear convergence of the quasi-error \(\Eta_\ell^{k}\)
from~\eqref{BMP:eq:quasi-error} as stated in 
Theorem~\ref{BMP:thm:full-linear-convergence}.
The quasi-error is reduced in each step of the 
adaptive algorithm
only up to a multiplicative constant.
Indeed, \(\Eta_\ell^{k}\) slightly increases in algebraic
solver steps, but this increase is limited and compensated
by the reduction of the mesh refinement and D\"orfler's marking strategy, 
leading to overall optimal convergence.

\begin{figure}
    \centering
    \input{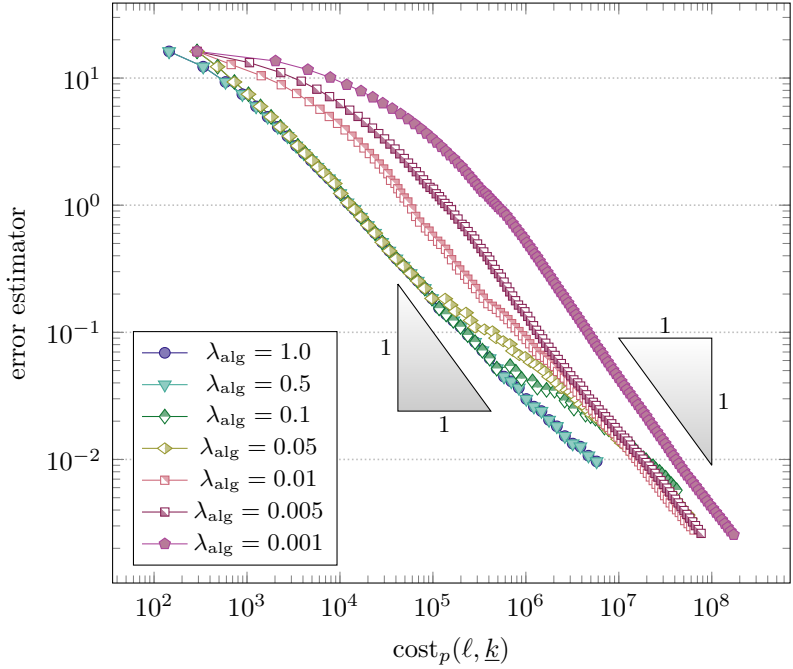}
    \caption{%
        Convergence plot of 
        Algorithm~\ref{BMP:algorithm:afem}
        to solve the Kellogg benchmark problem
        from Subsection~\ref{BMP:sec:numerics}
        for various choices
        of the stopping parameter \(\lambdaalg\).
        The remaining adaptivity parameter 
        is chosen as \(\theta = 0.3\)
        and the polynomial degree as \(p = 2\).
        The graphs display the values of
        the residual-based error estimator 
        \(\eta_\ell (u_\ell^{\kk})\) from~\eqref{BMP:eq:estimator}.
    }
    \label{BMP:fig:Kellogg:lambda}
\end{figure}

While Theorem~\ref{BMP:thm:optimal_complexity} ensures optimal 
complexity for sufficiently small adaptivity parameters, 
Figure~\ref{BMP:fig:Kellogg:theta} illustrates that also
moderate values of \( \theta\) lead to optimal convergence rates.
However, Figure~\ref{BMP:fig:Kellogg:lambda} indicates that smaller values of
\(\lambdaalg\) are required to achieve the optimal convergence rate
as the rates decrease for large \(\lambdaalg\) at about \(\cost(\ell, \kk) \geq 10^5\).
This observation might be caused by the strong singularity in this benchmark problem
and justifies the choice of \(\lambdaalg = 0.01\) in the remaining experiments.
The joint influence of the adaptivity parameters \(\theta\) and \(\lambdaalg\)
is investigated in Table~\ref{BMP:tab:Kellogg:parameters}, which 
displays the weighted cumulative runtime
\begin{equation}
    \label{BMP:eq:weighted_cumulative_runtime}
    \eta_\ell(u_\ell^{\kk})
    \sum_{\vert \ell', k' \vert \leq \vert \ell, \kk \vert}
    \operatorname{time}(\ell', k')
\end{equation}
with \(\eta_\ell(u_\ell^{\kk})/\eta_\ell(u_0^{0}) < 10^{-2}\) 
as stopping criterion.
The optimal choices with regard to this empirical balance of error estimator 
and computational runtime are moderate values 
\(0.5 \leq \theta \leq 0.7\) and \(0.5 \leq \lambdaalg = 0.9\).

\begin{figure}
    \centering
    \input{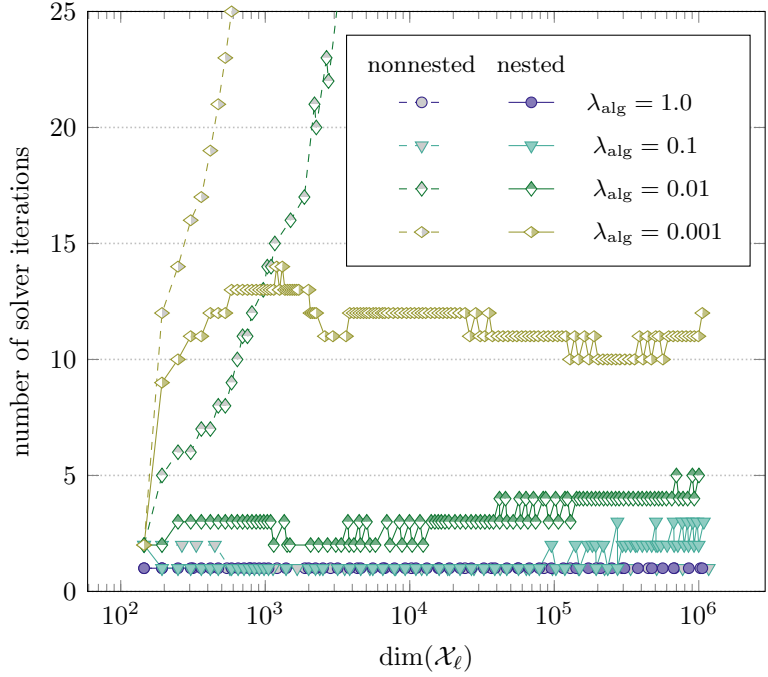}
    \caption{%
        Number of iteration steps
        of algebraic solver
        for the Kellogg benchmark problem from
        Subsection~\ref{BMP:sec:numerics}
        with nested (i.e., \(u_{\ell+1}^0 \coloneqq u_{\ell}^{\kk}\)) and 
        non-nested (i.e., \(u_{\ell+1}^0 \coloneqq 0 \)) iteration
        for various stopping parameters \(\lambdaalg\).
        The chosen polynomial degree reads \(p = 2\) and the 
        adaptivity parameter \(\theta = 0.5\).
    }
    \label{BMP:fig:Kellogg:iterations}
\end{figure}

Figure~\ref{BMP:fig:Kellogg:iterations} reveals that the use of nested
iteration in Algorithm~\ref{BMP:algorithm:afem}
is crucial to preserve a uniform number of algebraic solver
iterations.
Moreover, the plot illustrates how large parameters in the 
stopping criterion~\eqref{BMP:eq:j_stopping_criterion} lead to fewer algebraic 
iterations as the discretization error becomes the dominating
error source.

\begin{figure}
    \centering
    \input{chapter_afem_BMP/figures/Fig8_Kellogg_problem_theta.tex}
    \caption{%
        Convergence plot of Algorithm~\ref{BMP:algorithm:afem}
        to solve the Kellogg benchmark problem
        from Subsection~\ref{BMP:sec:numerics}
        for various choices
        of the bulk parameter \(\theta\).
        The remaining adaptivity parameter 
        is chosen as \(\lambdaalg = 0.01\)
        and the polynomial degree as \(p = 2\).
        The graphs display the values of
        the residual-based error estimator 
        \(\eta_\ell (u_\ell^{\kk})\) from~\eqref{BMP:eq:estimator}.
    }
    \label{BMP:fig:Kellogg:theta}
\end{figure}

\begin{table}
    \centering
    {
        %
%
\pgfplotstableset{%
    trim cells = true,%
    col sep = comma,%
    row sep = newline%
}
%
%
\pgfplotstableset{%
    fixed zerofill,%
    column type = {r},%
    every head row/.style = {before row = \toprule, after row = \midrule},%
    every last row/.style = {after row = \bottomrule},%
    columns/lambda/.style = {%
        fixed,%
        precision = 1,%
        column name={\diagbox[height=1.5\line]{$\lambdaalg$}{$\theta$}},%
    },%
    columns/theta0.9/.style = {column name = {$0.9$}},%
    columns/theta0.7/.style = {column name = {$0.7$}},%
    columns/theta0.5/.style = {column name = {$0.5$}},%
    columns/theta0.3/.style = {column name = {$0.3$}},%
    columns/theta0.1/.style = {column name = {$0.1$}},%
    rowmin/.style = {%
        postproc cell content/.append style={
            /pgfplots/table/@cell content/.add={\cellcolor{colorrowmin}}{}
        },
    },
    colmin/.style = {%
        postproc cell content/.append style={
            /pgfplots/table/@cell content/.add={\cellcolor{colorcolmin}}{}
        }
    },
    bothmin/.style = {%
        postproc cell content/.append style={
            /pgfplots/table/@cell content/.add={\cellcolor{colorbothmin}}{}
        }
    },
    highlight col min/.code 2 args = {%
        \pgfmathtruncatemacro\rowindex{#1-1}
        \edef\setstyles{%
            \noexpand\pgfplotstableset{%
                every row \rowindex\noexpand\space column #2/.style={colmin}
            }%
        }\setstyles
    },
    highlight row min/.code 2 args = {%
        \pgfmathtruncatemacro\rowindex{#1-1}
        \edef\setstyles{%
            \noexpand\pgfplotstableset{%
                every row \rowindex\noexpand\space column #2/.style={rowmin}
            }%
        }\setstyles
    },
    highlight both min/.code 2 args = {%
        \pgfmathtruncatemacro\rowindex{#1-1}
        \edef\setstyles{%
            \noexpand\pgfplotstableset{%
                every row \rowindex\noexpand\space column #2/.style={bothmin}
            }%
        }\setstyles
    },
}

%
%
%
\pgfplotstableread{%
lambda,theta0.9,theta0.7,theta0.5,theta0.3,theta0.1
0.9,1.35766828267154,0.963477191293004,0.946439624786913,1.06041721916581,1.02740347686486
0.7,1.53083861801114,1.03044702831798,1.0009901159443,1.05877851154672,1.0275336222503
0.5,1.49070430168122,0.963811882366322,1.14031230278502,1.0666662747084,1.03770229608879
0.3,1.54644444929028,0.975306990212493,1.15632641891119,1.05248610859973,1.05255535132248
0.1,1.67088643807793,0.998110760461873,1.38727957518827,1.07222616509858,1.06594810713307
}\dataTableAFEM
%
%
\pgfplotstabletypeset[%
    highlight col min = {1}{1},
    highlight col min = {1}{2},
    highlight col min = {4}{4},
    highlight col min = {1}{5},
    highlight row min = {2}{3},
    highlight row min = {4}{2},
    highlight row min = {5}{2},
    highlight both min = {1}{3},
    highlight both min = {3}{2},
]{\dataTableAFEM}

    }
    \medskip

    \caption{%
        \normalfont\footnotesize\quad
        Investigation of the influence of the adaptivity parameters
        \(\theta\) and \(\lambdaalg\) in Algorithm~\ref{BMP:algorithm:afem}
        to solve the Kellogg benchmark problem from
        Subsection~\ref{BMP:sec:numerics}.
        The polynomial degree is chosen as \(p = 2\).
        The best choice per column
        is marked in yellow, per row in blue, and for both in green.
        For all \( \theta\), the best performance is observed for
        \(\lambdaalg = 0.9\), while overall best results are achieved 
        for \( \theta \in \{ 0.5, 0.7\}\).
    }
    \label{BMP:tab:Kellogg:parameters}
\end{table}

%
%

\section{Convergence analysis}
\label{BMP:sec:proofs}

\subsection{Proof of full R-Linear convergence 
(Theorem~\ref*{BMP:thm:full-linear-convergence})}
\label{BMP:sec:full_linear_convergence}

The following result, implicitly found in~\cite[Lemma~4.9]{cfpp2014},
is the key tool for the characterization
of R-linear convergence of adaptive algorithms. 
We refer to~\cite[Lemma~10]{bfmps2023} for the present statement.
\begin{lemma}[tail summability vs.\ R-linear convergence]
    \label{BMP:lem:summability}
    Let \((\alpha_k)_{k \in \N_0}\)
    with \(\alpha_k \in \R_{\ge 0}\) for all \(k \in \N_0\).
    Then, the following statements are equivalent:
    \begin{enumerate}[label={\bf (\roman*)}]
        \item
            \emph{Tail summability:}
            There exists \(C_{\textup{sum}} > 0\)
            such that for all
            \(M,N\in \N_0\), 
            \begin{align}
            \label{BMP:eq:tail-summability}
                \sum_{k = M + 1}^{M + N} \alpha_k
                \le
                C_{\textup{sum}} \, \alpha_M.
            \end{align}
        \item
            \emph{R-linear convergence:}
            There exist \(0 < q_{\textup{conv}} < 1\)
            and \(C_{\textup{conv}} > 0\)
            such that for all
            \(M, N \in \N_0\)
            \begin{align}
            \label{BMP:eq:R-linear-convergence}
                \alpha_{M + N}
                \le
                C_{\textup{conv}}\, q_{\textup{conv}}^{N }\,
                \alpha_M. 
            \end{align}
    \end{enumerate}
\end{lemma}

\begin{proof}
{\bf (i)} \(\Rightarrow\) {\bf (ii)}: Note that for any 
\(M,N\in \N_0\), 
tail summability~\eqref{BMP:eq:tail-summability} yields
\[
    ( C_{\textup{sum}}^{-1} +1) \sum_{k = M + 1}^{M + N}   \alpha_k 
    = 
    C_{\textup{sum}}^{-1} \sum_{k = M + 1}^{M + N}   \alpha_k 
    + \! \sum_{k = M + 1}^{M + N}   \alpha_k  
    \le  \alpha_{M}  + \! \sum_{k = M + 1}^{M + N}    \alpha_k 
    =  \sum_{k = M }^{M + N}    \alpha_k.  
\]
Repeat this estimate inductively and use 
tail summability~\eqref{BMP:eq:tail-summability} to obtain
\begin{align*}
    \alpha_{M + N}
    &=
    \sum_{k = M + N}^{M + N}   \alpha_k 
    \le ( C_{\textup{sum}}^{-1} +1)^{-1}
    \sum_{k = M + N - 1}^{M + N} \alpha_k  
    \le \ldots
    \\
    &
    \le ( C_{\textup{sum}}^{-1} +1)^{-N}
    \sum_{k = M}^{M + N}   \alpha_k
    \le 
    ( C_{\textup{sum}}+1)
    ( C_{\textup{sum}}^{-1} +1)^{-N} \, \alpha_M.
\end{align*}
This proves the claim for 
\(q_{\textup{conv}} = ( C_{\textup{sum}}^{-1}+1)^{-1} \in (0,1)\)
and \(C_{\textup{conv}} = C_{\textup{sum}}+1\).

{\bf (ii)} \(\Rightarrow\) {\bf (i)}: Using R-linear 
convergence~\eqref{BMP:eq:R-linear-convergence}
for each term of the sum and the geometric series, we derive
\[
    \sum_{k = M + 1}^{M +N}   \alpha_k 
    =  \sum_{i = 1}^{N}   \alpha_{M + i} 
    \le C_{\textup{conv}} \sum_{i = 1}^{N} q_{\textup{conv}}^{i} \alpha_{M}
    \le\frac{C_{\textup{conv}} q_{\textup{conv}}}{1 
    - q_{\textup{conv}}} \alpha_{M}.
\]
This proves the claim for 
\(C_{\textup{sum}} = C_{\textup{conv}}q_{\textup{conv}}(1 
- q_{\textup{conv}})^{-1}\).
\end{proof}

\begin{proof}[Proof of Theorem~\ref{BMP:thm:full-linear-convergence}]
    The proof consists of five steps.

    \emph{Step~1} (contraction up to tail summable remainder).
    Fix 
    \(
        0
        <
        \mu
        <
        (1 - \qctr^2) \, \qctr^{-2}
    \)
    so that
    \(q_\mu \coloneqq (1 + \mu) \qctr^2 < 1\)
    and
    \(C_\mu \coloneqq (1 + \mu^{-1}) \qctr^2 > 0\). Then,
    the quasi-error contraction~\eqref{BMP:eq:contr-quasierr} and the Young 
    inequality prove,
    for all \((\ell+1, \underline k) \in \QQ\),
    that
    \begin{align}
        \label{BMP:eq:flc:contraction_w_remainder}
        &(\Eta_{\ell+1})^2
        \le
        q_\mu \, (\Eta_\ell)^2
        +
        C_\mu \,
        \vvvert u_{\ell+1}^\star - u_\ell^\star \vvvert^2.
    \end{align}

    \emph{Step~2} (tail summability with respect to \(\ell\)).
    For all \((\ell^\prime, \underline k) \in \QQ\)
    with \(\ell^\prime < \underline\ell - 1\),
    the Galerkin orthogonality in the form~\eqref{BMP:eq:Pythagoras}, the 
    telescoping sum, reliability~\eqref{BMP:axiom:reliability}, 
    and stability~\eqref{BMP:axiom:stability}, prove that
    \begin{align*}
        \sum_{\ell = \ell^\prime}^{\underline\ell - 1}
            \vvvert u_{\ell+1}^\star - u_\ell^\star \vvvert^2
        \;\;
        &\stackrel{\mathclap{\eqref{BMP:eq:Pythagoras}}}{=}\;\;
        \sum_{\ell = \ell^\prime}^{\underline\ell - 1}
        \big[
            \vvvert u^\star - u_{\ell}^\star \vvvert^2
            -
            \vvvert u^\star - u_{\ell+1}^\star \vvvert^2
        \big]
        \le
        \vvvert u^\star - u_{\ell^\prime}^\star \vvvert^2
        \\
        &
        \stackrel{\mathclap{\eqref{BMP:axiom:reliability}}}{\lesssim}
        \eta_{\ell^\prime}(u_{\ell^\prime}^\star)^2
        \stackrel{\eqref{BMP:axiom:stability}}{\lesssim}
        \eta_{\ell^\prime}(u_{\ell^\prime}^{\underline k})^2
        +
        \vvvert
            u_{\ell^\prime}^\star - u_{\ell^\prime}^{\underline k}
        \vvvert^2
        \stackrel{\eqref{BMP:eq:quasi-error-final-iterates}}{\eqsim} \,
        (\Eta_{\ell^\prime})^2.
    \end{align*}
    The combination with Step~1 and the geometric series
    results in \emph{quadratic tail summability} of the form,
    for all \((\ell^\prime, \underline k) \in \QQ\)
    with \(\ell^\prime < \underline\ell - 1\),
    \begin{align*}
            \sum_{\ell = \ell^\prime+1}^{\underline\ell-1}
            (\Eta_\ell)^2
            &\stackrel{\eqref{BMP:eq:flc:contraction_w_remainder}}{\le}
            \;
            (\Eta_{\ell^\prime})^{2} \,
            \sum_{\ell = \ell^\prime}^{\underline\ell-2}
            q_{\mu\delta}^{\ell-\ell^\prime}
            +
            C_\mu
            \sum_{\ell = \ell^\prime}^{\underline\ell-2}
                \vvvert
                    u_{\ell+1}^\star - u_\ell^\star
                \vvvert^2
            \lesssim
            (\Eta_{\ell^\prime})^2.
    \end{align*}
    According to Lemma~\ref{BMP:lem:summability},
    this quadratic tail summability
    is equivalent to R-linear convergence, i.e.,
    there exist \(C_1 > 0\) and \(0 < q_1 < 1\) such that
    \[
        (\Eta_{\ell})^{2}
        \le
        C_1 \, q_1^{\ell - \ell^\prime} \,
        (\Eta_{\ell^\prime})^{2}
        \quad \text{for all } 
                (\ell, \underline k),
        (\ell^\prime, \underline k)
        \in
        \QQ \text{ with }
        0 \le \ell^\prime \le \ell.
    \]
    Taking the square-root of this quasi-contraction estimate
    and another application of Lemma~\ref{BMP:lem:summability}
    result in \emph{linear tail summability}
    \begin{align}
        \label{BMP:eq:flc:2}
        \sum_{\ell = \ell^\prime + 1}^{\underline\ell - 1}
        \Eta_\ell
        \lesssim
        \Eta_{\ell^\prime}
        \quad \text{for all } (\ell^\prime, \underline k) \in \QQ.
    \end{align}
    This proves tail summability of the final iterates
    \((\ell, \underline k) \in \QQ\)
    with respect to \(\ell\).
    It remains to prove summability with respect to 
    the algebraic solver step index \(k\).

    \emph{Step~3} (quasi-contraction of quasi-error for
    every iterate).
    Recall the quasi-error with respect to 
    every iterate
    \(\Eta_\ell^k = \vvvert u_\ell^\star - u_\ell^k \vvvert + \eta_\ell(u_\ell^k)\)
    for all \((\ell, k) \in \QQ\) from~\eqref{BMP:eq:quasi-error}
    and note that
    \begin{align}
        \label{BMP:eq:flc:equivalence_quasi_errors}
        \Eta_\ell^{\underline k}
        \eqsim
        \Eta_\ell,
    \end{align}
    where the hidden equivalence constants depend only on \(\gamma\).
    The contraction~\eqref{BMP:eq:algebra} of the algebraic
    solver in estimate~\eqref{BMP:eq:algebra:apost} establishes,
    for all \(1 \le k \le \underline k\),
    \begin{equation}
        \label{BMP:eq:flc:stability_iterates}
        \vvvert
            u_\ell^k - u_\ell^{k-1}
        \vvvert
        \stackrel{\eqref{BMP:eq:algebra:apost}}{\le}
        (1 + \qalg) \,
        \vvvert
            u_\ell^\star -
            u_\ell^{k-1}
        \vvvert
        \le
        2 \,
        \vvvert
            u_\ell^\star -
            u_\ell^{k-1}
        \vvvert.
    \end{equation}
    For all \((\ell,\underline k) \in \QQ\),
    this leads to
    \begin{align}
        \label{BMP:eq:flc:stability_quasi_error}
        \begin{split}
            \Eta_\ell^{\underline k}
            &=
            \vvvert
                u_\ell^\star - u_\ell^{\underline k}
            \vvvert
            +
            \eta_\ell(u_\ell^{\underline k})
            \stackrel{\eqref{BMP:axiom:stability}}{\lesssim}
            \vvvert
                u_\ell^\star - u_\ell^{\underline k}
            \vvvert
            +
            \vvvert
                u_\ell^{\underline k} - u_\ell^{\underline k-1}
            \vvvert
            +
            \eta_\ell(u_\ell^{\underline k-1})
            \\
            &\stackrel{\mathclap{\eqref{BMP:eq:quasi-error}}}{\le}
            \;
            \Eta_\ell^{\underline k-1}
            +
            2 \,
            \vvvert
                u_\ell^{\underline k} - u_\ell^{\underline k-1}
            \vvvert
            \stackrel{\eqref{BMP:eq:flc:stability_iterates}}{\le}
            5 \, \Eta_\ell^{\underline k-1}.
        \end{split}
    \end{align}
    For the iterates with
    \(0 \le k^\prime < k < \underline k[\ell]\),
    the failure of the stopping criterion~\eqref{BMP:eq:j_stopping_criterion}
    for the algebraic solver in
    Algorithm~\ref{BMP:algorithm:afem}~(I.b)
    guarantees
    \begin{align}
        \label{BMP:eq:flc:quasi_contraction_iterates}
        \begin{split}
            \Eta_\ell^k
            &=
            \vvvert u_\ell^\star - u_\ell^k \vvvert
            +
            \eta_\ell(u_\ell^k)
            \stackrel{\eqref{BMP:eq:j_stopping_criterion}}
            <
            \vvvert u_\ell^\star - u_\ell^{k} \vvvert
            +
            \lambdaalg^{-1} \,
            \vvvert u_\ell^k - u_\ell^{k-1} \vvvert
            \\
            &\stackrel{\mathclap{\eqref{BMP:eq:flc:stability_iterates}}}\le
            \vvvert u_\ell^\star - u_\ell^k \vvvert
            +
            2 \, \lambdaalg^{-1}
            \vvvert u_\ell^\star - u_\ell^{k-1} \vvvert
            \stackrel{\mathclap{\eqref{BMP:eq:algebra}}}{\le}
            (\qalg + 2 \, \lambdaalg^{-1}) \,
            \vvvert u_\ell^\star - u_\ell^{k-1} \vvvert
            \\
            &\stackrel{\eqref{BMP:eq:algebra}}{\le}
            \frac{
                \qalg + 2 \, \lambdaalg^{-1}
            }{
                \qalg
            } \,
            \qalg^{k - k^\prime} \,
            \vvvert u_\ell^\star - u_\ell^{k^\prime} \vvvert.
        \end{split}
    \end{align}
    For \(0 \le k^\prime < k = \underline k[\ell]\),
    the combination
    of~\eqref{BMP:eq:flc:stability_quasi_error}--\eqref{BMP:eq:flc:quasi_contraction_iterates}
    results in
    \begin{align}
        \label{BMP:eq:flc:quasi_contraction_final}
        \Eta_\ell^{\underline k}
        \stackrel{\eqref{BMP:eq:flc:stability_quasi_error}}{\lesssim}
        \Eta_\ell^{\underline k-1}
        \stackrel{\eqref{BMP:eq:flc:quasi_contraction_iterates}}{\lesssim}
        \qalg^{(\underline k-1) - k^\prime} \,
        \vvvert u_\ell^\star - u_\ell^{k^\prime} \vvvert
        =
        \frac{1}{\qalg} \,
        \qalg^{\underline k - k^\prime} \,
        \vvvert u_\ell^\star - u_\ell^{k^\prime} \vvvert.
    \end{align}
    Since
    \(
        \vvvert
            u_\ell^\star - u_\ell^{k^\prime}
        \vvvert
        \le
        \Eta_\ell^{k^\prime}
    \),
    the
    estimates~\eqref{BMP:eq:flc:quasi_contraction_iterates}--%
    \eqref{BMP:eq:flc:quasi_contraction_final} prove
    \begin{align}
        \label{BMP:eq:flc:quasi_contraction}
        \Eta_\ell^k
        \lesssim
        \qalg^{k - k^\prime}
        \Eta_\ell^{k^\prime}
        \quad \text{for all } (\ell, k) \in \QQ
    \text{ with } 0 \le k^\prime \le k \le \underline k[\ell].
    \end{align}

    In the remaining case of \(k = 0\), Céa-type estimate~\eqref{BMP:eq:cea},
    reliability~\eqref{BMP:axiom:reliability}, and stability~\eqref{BMP:axiom:stability} imply
    \begin{align*}
        \vvvert u_\ell^\star - u_{\ell-1}^\star \vvvert
        &\le
        \vvvert u_\ell^\star - u^\star \vvvert
        +
        \vvvert u^\star - u_{\ell-1}^\star \vvvert
        \stackrel{\eqref{BMP:eq:cea}}{\lesssim}
        \vvvert u^\star - u_{\ell-1}^\star \vvvert
        \stackrel{\eqref{BMP:axiom:reliability}}{\lesssim}
        \eta_{\ell-1}(u_{\ell-1}^\star)
        \\
        &\stackrel{\mathclap{\eqref{BMP:axiom:stability}}}{\lesssim}
        \eta_{\ell-1}(u_{\ell-1}^{\underline k})
        +
        \vvvert
            u_{\ell-1}^\star -
            u_{\ell-1}^{\underline k}
        \vvvert
        =
        \Eta_{\ell-1}^{\underline k}
        \stackrel{\eqref{BMP:eq:flc:equivalence_quasi_errors}}{\eqsim}
        \Eta_{\ell-1}.
    \end{align*}
    This and nested iteration
    \(
        u_\ell^0
        =
        u_{\ell-1}^{\underline k}
    \)
    yield
    \begin{equation}
        \begin{split}
            \label{BMP:eq3:step7}
            \Eta_\ell^0
            \;\;
            &\stackrel{\mathclap{\eqref{BMP:eq:quasi-error}}}{=}
            \;\;
            \vvvert
                u_\ell^\star - u_{\ell-1}^{\underline k}
            \vvvert
            +
            \eta_\ell(u_{\ell-1}^{\underline k})
            \stackrel{\eqref{BMP:eq:reduction_est}}\lesssim
            \vvvert
            u_\ell^\star - u_{\ell-1}^{\underline k}
            \vvvert
            +
            \eta_{\ell-1}(u_{\ell-1}^{\underline k})
            \\
            &
            \le
            \vvvert u_\ell^\star - u_{\ell-1}^\star \vvvert
            +
            \Eta_{\ell-1}^{\underline k}
            \stackrel{\eqref{BMP:eq:flc:equivalence_quasi_errors}}{\lesssim}
            \Eta_{\ell-1}
            \text{ for all }(\ell, 0) \in \QQ 
            \text{ with } \ell >0 .
        \end{split}
    \end{equation}

    \emph{Step~4} (tail summability with respect to \(\ell\) and \(k\)).
    Let \((\ell^\prime, k^\prime) \in \QQ\).
    The quasi-contraction estimates from Step~3
    and the geometric series prove that
    \begin{align}
        \label{BMP:eq:flc:summability_ell_k}
        \begin{split}
            \sum_{\substack{
                (\ell, k) \in \QQ
                \\
                \vert \ell, k \vert
                >
                \vert \ell^\prime, k^\prime \vert
            }}
            \!\!
            \Eta_\ell^k
            &=
            \sum_{k = k^\prime+1}^{\underline k[\ell^\prime]}
            \Eta_{\ell^\prime}^k
            +
            \sum_{\ell = \ell^\prime + 1}^{\underline\ell}
            \sum_{k=0}^{\underline k[\ell]}
            \Eta_\ell^k
            \stackrel{\eqref{BMP:eq:flc:quasi_contraction}}{\lesssim}
            \Eta_{\ell^\prime}^{k^\prime}
            +
            \sum_{\ell = \ell^\prime+1}^{\underline\ell}
            \Eta_\ell^0
            \\
            &\stackrel{\mathclap{\eqref{BMP:eq3:step7}}}{\lesssim}
            \;
            \Eta_{\ell^\prime}^{k^\prime}
            + \!
            \sum_{\ell = \ell^\prime}^{\underline\ell - 1}
            \Eta_\ell
            \!
            \stackrel{\eqref{BMP:eq:flc:2}}{\lesssim}
            \!
            \Eta_{\ell^\prime}^{k^\prime}
            +
            \Eta_{\ell^\prime}
            \!
            \stackrel{\eqref{BMP:eq:flc:equivalence_quasi_errors}}{\eqsim}
            \!
            \Eta_{\ell^\prime}^{k^\prime}
            +
            \Eta_{\ell^\prime}^{\underline k}
            \!
            \stackrel{\eqref{BMP:eq:flc:quasi_contraction}}{\lesssim}
            \!
            \Eta_{\ell^\prime}^{k^\prime}.
        \end{split}
    \end{align}

    \emph{Step~5} (full R-linear convergence).
    Since the index set \(\QQ\) is
    linearly ordered with respect to the total step counter
    \(| \cdot, \cdot |\),
    the tail summability~\eqref{BMP:eq:flc:summability_ell_k} and
    Lemma~\ref{BMP:lem:summability} conclude the proof, i.e.,
    there exist \(C_2 > 0\) and \(0 < q_2 < 1\) such that
    \[
        \Eta_\ell^{k}
        \le
        C_2 \,
        q_2^{
            \vert \ell, k \vert
            -
            \vert \ell^\prime, k^\prime \vert
        } \,
        \Eta_{\ell^\prime}^{k^\prime}
        \quad \text{for all } (\ell^\prime, k^\prime), (\ell, k) \in \QQ
        \text{ with }  \vert \ell, k \vert
        \ge
        \vert \ell^\prime, k^\prime \vert.
    \]

    Finally, we show the upper bound on the overall error by 
    the quasi-error in~\eqref{BMP:eq:full-linear-convergence-robust} 
    yielding parameter-robust convergence. Using 
    reliability~\eqref{BMP:axiom:reliability} and
    stability~\eqref{BMP:axiom:stability},
    there holds
    \begin{align*}
        \vvvert u^\star - u_\ell^k \vvvert
        &\le
        \vvvert u^\star - u_\ell^\star \vvvert 
        +
        \vvvert u_\ell^\star  - u_\ell^k \vvvert
        \stackrel{\eqref{BMP:axiom:reliability}}\lesssim 
        \eta_{\ell}(u_{\ell}^\star) 
        + \vvvert u_\ell^\star  - u_\ell^k \vvvert
        \\
        &
        \stackrel{\mathclap{\eqref{BMP:axiom:stability}}}\lesssim 
        \eta_{\ell}(u_{\ell}^k)  + \vvvert u_\ell^\star  - u_\ell^k \vvvert
        =
        \Eta_{\ell}^k. 
        \qedhere
    \end{align*}
\end{proof}

\subsection{Proof of optimal complexity \texorpdfstring{(Theorem~\ref{BMP:thm:optimal_complexity})}{(Theorem~\ref*{BMP:thm:optimal_complexity})}}
\label{BMP:sec:optimal_complexity}

The following corollary asserts that full R-linear convergence from
Theorem~\ref{BMP:thm:full-linear-convergence} implies
the coincidence of the convergence rates of the quasi-error
from~\eqref{BMP:eq:quasi-error} with respect to
the number of degrees of freedom and 
with respect to 
the total computational work.
Thus, full R-linear convergence is a key argument
to prove optimal complexity of Algorithm~\ref{BMP:algorithm:afem}.
\begin{corollary}[{\cite[Corollary~11]{bfmps2023}}]
    \label{BMP:cor:equivalence_rates}
    For \(s > 0\), abbreviate
    \begin{align}
    \label{BMP:eq:def:Ms}
        M(s)
        \coloneqq
        \sup_{\substack{(\ell, k) \in \QQ}}
        (\# \TT_\ell)^s \Eta_\ell^k.
    \end{align}
    Suppose that the assumptions
    of Theorem~\ref{BMP:thm:full-linear-convergence} hold.
    Then, Algorithm~\ref{BMP:algorithm:afem} guarantees 
    the equivalence
    \begin{align}
    \label{BMP:eq:equivalence_rates}
        M(s)
        \le
        \sup_{\substack{(\ell, k) \in \QQ}}
        \cost(\ell, k)^s\,
        \Eta_\ell^k
        \le
        \frac{\Clin}{[1 - \qlin^{1/s}]^s}
        \, M(s).
    \end{align}
\end{corollary}

\begin{proof}
    Recall from~\eqref{BMP:eq:cost} that 
    \[
        \cost(\ell, k)
        =
        \sum_{\substack{
            (\ell^\prime, k^\prime) \in \QQ
            \\
            \vert \ell^\prime, k^\prime \vert
            \le \vert \ell, k\vert
        }}
        \# \TT_{\ell'}.
    \]
    The first estimate in~\eqref{BMP:eq:equivalence_rates} is thus immediate 
    as it 
    suffices to bound the sum of positive terms from below 
    by one of its components.
    Next, note that 
    for all \(( \ell', k') \in \QQ \),
    the definition~\eqref{BMP:eq:def:Ms} provides
    \[
        \# \TT_{\ell' }\le M(s)^{1/s} (\Eta_{\ell'}^{k'})^{-1/s}.
    \]
    Summing and using the full R-linear 
    convergence~\eqref{BMP:eq:full-linear-convergence} together with 
    the geometric series gives
    \begin{align*}
        \cost(\ell, k)
        &=
        \sum_{\substack{
            (\ell^\prime, k^\prime) \in \QQ
            \\
            \vert \ell^\prime, k^\prime \vert
            \le \vert \ell, k\vert
        }} \! \!
        \# \TT_{\ell'}
        \le M(s)^{1/s}  \! \!
        \sum_{\substack{
            (\ell^\prime, k^\prime) \in \QQ
            \\
            \vert \ell^\prime, k^\prime \vert
            \le \vert \ell, k\vert
        }} \! \!
        (\Eta_{\ell'}^{k'})^{-1/s} 
        \\
        &
        \stackrel{\mathclap{\eqref{BMP:eq:full-linear-convergence}}}{\le}
        M(s)^{1/s}  \,\Clin^{1/s}  
        \sum_{\substack{
            (\ell^\prime, k^\prime) \in \QQ
            \\
            \vert \ell^\prime, k^\prime \vert
            \le \vert \ell, k\vert
        }} \! \!
            (\qlin^{1/s})^{
                \vert \ell, k \vert
                - \vert \ell^\prime, k^\prime \vert
            }  (\Eta_{\ell}^{k})^{-1/s} 
        \\
        &
        \le  M(s)^{1/s} 
        \frac{\Clin^{1/s}}{1 - \qlin^{1/s}}  \,
        (\Eta_{\ell}^{k})^{-1/s}. 
    \end{align*}  
    Rearranging the terms and taking the supremum over \(( \ell, k) \in \QQ \) 
    proves the second estimate in~\eqref{BMP:eq:equivalence_rates}.
\end{proof}

The following result relates the error estimator for the
final algebraic iterates to the one for the exact solution.
It imposes a restriction on the algebraic stopping parameter
\(\lambdaalg\) to ensure that the final iterate is 
\emph{sufficiently close} to the exact solution.
\begin{lemma}[estimator equivalence]
    \label{BMP:lem:estimator_equivalence}
    There exists
    \(
        \lambdaalg^\star 
        \coloneqq
        (1 - \qalg) / ( \Cstab\, \qalg)
    > 0
    \)
    such that every
    \(0 < \theta \leq 1\)
    and
    \(
        0
        <
        \lambdaalg
        <
        \lambdaalg^\star
    \)
    guarantee
    \begin{equation}
        \label{BMP:eq:estimator_equivalence}
        (1 - \lambdaalg /
        \lambdaalg^\star)
        \,
        \eta_\ell(u_\ell^{\kk})
        \leq
        \eta_{\ell}(u_\ell^\star)
        \leq
        (1 + \lambdaalg / \lambdaalg^\star)
        \,
        \eta_\ell(u_\ell^{\kk}).
    \end{equation}
\end{lemma}

\begin{proof}
    The a-posteriori error control in~\eqref{BMP:eq:algebra:apost} 
    and
    the stopping
    criterion~\eqref{BMP:eq:j_stopping_criterion}
    establish
    \begin{equation}
        \label{BMP:eq:stopping_estimate}
        \Cstab \,
        \vvvert u_\ell^\star - u_\ell^{\kk} \vvvert
        \leq  
        \frac{ \Cstab\, \qalg}{1 - \qalg}\,
        \vvvert u_\ell^{\kk} - u_\ell^{\kk - 1} \vvvert
        \leq
        \frac{\lambdaalg}{\lambdaalg^\star}\,
        \eta_\ell(u_\ell^{\kk}).
    \end{equation}
    For any $\lambdaalg > 0$, this and stability~\eqref{BMP:axiom:stability}
    prove the upper bound
    \[
        \eta_\ell(u_\ell^\star)
        \leq
        \eta_\ell(u_\ell^{\kk})
        +
        \Cstab\,
        \vvvert u_\ell^\star - u_\ell^{\kk} \vvvert
        \leq
        \big(
            1 + \lambdaalg / \lambdaalg^\star
        \big)\,
        \eta_\ell(u_\ell^{\kk}).
    \]
    Analogously, we prove
    \[
        \eta_\ell(u_\ell^{\kk})
        \leq
        \eta_\ell(u_\ell^\star)
        +
        \Cstab\,
        \vvvert u_\ell^\star - u_\ell^{\kk} \vvvert
        \leq
        \eta_\ell(u_\ell^\star)
        + \frac{\lambdaalg}{\lambdaalg^\star} \,
        \eta_\ell(u_\ell^{\kk}).
    \]
    Rearrangement yields,
    for all
    \(
        0 < \lambdaalg
        < \lambdaalg^\star
    \),
    that
    \[
        \big(1- \lambdaalg/\lambdaalg^\star\big)  \, 
        \eta_\ell(u_\ell^{\kk})
        \leq
        \eta_\ell(u_\ell^\star).
        \qedhere
    \]
\end{proof}

Given a particular triangulation \(\TT_\ell\),
the following lemma provides the existence of an optimal set for refinement
satisfying the D\"orfler criterion for a sufficiently small
bulk parameter~\(\theta_{\textup{mark}}\).
In the subsequent proof of
optimal complexity of Algorithm~\ref{BMP:algorithm:afem}, 
this allows to relate the actual marked set with these theoretical optimal sets
of refined elements.
\begin{lemma}[equivalence with optimal-practical 
    refinement]
    \label{BMP:lem:optimal_refinement}
    Given \(\lambdaalg^\star > 0\) from Lemma~\ref{BMP:lem:estimator_equivalence},
    choose \(0 < \theta^\star \leq 1\) sufficiently small to ensure
    that the parameters 
    \(0<\theta<\theta^\star\) and 
    \(0<\lambdaalg < \lambdaalg^\star\)
    satisfy
    \[
        0
        <
        \theta_{\textup{mark}}
        \coloneqq
        (\theta^{1/2}
        + \lambdaalg /
        \lambdaalg^\star)^2
        /
        (1 - \lambdaalg /
        \lambdaalg^\star)^2
        <
        (1 + \Cstab^2 \Cdrel^2)^{-1}.
    \]
    Then, there exists a set \(\RR_\ell \subseteq \TT_\ell\)
    such that
    \begin{equation}
        \label{BMP:eq:optimal_refinement}
        \# \RR_{\ell}
        \lesssim
        \Vert u^\star \Vert_{\A_s}^{1/s}
        \eta_{\ell}(u_{\ell}^\star)^{-1/s}
        \quad\text{and}\quad
        \theta_{\textup{mark}} \,
        \eta_{\ell}(u_{\ell}^\star)^2
        \leq
        \eta_{\ell}(\RR_{\ell}, u_{\ell}^\star)^2.
    \end{equation}
    Moreover, the above D\"orfler marking~\eqref{BMP:eq:optimal_refinement} 
    for \( (\theta_{\textup{mark}} , u_{\ell}^\star) \) implies 
    that \(\RR_{\ell}\) also satisfies the 
    D\"orfler marking~\eqref{BMP:eq:marking_criterion} for 
    \( (\theta, u_{\ell}^\kk) \), meaning that
    \begin{equation}
        \label{BMP:eq:practical_refinement}
        \theta \,
        \eta_{\ell}(u_{\ell}^\kk)^2
        \leq
        \eta_{\ell}(\RR_{\ell}, u_{\ell}^\kk)^2.
    \end{equation}
\end{lemma}

\begin{proof}
    The existence of the set \(\RR_{\ell}\) satisfying~\eqref{BMP:eq:optimal_refinement} 
    follows from~\cite[Lem.~4.14]{cfpp2014}, 
    because \(\theta_{\textup{mark}} < (1 + \Cstab^2 \Cdrel^2)^{-1}\) is 
    guaranteed by the choice of the parameters.
    To show~\eqref{BMP:eq:practical_refinement}, 
    we use the estimator equivalence from Lemma~\ref{BMP:lem:estimator_equivalence}
    and the marking criterion~\eqref{BMP:eq:optimal_refinement} to obtain
    \begin{align*}
        (1 - \lambdaalg / \lambdaalg^\star) \, \theta_{\textup{mark}}^{1/2}
        \,
        \eta_\ell(u_\ell^{\kk})
        \stackrel{\eqref{BMP:eq:estimator_equivalence}}{\le}
        \theta_{\textup{mark}}^{1/2} \,
        \eta_{\ell}(u_\ell^\star)
        \stackrel{\eqref{BMP:eq:optimal_refinement}}{\le}
        \eta_{\ell}(\RR_{\ell}, u_{\ell}^\star).
    \end{align*}
    Next, using stability~\eqref{BMP:axiom:stability}, a-posteriori control~\eqref{BMP:eq:algebra:apost}, and the stopping criterion~\eqref{BMP:eq:j_stopping_criterion}, we deduce
    \begin{align*}
        \eta_{\ell}(\RR_{\ell}, u_{\ell}^\star)
        &
        \stackrel{\mathclap{\eqref{BMP:axiom:stability}}}{\le}
        \eta_{\ell}(\RR_{\ell}, u_\ell^{\kk})
        + \Cstab  \vvvert u_{\ell}^\star - u_\ell^{\kk} \vvvert
        \\
        &
        \stackrel{\mathclap{\eqref{BMP:eq:algebra:apost}}}{\le}
        \eta_{\ell}(\RR_{\ell}, u_\ell^{\kk}) 
        + \Cstab\,  \frac{\qalg}{1 - \qalg} \,  
        \vvvert u_{\ell}^{\kk} - u_\ell^{\kk-1} \vvvert 
        \\
        &
        \stackrel{\mathclap{\eqref{BMP:eq:j_stopping_criterion}}}{\le}
        \eta_{\ell}(\RR_{\ell}, u_\ell^{\kk})
        + \lambdaalg/\lambdaalg^\star \, \eta_{\ell}(u_\ell^{\kk}).
    \end{align*}
    The combination of the two previous estimates
    with the definition of \(\theta_{\textup{mark}}\) 
    yields
    \[
        (1 - \lambdaalg /
        \lambdaalg^\star) \, \theta_{\textup{mark}}^{1/2}
        \,
        \eta_\ell(u_\ell^{\kk})
        \le
        \eta_{\ell}(\RR_{\ell}, u_\ell^{\kk}) 
        + 
        \big[ 
            \theta_{\textup{mark}}^{1/2} 
            (
                1 - \lambdaalg/\lambdaalg^\star
            )
            - \theta^{1/2}  
        \big]
        \,\eta_{\ell}(u_\ell^{\kk}).
    \]
    Rearranging the terms gives the practical
    marking criterion~\eqref{BMP:eq:practical_refinement}.
\end{proof}

\begin{remark}[optimality of D\"orfler marking]
    Though the details on the set \(\RR_{\ell}\) 
    satisfying~\eqref{BMP:eq:optimal_refinement} are not fully discussed 
    in the previous proof, note 
    that this provides a central argument for using D\"orfler marking 
    criterion~\eqref{BMP:eq:marking_criterion} in adaptive FEMs. 
    Indeed, referring to, e.g., \cite[Lem.~4.14]{cfpp2014} and the seminal 
    paper~\cite{stevenson2007}, 
    let us suppose that, for a given sequence of successively refined 
    triangulations \( (\TT_\ell )_{\ell \in \N_0}\), we compute the exact 
    discrete solutions \( u_\ell^\star \) and that there holds linear 
    convergence 
    \( \eta_{\ell+n}(u_{\ell+n}^\star ) 
    \le C \qlin^n  \eta_\ell(u_\ell^\star) \) for all \( \ell,n \in \N_0 \), 
    where \( C > 0 \) 
    and \( 0 < \qlin < 1 \) are \(\ell\)-independent. 
    Then, there exists \(0<\theta_\star \le 1 \) 
    such that, for any \(0<\theta < \theta_\star \),  
    the set \(\RR_\ell \coloneqq \TT_\ell \setminus  \TT_{\ell+1} \) 
    satisfies the D\"orfler marking 
    \begin{align}\label{BMP:rem:marking_criterion}
        \theta \, \eta_{\ell}(u_\ell^\star)^2 
        \le \eta_{\ell}(\RR_{\ell}, u_\ell^\star)^2
        \quad \text{for all } \ell \in \N_0.
    \end{align}
    Independently of how the marked elements are chosen, D\"orfler 
    marking is thus satisfied along the sequence of refined elements if 
    there holds R-linear convergence. 
    In this sense, D\"orfler marking is 
    not only sufficient for (full) R-linear convergence, but also 
    necessary.
\end{remark}

\begin{proof}[Proof of Theorem~\ref{BMP:thm:optimal_complexity}]
    By Corollary~\ref{BMP:cor:equivalence_rates},
    it suffices to prove, for every \(s > 0\),
    \[
        \sup_{(\ell, k) \in \QQ}
        (\# \TT_\ell)^s
        \Eta_\ell^k
        \lesssim
        \max\{
            \Vert u^\star \Vert_{\A_s},
            \Eta_0^0
        \}.
    \]
    Let \(s>0\) and assume that the right-hand side of the previous 
    estimate is finite, because otherwise the claim is trivial.
    The proof consists of three steps.

    \emph{Step~1.}\quad
    For all \(0 \leq \ell' < \elll\),
    Lemma~\ref{BMP:lem:optimal_refinement} guarantees existence
    of a set \(\RR_{\ell'} \subseteq \TT_{\ell'}\)
    satisfying~\eqref{BMP:eq:optimal_refinement}
    for \(\ell'\) replacing \(\ell\).
    This, the estimate~\eqref{BMP:eq3:step7}, a-posteriori control 
    of the algebraic error~\eqref{BMP:eq:algebra:apost},
    the stopping criterion~\eqref{BMP:eq:j_stopping_criterion},
    and the estimator equivalence from
    Lemma~\ref{BMP:lem:estimator_equivalence}
    imply
    \[
        \Eta_{\ell'+1}^0
        \stackrel{\eqref{BMP:eq3:step7}}{\lesssim}
        \Eta_{\ell'}
        \stackrel{\eqref{BMP:eq:algebra:apost}}{\lesssim}
        \vvvert u_{\ell'}^{\kk} - u_{\ell'}^{\kk-1} \vvvert
        +
        \eta_{\ell'}(u_{\ell'}^{\kk})
        \stackrel{\eqref{BMP:eq:j_stopping_criterion}}{\lesssim}
        \eta_{\ell'}(u_{\ell'}^{\kk}) \,
        \stackrel{\eqref{BMP:eq:estimator_equivalence}}{\eqsim} \,
        \eta_{\ell'}(u_{\ell'}^\star).
    \]
    The combination with~\eqref{BMP:eq:optimal_refinement} proves
    \[
        \# \RR_{\ell'}
        \lesssim
        \Vert u^\star \Vert_{\A_s}^{1/s}
        (\Eta_{\ell'+1}^0)^{-1/s}.
    \]

    \emph{Step~2}.\quad
    Since \(\RR_{\ell'} \subseteq \TT_{\ell'}\),
    the quasi-minimality of \(\MM_{\ell'}\)
    implies \(\#\MM_{\ell'} \lesssim
    \#\RR_{\ell'}\) due to the comparison of optimal and practical marking~\eqref{BMP:eq:practical_refinement}. 
    This and the mesh-closure estimate~\eqref{BMP:axiom:closure}
    show
    \[
        \#\TT_\ell - \#\TT_0
        \lesssim
        \sum_{\ell'=0}^\ell
        \#\MM_{\ell'}
        \lesssim
        \sum_{\ell'=0}^\ell
        \#\RR_{\ell'}.
    \]

    \emph{Step~3}.\quad
    Full R-linear
    convergence~\eqref{BMP:eq:full-linear-convergence}
    proves
    \[
        \sum_{\ell'=0}^{\ell-1}
        (\Eta_{\ell'+1}^0)^{-1/s}
        \leq
        \sum_{ \substack{(\ell', k') \in \QQ \\
                \vert \ell', k' \vert \leq \vert \ell, k \vert
        } }
        (\Eta_{\ell'}^{k'})^{-1/s}
        \lesssim
        (\Eta_{\ell}^{k})^{-1/s}
        \sum_{ \substack{(\ell', k') \in \QQ \\
                \vert \ell', k' \vert \leq \vert \ell, k \vert
        } }
        (\qlin)^{\vert \ell, k \vert - \vert \ell', k'\vert}
        \lesssim
        \Vert u^\star \Vert_{\A_s}^{1/s}
        (\Eta_{\ell}^{k})^{-1/s}.
    \]
    The combination of this with the Steps~1--2 reads
    \[
        \#\TT_\ell - \#\TT_0
        \lesssim
        \Vert u^\star \Vert_{\A_s}^{1/s}
        (\Eta_{\ell}^{k})^{-1/s}.
    \]
    For \(\ell > 0\),
    the lower bound
    \(1 \leq \#\TT_\ell - \#\TT_0\) implies
    \(\#\TT_\ell - \#\TT_0 + 1 \leq 2(\#\TT_\ell - \#\TT_0)\)
    and thus
    \[
        (\#\TT_\ell - \#\TT_0 + 1)^s\,
        \Eta_{\ell}^{k}
        \lesssim
        \Vert u^\star \Vert_{\A_s}.
    \]
    For \(\ell = 0\), full R-linear
    convergence~\eqref{BMP:eq:full-linear-convergence}
    yields
    \[
        (\#\TT_0 - \#\TT_0 + 1)^s\,
        \Eta_{0}^{k}
        \lesssim
        \Eta_{0}^{0}.
    \]
    The combination with the elementary estimate
    \(
        \#\TT_\ell
        \leq
        \#\TT_0(\#\TT_\ell - \#\TT_0 + 1)
    \)
    from
    \cite[Lemma~22]{bhp2017}
    concludes the proof
    \[
        (\#\TT_\ell)^s\, \Eta_\ell^k
        \lesssim
        \max\{
            \Vert u^\star \Vert_{\A_s},
            \Eta_0^0
        \}
        \quad \text{for all } (\ell, k) \in \QQ.
        \qedhere
    \]
\end{proof}

%
%

\section{Extension 1: Goal-oriented adaptivity}
\label{BMP:sec:goafem}

The Algorithm~\ref{BMP:algorithm:afem} from Subsection~\ref{BMP:sec:AFEM}
aims at the efficient approximation
of the solution \(u^\star\) with respect to
the energy norm \(\vvvert \cdot \vvvert\).
In many applications, one rather looks for
a particular quantity of interest \(G(u^\star)\)
in terms of a linear goal functional
\(G: \XX \to \R\).
A non-standard marking procedure allows
to approximate the goal functional with the doubled
optimal convergence rate.
In this section, we present some results 
from~\cite{bgip2023} on cost-optimal 
goal-oriented adaptivity with inexact solver.

\subsection{Adaptive FEM with linear goal functional}

Let \(u^\star \in H^1_0 (\Omega)\) be the weak solution
to~\eqref{BMP:eq:weak}.
Given \(g \in L^2(\Omega)\)
and \(\boldsymbol{g} \in [L^2(\Omega)]^d\), we consider the linear
functional
\begin{equation}
    \label{BMP:eq:goal_functional}
    G(u^\star)
    \coloneqq
    \int_\Omega
    (g\, u^\star - \boldsymbol{g} \cdot \nabla u^\star)
    \d x.
\end{equation}
The goal-oriented adaptive finite element method
(GOAFEM) seeks for an efficient approximation 
of \(G(u^\star)\) using an (approximate) solution \(u_H \in \XX_H\)
to the discrete formulation~\eqref{BMP:eq:discrete}.

A duality approach from \cite{gs2002} employs
the solution \(z^\star \in \XX\)
to the so-called dual problem
\begin{equation}
    \label{BMP:eq:dual}
    a(v, z^\star)
    =
    G(v)
    \quad\text{for all }
    v \in \XX.
\end{equation}
While \(\eta_H(u_H^\star)\) from~\eqref{BMP:eq:estimator}
controls the energy error of the primal variable \(u_H^\star \approx u^\star\),
GOAFEM requires another a-posteriori
error estimator
for the
discrete dual problem with solution \(z^\star_H \in \XX_H\) to the problem
\begin{equation}
    \label{BMP:eq:discrete_dual}
    a(v_H, z^\star_H)
    =
    G(v_H)
    \quad\text{for all }
    v_H \in \XX_H.
\end{equation}
To this end, note that~\eqref{BMP:eq:dual} corresponds to the
elliptic PDE
\[ 
    -\div (\boldsymbol{A}\nabla z^\star)
    =
    g - \div \boldsymbol{g} 
    \text{ in } \Omega
    \quad\text{and}\quad
    z^\star
    =
    0
    \text{ on } \partial\Omega.
\]
    For \(z_H \in \XX_H\),
we thus define
\(
    \zeta_H(z_H)
    \coloneqq
    \big(
    \sum_{T \in \TT_{H}}
    \zeta_H(T, z_H)^2
    \big)^{1/2}
\)
by
\[
    \zeta_H(T, z_H)^2
    \coloneqq
    \vert T \vert^{2/d} \,
    \Vert g \, + \, \div(A \nabla z_H \, - \, \boldsymbol{g})
    \Vert_{L^2(T)}^2
    +
    \vert T \vert^{1/d}  \,
    \Vert \lbracket (A \nabla z_H  - \boldsymbol{g}) \cdot n
    \rbracket \Vert_{L^2(\partial T \cap \Omega)}^2
\]
and employ the same abbreviation \(\zeta_H(z_H)\) for the global values
analogously to the primal estimator with~\eqref{BMP:eq:estimator:b}.
It is well-known that \(\zeta_H\) satisfies
\eqref{BMP:axiom:stability}--\eqref{BMP:axiom:discrete_reliability}
with \(\eta_H\) replaced by \(\zeta_H\)
and all primal variables replaced by their dual
counterparts \cite{fpz2016}.
Hence, the dual estimator \(\zeta_H(z_H)\) allows for a
quasi-optimal adaptive resolution of the dual solution \(z^\star\).
Using the linearity of the goal functional \(G\) and the
dual problem~\eqref{BMP:eq:discrete_dual}, for any \(u_H, z_H \in \XX_H\),
recall from \cite{gs2002} that there holds
\begin{align*}
    G(u^\star) - G(u_H)
    &=
    G(u^\star - u_H)
    =
    a(u^\star - u_H, z^\star)
    =
    a(u^\star - u_H, z^\star - z_H)
    +
    a(u^\star - u_H, z_H)
    \\
    &=
    a(u^\star - u_H, z^\star - z_H)
    +
    \big[
        F(z_H)
        -
        a(u_H, z_H)
    \big].
\end{align*}
This motivates the definition
of the discrete goal functional
\begin{equation}
    \label{BMP:eq:discrete_goal_functional}
    G_H(u_H, z_H)
    \coloneqq
    G(u_H)
    +
    \big[
        F(z_H)
        -
        a(u_H, z_H)
    \big]
\end{equation}
and ensures the estimate
\[
    \vert
    G(u^\star)
    -
    G_H(u_H, z_H)
    \vert
    \leq
    \vert
    a(u^\star - u_H, z^\star - z_H)
    \vert
    \leq
    \vvvert u^\star - u_H \vvvert
    \,
    \vvvert z^\star - z_H \vvvert
\]
for the \(a (\cdot, \cdot)\)-induced energy norm.
Consequently,
for sufficiently good approximations
\(u_H\) to \(u^\star\) and \(z_H\) to \(z^\star\),
the discrete goal value \(G_H(u_H, z_H)\) converges to \(G(u^\star)\)
with the sum of the convergence rates
for the primal and the dual problem.
The reliability~\eqref{BMP:axiom:reliability}
of the primal and the dual estimator
result in
\begin{equation}
    \label{BMP:eq:goal_error_estimate}
    \vert
    G(u^\star)
    -
    G_H(u_H^\star, z_H^\star)
    \vert
    \lesssim
    \eta_H(u_H^\star)
    \,
    \zeta_H(z_H^\star).
\end{equation}
This reveals that GOAFEM needs to ensure
optimal convergence of the product of the primal and the
dual estimator.
To this end, \cite{fpz2016} introduces
the marking strategy in step~(III)
of Algorithm~\ref{BMP:algorithm:goafem}
below that avoids over-refinement for both primal and dual solutions.
Since the bilinear form \(a (\cdot, \cdot)\) is symmetric,
the iterative solver \(\Psi\) from Section~\ref{BMP:sec:solve}
may be also applied for the solution
of the discrete dual problem~\eqref{BMP:eq:discrete_dual}.

\subsection{Algorithm 2: Goal-oriented adaptive FEM}
\labeltext{2}{BMP:algorithm:goafem}

\noindent
{\bfseries Input:}
Initial mesh \(\TT_0\),
polynomial degree \(p \in \N\),
initial iterates
\(u_0^0 \coloneqq u_0^{\mm} \coloneqq z_0^0 \coloneqq z_0^{\mmu} \coloneqq 0 \),
marking parameter \(0 < \theta \le 1\),
solver parameter \(\lambdaalg > 0\).

\medskip
\noindent
{\bfseries Adaptive loop:}
For all \(\ell = 0, 1, 2, \dots\),
repeat the following steps {\rmfamily(I)--(IV)}:
\begin{enumerate}[leftmargin=2.5em]
    \item[\rmfamily(I)]
        {\ttfamily SOLVE \& ESTIMATE (PRIMAL).}
        For all \(m = 1, 2, 3, \dots\),
        repeat {\rmfamily(a)--(b)}:
        \begin{itemize}[leftmargin=2em]
            \item[\rmfamily(a)]
                Compute
                \(
                    u_\ell^m
                    \coloneqq
                    \Psi_\ell(F; u_\ell^{m-1})
                \)
                and the corresponding refinement indicators
                \\
                \(\eta_\ell(T; u_\ell^m)\)
                for all \(T \in \TT_\ell\).
            \item[\rmfamily(b)]
                Terminate \(m\)-loop and
                define \(\mm[\ell] \coloneqq m\)
                provided that
                \[
                    \vvvert u_\ell^m - u_\ell^{m-1} \vvvert
                    \le
                    \lambdaalg \,
                    \eta_\ell(u_\ell^m).
                \]
        \end{itemize}
    \item[\rmfamily(II)]
        {\ttfamily SOLVE \& ESTIMATE (DUAL).}
        For all \(\mu = 1, 2, 3, \dots\),
        repeat {\rmfamily(a)--(b)}:
        \begin{itemize}[leftmargin=2em]
            \item[\rmfamily(a)]
                Compute
                \(
                    z_\ell^\mu
                    \coloneqq
                    \Psi_\ell(G; z_\ell^{\mu-1})
                \)
                and the corresponding refinement indicators
                \\
                \(\zeta_\ell(T; z_\ell^\mu)\)
                for all \(T \in \TT_\ell\).
            \item[\rmfamily(b)]
                Terminate \(\mu\)-loop and
                define
                \(\mmu[\ell] \coloneqq \mu\)
                provided that
                \[
                    \vvvert z_\ell^\mu - z_\ell^{\mu-1} \vvvert
                    \le
                    \lambdaalg \,
                    \zeta_\ell(z_\ell^\mu).
                \]
        \end{itemize}
    \item[\rmfamily(III)]
        {\ttfamily COMBINED MARK.}
        Determine sets
        \(
            \overline{\MM}_\ell^u,
            \overline{\MM}_\ell^z
            \subseteq \TT_\ell
        \)
        of minimal cardinality that satisfy
        the Dörfler marking criteria
        \[
            \theta \, \eta_\ell(u_\ell^{\mm})^2
            \le
            \eta_\ell(\overline{\MM}_\ell^u, u_\ell^\mm)^2
            \quad\text{and}\quad
            \theta \, \zeta_\ell(z_\ell^\mmu)^2
            \le
            \zeta_\ell(\overline{\MM}_\ell^z, z_\ell^\mmu)^2.
        \]
        Following \cite{fpz2016},
        let the marked elements
        \(
            \MM_\ell
            \coloneqq
            \MM_\ell^u \cup \MM_\ell^z
        \)
        consist of subsets
        \(
            \MM_\ell^u
            \subseteq
            \overline{\MM}_\ell^u
        \)
        and
        \(
            \MM_\ell^z
            \subseteq
            \overline{\MM}_\ell^z
        \)
        of same cardinality
        \(
            \#\MM_\ell^u
            =
            \# \MM_\ell^z
            =
            \min\{
                \#\overline{\MM}_\ell^u,
                \#\overline{\MM}_\ell^z
            \}
        \).
    \item[\rmfamily(IV)]
            {\ttfamily REFINE.}
            Generate the new mesh
            \(
                \TT_{\ell+1}
                \coloneqq
                \texttt{refine}(\TT_\ell, \MM_\ell)
            \)
            by NVB and employ nested iteration
            \(u_{\ell+1}^0 \coloneqq u_\ell^{\mm}\) and
            \(z_{\ell+1}^0 \coloneqq z_\ell^{\mmu}\).
\end{enumerate}

\noindent
{\bfseries Output:}
Sequences of successively
refined triangulations \(\TT_\ell\),
discrete approximations
\(u_\ell^m\), \(z_\ell^\mu\),
and corresponding error estimators
\(\eta_\ell(u_\ell^k)\),
\(\zeta_\ell(z_\ell^\kappa)\).
\medskip

Note that, in practice, the primal~(I) and dual~(II)
iterations can be performed in parallel: The underlying matrix is
indeed the same, one just needs to accommodate the
implementation for two different right-hand sides.
Algorithm~\ref{BMP:algorithm:goafem} essentially goes 
back to~\cite{MS09} for the Poisson model problem, while the 
analysis for general second-order elliptic PDEs is given in~\cite{fpz2016}.
While the earlier work~\cite{MS09} employed the combined marking strategy 
with
\(
    \MM_\ell
    \in
    \{
        \overline{\MM}_\ell^u,
        \overline{\MM}_\ell^z
    \}
\)
such that
\(
    \# \MM_\ell
    =
    \min\{
        \#\overline{\MM}_\ell^u,
        \#\overline{\MM}_\ell^z
    \}
\),
it empirically turned out that the marking strategy \textbf{(III)} 
from \cite{fpz2016}
is computationally more efficient.
Note, however, that~\cite{MS09,fpz2016} employ the exact computation of the 
primal and the dual solutions, while the extension to inexact solutions 
is done in~\cite{bgip2023} for symmetric PDEs.

The simultaneous consideration of primal and dual solution
requires an adaptation of the countably infinite index
set which reads, for the function \(v \in \{u,z\}\),
\[
    \QQ^v
    \coloneqq
    \{
        (\ell, k) \in \N_0^2 \;:\;
        v_\ell^k \text{ is defined in
        Algorithm~\ref{BMP:algorithm:goafem}}
    \}.
\]
The following definitions of the final indices
coincide with those defined in
Algorithm~\ref{BMP:algorithm:goafem}
\begin{align*}
    \elll
    &\coloneqq
    \sup\{
        \ell \in \N_0 \;:\;
        (\ell, 0) \in \QQ^u
        \text{ or }
        (\ell, 0) \in \QQ^z
    \}
    \in \N_0 \cup \{\infty\},
    \\
    \mm[\ell]
    &\coloneqq
    \sup\{
        m \in \N \;:\;
        (\ell, m) \in \QQ^u
    \},
    \quad
    \mmu[\ell]
    \coloneqq
    \sup\{
        \mu \in \N \;:\;
        (\ell, \mu) \in \QQ^z
    \}.
\end{align*}
In addition to the quasi-error \(\Eta_\ell^k\) and \(\Eta_\ell\) 
from~\eqref{BMP:eq:quasi-error}--\eqref{BMP:eq:quasi-error-final-iterates} for 
the primal problem, the corresponding quasi-errors
for the dual problem read
\begin{align*}
    \textrm{Z}_\ell^{k}
    &\coloneqq
    \vvvert z_\ell^\star - z_\ell^{k} \vvvert
    +
    \zeta_\ell(z_\ell^{k})
    \quad\text{for all }
    (\ell, {k}) \in \QQ^z,
    \\
    \textrm{Z}_\ell
    &\coloneqq
    \vvvert z_\ell^\star - z_\ell^{\kk} \vvvert
    +
    \gamma\, \zeta_\ell(z_\ell^{\kk})
    \quad\text{for all }
    (\ell, {\kk}) \in \QQ^z.
\end{align*}
The quasi-errors naturally extend to the
full index set
\(
    (\ell, k)
    \in
    \QQ
    \coloneqq
    \QQ^u
    \cup
    \QQ^z
\)
with the maximum
\(
    \kk[\ell] \coloneqq
    \max\{\mm[\ell], \mmu[\ell]\}
\) via
\[
    \Eta_\ell^k
    \coloneqq
    \Eta_\ell^\mm
    \text{ if }
    (\ell, k) \not\in \QQ^u
    \quad\text{and}\quad
    \Zeta_\ell^k
    \coloneqq
    \Zeta_\ell^\mmu
    \text{ if }
    (\ell, k) \not\in \QQ^z.
\]

Arguing as for Lemma~\ref{BMP:prop:overall-error-estimator}, one obtains 
the following a-posteriori control on the goal-error, which can be 
used to terminate Algorithm~\ref{BMP:algorithm:goafem}.
\begin{lemma}
    For \(\ell \in \N_0\) and \(m, \mu \ge 1\), let 
    \(u_\ell^m, z_\ell^\mu \in \XX_\ell\) be computed by the adaptive loop of 
    Algorithm~\ref{BMP:algorithm:goafem}. Then, there holds 
    \[
        \vert
        G(u^\star)
        -
        G_H(u_\ell^m, z_\ell^\mu)
        \vert
        \lesssim 
        \bigl[ 
            \vvvert u_\ell^m - u_\ell^{m-1} \vvvert 
            + \eta_\ell(u_\ell^m) 
        \bigr] 
        \,
        \bigl[ 
            \vvvert z_\ell^\mu - z_\ell^{\mu-1} \vvvert 
            + \zeta_\ell(z_\ell^\mu)
        \bigr]
    \]
    and, for \( m = \mm[\ell]\), \( \mu = \mmu[\ell] \),
    \begin{align}
        \label{BMP:eq:goal_error_aposteriori_estimate_final}
        \vert
        G(u^\star)
        -
        G_H(u_\ell^\mm, z_\ell^\mmu)
        \vert
        \lesssim 
        \eta_\ell(u_\ell^\mm)  \zeta_\ell(z_\ell^\mmu).
        \hspace{3.5em}\qed\hspace{-3.5em}
    \end{align}
\end{lemma}

\subsection{Convergence and complexity}

Algorithm~\ref{BMP:algorithm:goafem}
satisfies the following convergence
results from \cite{bgip2023}
generalizing the statements from
Subsection~\ref{BMP:sec:convergence}
to goal-oriented adaptivity.
The proof of full R-linear convergence follows the summability-based 
proof given for Theorem~\ref{BMP:thm:full-linear-convergence} above. 
Additional difficulties arise from the inherent nonlinear product structure 
of \(\Eta_\ell^k \Zeta_\ell^k\), when compared to the linear quasi-error 
\(\Eta_\ell^k \) for standard AFEM.
\begin{theorem}[parameter-robust full R-linear convergence {\cite[Thm.~6]{bgip2023}}]
    \label{BMP:thm:full-linear-convergence-goafem}
    Let \(0 < \theta \leq 1\)
    and \(\lambdaalg > 0\).
    Then, Algorithm~\ref{BMP:algorithm:goafem}
    guarantees the existence of constants
    \(\Clin > 0\) and
    \(0 < \qlin < 1\) such that
    \[
        \Eta_\ell^k \Zeta_\ell^k
        \leq
        \Clin\,
        \qlin^{\vert \ell,k \vert - \vert
        \ell', k' \vert}\,
        \Eta_{\ell'}^{k'} \Zeta_{\ell'}^{k'}
        \quad\text{for all }
        (\ell, k), (\ell^\prime, k^\prime) \in \QQ
        \text{ with }
        \vert \ell, k \vert > \vert \ell^\prime, k^\prime \vert.
    \]
    In particular, this yields parameter-robust convergence 
    \[
        \vert
        G(u^\star)
        -
        G_\ell(u_\ell^k, z_\ell^k)
        \vert
        \lesssim 
        \Eta_\ell^k \Zeta_\ell^k
        \leq
        \Clin\,
        \qlin^{\vert \ell,k \vert}\,
        \Eta_{0}^{0} \Zeta_{0}^{0} 
        \rightarrow 0 
        \text{ as } \vert \ell,k \vert \rightarrow \infty.
        \qed
    \]
\end{theorem}
The proof of optimal complexity follows as for standard AFEM. The only 
new ingredient is the comparison lemma 
(see Lemma~\ref{BMP:lem:optimal_refinement} 
for standard AFEM), which needs to be adapted to the product structure of 
GOAFEM.
\begin{theorem}[optimal complexity {\cite[Thm.~8]{bgip2023}}]
    \label{BMP:thm:optimal-complexity-goafem}
    There exist upper bounds
    \(0 < \theta^\star \leq 1\)
    and \(\lambdaalg^\star > 0\)
    such that, for sufficiently small
    \(0 < \theta < \theta^\star\)
    and \(0 < \lambdaalg < \lambdaalg^\star\),
    the following holds:
    Algorithm~\ref{BMP:algorithm:goafem}
    guarantees that
    \begin{align*}
        \sup_{(\ell, k) \in \QQ}
        \cost(\ell, k)^{s+t}\,
        \Eta_\ell^k\,
        \Zeta_\ell^k
        \lesssim
        \max\{
            \Vert u^\star \Vert_{\A_s}\,
            \Vert z^\star \Vert_{\A_t},\:
            \Eta_0^0\,
            \Zeta_0^0
        \}
        \text{ for all }
        s, t > 0.
        \qed
    \end{align*}
\end{theorem}

\subsection{Numerical example}
\label{BMP:sec:goafem:experiment}

For the empirical investigation of
Algorithm~\ref{BMP:algorithm:goafem},
consider the Z-shaped domain
\(
    \Omega
    \coloneqq
    (-1, 1)^2 \setminus
    \operatorname{conv}\{(0,0), (-1,0), (-1-1)\}
    \subset \R^2
\)
with the Dirichlet boundary
\(
    \Gamma_{\textup{D}}
    \coloneqq
    \operatorname{conv}\{(-1,0), (0,0)\}
    \cup
    \operatorname{conv}\{(-1,-1), (0,0)\}
\)
and the Neumann boundary
\(
    \Gamma_{\textup{N}}
    \coloneqq
    \partial\Omega \setminus \Gamma_{\textup{D}}
\).
The benchmark problem seeks \(u^\star \in H^1(\Omega)\)
such that
\[
    - \Delta u^\star
    =
    1
    \text{ in } \Omega
    \quad\text{with}\quad
    u^\star = 0
    \text{ on } \Gamma_{\textup{D}}
    \text{ and }
    \partial u / \partial n = 0
    \text{ on } \Gamma_{\textup{N}}.
\]
Given the subdomain
\(
    S \coloneqq (-0.5, 0.5)^2 \cap \Omega
\),
the goal functional is determined
by \(g \coloneqq 0\)
and \(\boldsymbol{g} \coloneqq \chi_S (1,1)^\top\)
and reads
\begin{equation}
    \label{BMP:eq:goal_value_extrapolated}
    G(u^\star)
    =
    \int_S (\partial_1 u^\star + \partial_2 u^\star)
    \d{x}
    \approx
    1.015559272415834.
\end{equation}
The exact goal value~\eqref{BMP:eq:goal_value_extrapolated} 
is approximated by a reference solution with sextic polynomial functions
with adaptive mesh refinement until \(\dim(\XX_\ell)\) exceeds
\(5\times 10^6\).

\begin{figure}
    \centering

    \subfloat[Initial mesh \(\TT_0\).]{%
        \label{BMP:fig:GOAFEM:initial_mesh}
        \begin{tikzpicture}
    \begin{axis}[%
        axis equal image,%
        width=4.6cm,%
        xmin=-1.15, xmax=1.15,%
        ymin=-1.15, ymax=1.15,%
        font=\footnotesize%
    ]
        \addplot graphics [xmin=-1, xmax=1, ymin=-1, ymax=1]
        {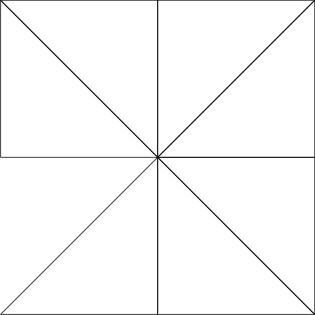};
    \end{axis}
\end{tikzpicture}
    }
    \hfil
    \subfloat[Adaptively refined mesh on level \(\ell = 11\) with \(1\,015\) triangles.]{%
        \label{BMP:fig:GOAFEM:mesh}
        \begin{tikzpicture}
    \begin{axis}[%
        axis equal image,%
        width=4.6cm,%
        xmin=-1.15, xmax=1.15,%
        ymin=-1.15, ymax=1.15,%
        font=\footnotesize%
    ]
        \addplot graphics [xmin=-1, xmax=1, ymin=-1, ymax=1]
        {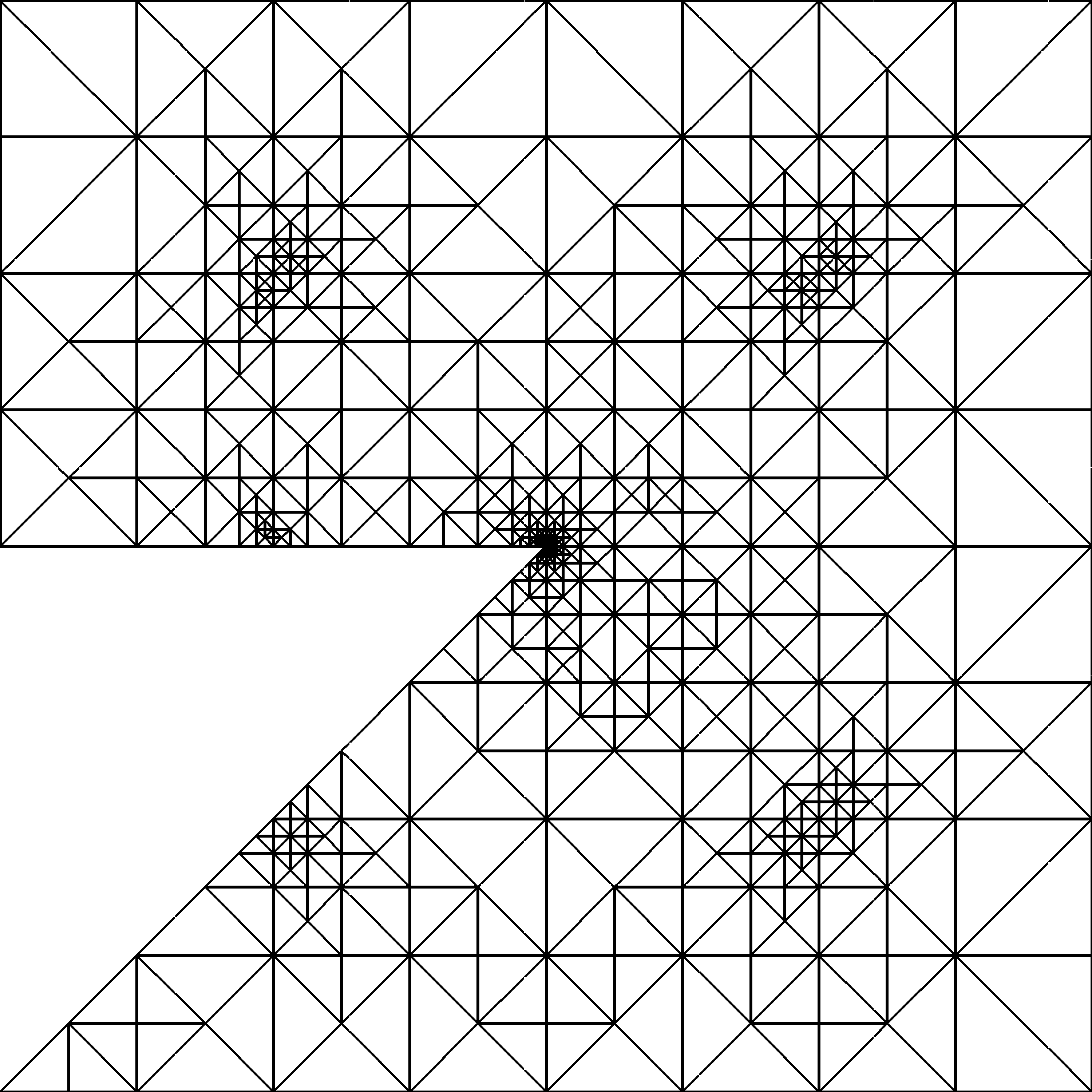};
    \end{axis}
\end{tikzpicture}
    }
    \medskip

    \subfloat[Discrete primal solution \(u_{11}^{\kk} \in \XX_{11}\).]{%
        \label{BMP:fig:GOAFEM:solution_primal}
        \begin{tikzpicture}
    \pgfplotsset{/pgf/number format/fixed}
    \begin{axis}[%
        width=4.2cm,%
        xmin=-1.1, xmax=1.1,%
        ymin=-1.1, ymax=1.1,%
        zmin=-0.1, zmax=1.17,%
        font=\footnotesize,%
    ]
        \addplot3 graphics [%
            points={%
                (-1,1,0.671699) => (1.4,339.7-143.0)
                (-1,0,0) => (120.5,339.7-318.1)
                (1,-1,0.964961) => (480.2,339.7-75.6)
                (1,1,1.11083) => (241.1,339.7-0.5)
            }%
            ]
            {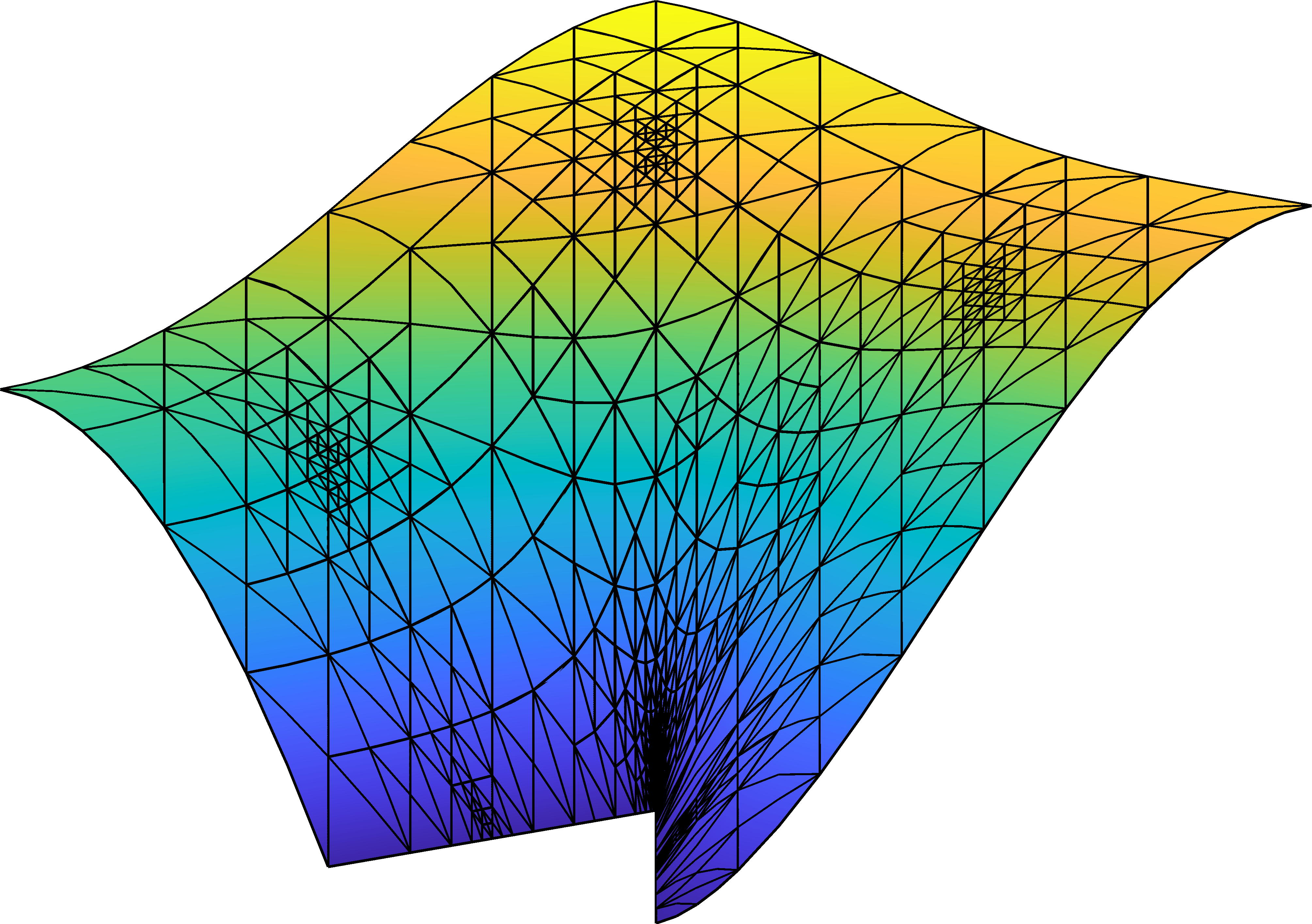};
    \end{axis}
\end{tikzpicture}
    }
    \hfil
    \subfloat[Discrete dual solution \(z_{11}^{\kk} \in \XX_{11}\).]{%
        \label{BMP:fig:GOAFEM:solution_dual}
        \begin{tikzpicture}
    \pgfplotsset{/pgf/number format/fixed}
    \begin{axis}[%
        width=4.0cm,%
        xmin=-1.1, xmax=1.1,%
        ymin=-1.1, ymax=1.1,%
        zmin=-0.14, zmax=0.775,%
        font=\footnotesize,%
    ]
        \addplot3 graphics [%
            points={%
                (-1,1,0.23535) => (1.7,305.9-178.3)
                (-1,0,0) => (120.6,305.9-280.2)
                (1,-1,0.240754) => (480.4,305.9-176.3)
                (0.5,0.5,0.688078) => (241.0,305.9-1.0)
            }%
            ]
            {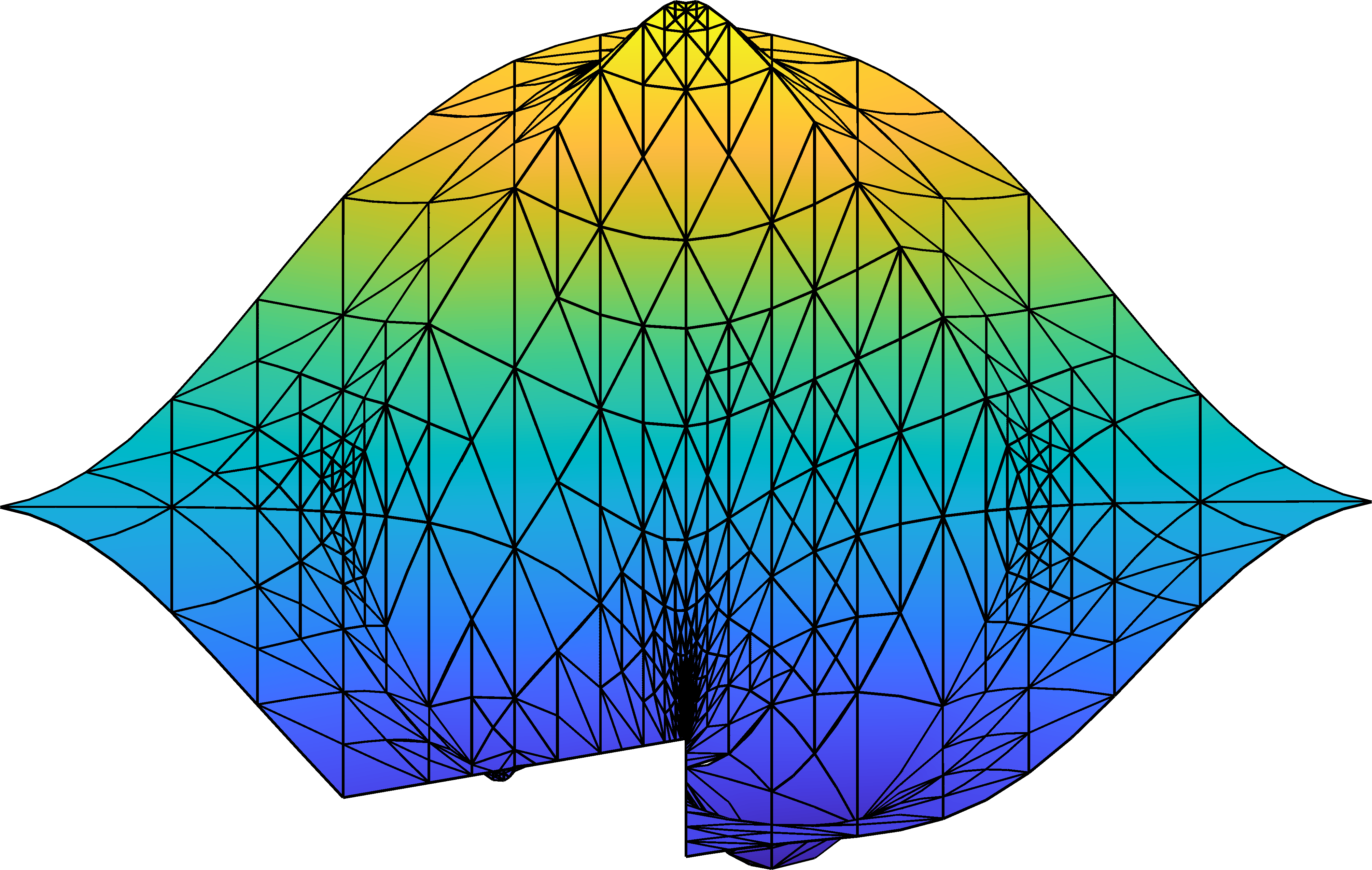};
    \end{axis}
\end{tikzpicture}
    }
    \caption{%
        Mesh and corresponding discrete primal and dual solution
        for the benchmark problem from
        Subsection~\ref{BMP:sec:goafem:experiment}.
        The results are generated by Algorithm~\ref{BMP:algorithm:goafem} with
        polynomial degree \(p = 2\), bulk parameter \(\theta = 0.3\),
        and stopping parameter \(\lambdaalg = 0.7\).
    }
    \label{BMP:fig:goafem:mesh_solution}
\end{figure}

\begin{figure}
    \centering
    \input{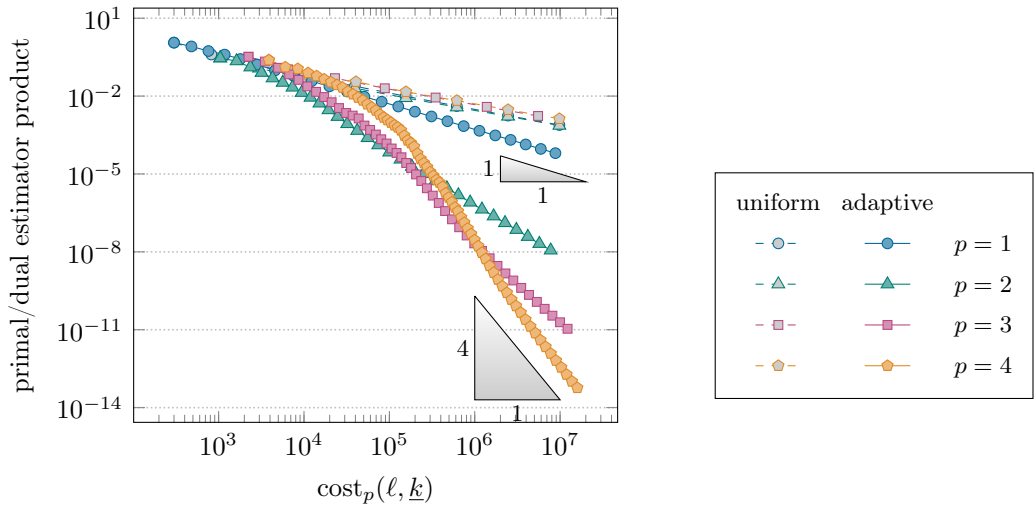}
    \caption{
        Convergence plot of the
        adaptive Algorithm~\ref{BMP:algorithm:goafem}
        to solve the benchmark problem from
        Subsection~\ref{BMP:sec:goafem:experiment}
        for various polynomial degrees.
        The adaptivity parameters read \(\theta = 0.3\)
        and \(\lambdaalg = 0.7\).
        All graphs display values of
        the residual-based error estimator product
        \(\eta_\ell (u_\ell^\mm )\, \zeta_\ell  (z_\ell^\mmu )\)
        from~\eqref{BMP:eq:goal_error_aposteriori_estimate_final}.
    }
    \label{BMP:fig:goafem:convergence}
\end{figure}

\begin{figure}
    \centering
    \input{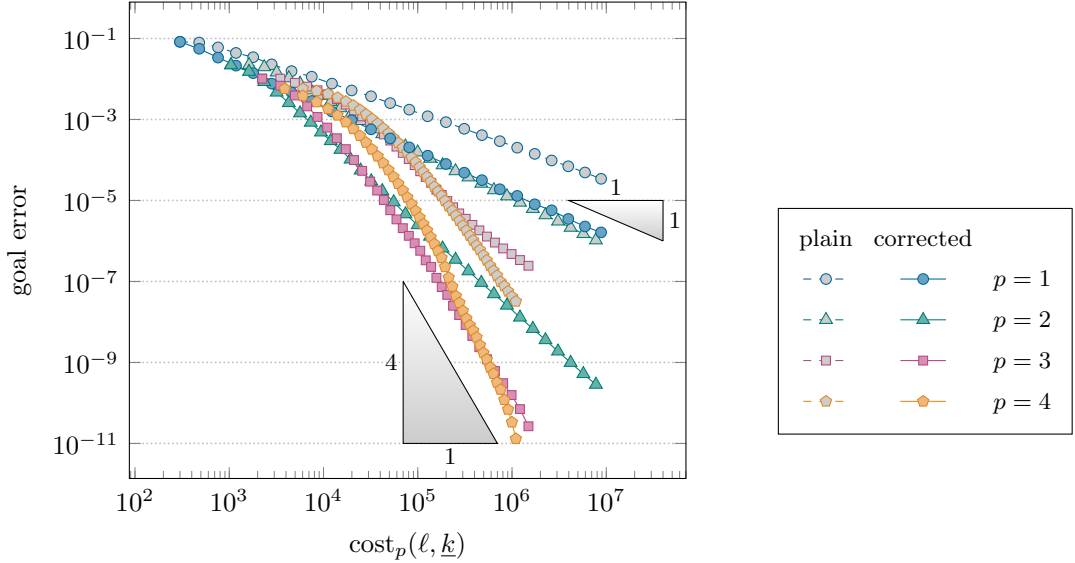}
    \caption{
        Convergence plot of the goal error
        for the adaptive Algorithm~\ref{BMP:algorithm:goafem}
        to solve the benchmark problem from
        Subsection~\ref{BMP:sec:goafem:experiment}
        for various polynomial degrees.
        The adaptivity parameters read \(\theta = 0.3\)
        and \(\lambdaalg = 0.7\).
        The graphs display the error of the
        discrete goal value \(G(u_\ell^\mm )\)
        from~\eqref{BMP:eq:goal_functional} using only the primal problem (plain)
        and the discrete
        goal value \(G_H(u_\ell^\mm , z_\ell^\mmu )\)
        from~\eqref{BMP:eq:discrete_goal_functional}
        using the duality technique to double the convergence rates (corrected).
        The error is computed with respect to the
        reference goal value
        from~\eqref{BMP:eq:goal_value_extrapolated}.
        The reference value only allowed reasonable
        results up to error values of about \(10^{-11}\) and
        thus the remaining results of the graphs are omitted.
    }
    \label{BMP:fig:goafem:goal_error}
\end{figure}

Some discrete primal and dual solutions are depicted
in Figure~\ref{BMP:fig:goafem:mesh_solution}.
The corresponding adaptively generated mesh exhibits
increased refinement at the reentrant corner
(singularity of \(u^\star\))
as well as at the vertices of the support \(S\) of
\(\boldsymbol{g}\) (singularities of \(z^\star\)).
This adaptive refinement ensures the (doubled) optimal
convergence rates for the primal and dual estimator product
in Figure~\ref{BMP:fig:goafem:convergence}.
The approximated goal error with respect to the extrapolated
goal value from~\eqref{BMP:eq:goal_value_extrapolated}
in Figure~\ref{BMP:fig:goafem:goal_error}
confirms these convergence rates.

%
%

\section{Extension 2: Non-symmetric problems}
\label{BMP:sec:aisfem}

While iterative algebraic solvers for non-symmetric problems are
widely used in practice (e.g., GMRES, BiCGStab),
obtaining a uniform contraction property~\eqref{BMP:eq:algebra}
in the energy norm remains an open question.
As a remedy, and inspired by the state-of-the-art proof of the famous 
Lax--Milgram lemma, \cite{aisfem,aisfem-corr} introduced an additional
Zarantonello iteration for the symmetrization of the discrete equation.
This allows for the application of the 
previously used optimal iterative
solvers within
Algorithm~\ref{BMP:algorithm:afem} 
for symmetric linear systems to compute the inexact
Zarantonello update.
This section addresses the generalization
of the adaptive algorithm to non-symmetric problems
by using such a two-fold nested iteration
with a symmetrization step.

\subsection{General second-order linear elliptic model problem}

The (non-symmetric) second-order linear elliptic problem
seeks a solution \(u^\star \in \XX\) with
\begin{equation}
    \label{BMP:eq:general-second-order}
    -\div(\boldsymbol{A} \nabla u^\star)
    + \boldsymbol{b} \cdot \nabla u^\star
    + cu^\star
    =
    f - \div \boldsymbol{f}
    \text{ in } \Omega
    \quad\text{and}\quad
    u^\star = 0
    \text{ on } \partial\Omega,
\end{equation}
where
\(\boldsymbol{b} \in \bigl[L^\infty(\Omega)\bigr]^{d}\) is 
the convection coefficient
and \(c \in L^\infty(\Omega)\) the reaction coefficient, 
in addition to the diffusion matrix \(\boldsymbol{A}\) and
the right-hand sides \(f\) and \(\boldsymbol{f}\) 
from Subsection~\ref{BMP:sec:model_problem}.
Using the bilinear form \(a(\cdot,\cdot)\) from~\eqref{BMP:eq:weak},
the weak formulation of~\eqref{BMP:eq:general-second-order}
seeks \(u^\star \in \XX\) such that
\begin{equation}
    \label{BMP:eq:weak_nonsymmetric}
    b(u^\star, v)
    \coloneqq
    a(u^\star, v)
    +
    \int_\Omega 
    (\boldsymbol{b} \cdot \nabla u^\star + c u^\star) \, v
    \d x
    =
    F(v)
    \quad\text{for all }
    v \in \XX.
\end{equation}
We suppose that \(b(\cdot,\cdot)\) is elliptic
with respect to the \(a(\cdot,\cdot)\)-induced norm \(\vvvert \cdot \vvvert\)
in the sense of 
\begin{equation}
    \label{BMP:eq:AISFEM:ellipticity}
    b(v,v)
    \gtrsim \vvvert v \vvvert^2 
    \quad \text{for all } v \in \XX.
\end{equation}
For instance, this is satisfied if \(-\frac{1}{2} \div \boldsymbol{b} +c \ge 0\). 
Under the assumption~\eqref{BMP:eq:AISFEM:ellipticity},
the Lax--Milgram lemma applies and indeed proves
existence and uniqueness of the weak solution \(u^\star \in \XX \) 
to~\eqref{BMP:eq:weak_nonsymmetric} and likewise that of the discrete 
solution \(u^\star_H \in \XX_H \) to 
\[
    b(u^\star_H, v_H)
    =
    F(v_H)
    \quad\text{for all }
    v_H \in \XX_H.
\] 
Moreover, the corresponding residual-based error estimator reads
\begin{equation}
    \label{BMP:eq:aisfem:estimator}
    \begin{split}
        \eta_H(T, v_H)^2
        &\coloneqq
        \vert T \vert^{2/d} 
        \,
        \Vert
        f + \div(\boldsymbol{A} \nabla u^\star - \boldsymbol{f})
        - \boldsymbol{b} \cdot \nabla u^\star
        - c u^\star
        \Vert_{L^2(T)}^2
        \\
        &\phantom{{}\coloneqq{}}
        + \, \vert T \vert^{1/d}  \,
        \Vert
        \lbracket (\boldsymbol{A} \nabla u^\star - \boldsymbol{f})
        \cdot n \rbracket
        \Vert_{L^2(\partial T \cap \Omega)}^2.
    \end{split}
\end{equation}
It is well-known from, e.g., \cite[Sect.~6]{cfpp2014} that \(\eta_H\)
satisfies the axioms of adaptivity
\eqref{BMP:axiom:stability}--\eqref{BMP:axiom:discrete_reliability}.

Standard iterative solvers for non-symmetric linear systems
lack the contraction property~\eqref{BMP:eq:algebra}
in the energy norm.
This prevents the straight-forward application
of the full R-linear convergence proof from
Subsection~\ref{BMP:sec:full_linear_convergence}
to problem~\eqref{BMP:eq:general-second-order}
and motivated the use of
the additional Zarantonello iteration
\cite{Zarantonello1960}
for the symmetrization of the discrete equation
in \cite{aisfem}.
For a given damping parameter \(\delta > 0\)
and the current iterate \(u_H \in \XX_H\),
the Zarantonello mapping
\(\Phi_H(\delta; \cdot): \XX_H \to \XX_H\)
is defined by
\begin{equation}
    \label{BMP:eq:Zarantonello}
    a(\Phi_H(\delta; u_H), v_H)
    =
    a(u_H, v_H)
    +
    \delta
    [F(v_H) - b(u_H, v_H)]
    \text{ for all }
    v_H \in \XX_H.
\end{equation}
Since \(a (\cdot,\cdot)\) is a scalar product,
the Riesz--Fischer theorem guarantees well-defined\-ness
of \(\Phi_H(\delta; u_H) \in \XX_H \).
For sufficiently small \(\delta > 0\),
the mapping \(\Phi_H(\delta; \cdot)\)
is a contraction \cite{hw2020,hw2020:conv},
i.e., there exists \(0<\qsym^\star <1\) such that
\begin{equation}
    \label{BMP:eq:Zarantonello_contraction}
    \vvvert u_H^\star - \Phi_H(\delta; u_H) \vvvert
    \leq
    \qsym^\star \,
    \vvvert u_H^\star - u_H \vvvert
    \quad
    \text{for all }
    u_H \in \XX_H.
\end{equation}
As in Subsection~\ref{BMP:sec:solve},
the linear SPD system~\eqref{BMP:eq:Zarantonello}
can be treated using the algebraic solver
\(\Psi_H: \XX^* \times \XX_H \to \XX_H\)
from Subsection~\ref{BMP:sec:solve}
with the contraction property~\eqref{BMP:eq:algebra}.
In case of a complicated right-hand side as
in~\eqref{BMP:eq:Zarantonello},
the first argument of \(\Psi_H\)
is replaced by the exact Riesz representative
\(\Phi_H(\delta; u_H)\) of the right-hand side
to shorten the notation.
Importantly, the solution 
\(\Phi_H(\delta; u_H)\) of~\eqref{BMP:eq:Zarantonello}
will never be computed in practice, but only approximated via 
this iterative 
solver. 

Overall, a nested solver iteration is necessary
and thus a triple index \( (\ell,k,j) \) 
will be used in the adaptive algorithm: The index \(\ell \in \N_0\) 
denotes the mesh-level, \(k \in \N_0\) denotes the Zarantonello step, and 
\(j \in \N_0\) denotes the step of the algebraic solver, i.e., 
\[
    u_\ell^{k,j} = \Psi_\ell(u_\ell^{k,\star},u_\ell^{k,j-1}) 
    \approx u_\ell^{k,\star} 
    \coloneqq \Phi_{\ell} (\delta; u_\ell^{k-1,\jj}) 
    \approx u_\ell^\star,
\] 
where neither \(u_\ell^{k,\star}\) nor 
\(u_\ell^\star\) are computed, but only \(u_\ell^{k,j}\).

\subsection{Algorithm 3: Adaptive iteratively symmetrized FEM}
\label{BMP:sec:aisfem:algorithm}
\labeltext{3}{BMP:algorithm:aisfem}

\noindent
{\bfseries Input:}
Initial mesh \(\TT_0\),
polynomial degree \(p \in \N\),
initial iterate $u_0^{0,0} \coloneqq u_0^{0,\underline{j}} \coloneqq 0$,
the Zarantonello damping parameter $\delta > 0$,
marking parameter \(0 < \theta \le 1\),
stopping parameters $\lambdasym, \lambdaalg > 0$.

\medskip
\noindent
{\bfseries Adaptive loop:}
For all \(\ell = 0, 1, 2, \dots\),
repeat the steps~(I)--(III):
\begin{enumerate}[leftmargin=2.5em]
    \item[\rmfamily(I)]
        {\ttfamily SOLVE \& ESTIMATE.}
        For all \(k = 1, 2, 3, \dots\),
        repeat the steps (a)--(b):
        \begin{itemize}[leftmargin=1.7em]
            \item[\rmfamily(a)]
                \textbf{Algebraic solver loop:}
                For all \(j = 1, 2, 3, \dots\),
                repeat the steps~(i)--(ii):
                \begin{itemize}[leftmargin=2em]
                    \item[\rmfamily(i)]
                        Compute
                        \(
                            u_\ell^{k,j}
                            \coloneqq
                            \Psi_\ell(u_\ell^{k,\star}, u_\ell^{k,j-1})
                        \)
                        (for the approximation of the theoretical
                        quantity
                        $u_\ell^{k, \star}
                        \coloneqq
                        \Phi_\ell(\delta;u_\ell^{k-1,\underline j})
                        \in \XX_\ell$ solving~\eqref{BMP:eq:Zarantonello} 
                        with \(u_H = u_\ell^{k-1,\underline j}\))
                        and refinement indicators
                        \(\eta_\ell(T, u_\ell^{k,j})\) for all $T \in \TT_\ell$.
                    \item[\rmfamily(ii)]
                        Terminate \(j\)-loop 
                        with \(\jj[\ell,k] \coloneqq j\)
                        and employ nested iteration \(u_\ell^{k, 0} \coloneqq u_\ell^{k-1,\jj}\)
                        provided that
                        \begin{align}\label{BMP:eq:aisfem:algstop}
                            \! \! \! \! \!  \! \! \! \!
                            \vvvert u_\ell^{k,j} - u_\ell^{k,j-1} \vvvert
                            \le
                            \lambdaalg \, \big[ \lambdasym
                            \eta_\ell(u_\ell^{k,j}) +
                            \vvvert u_\ell^{k,j} 
                            - u_\ell^{k-1,\underline{j}} \vvvert \big].
                            \! \!
                        \end{align}
                \end{itemize}
            \item[\rmfamily(b)]
                Terminate \(k\)-loop and
                define \(\kk[\ell] \coloneqq k\)
                provided that
                \begin{align}\label{BMP:eq:aisfem:symstop}
                    \vvvert u_\ell^{k, \jj} - u_\ell^{k-1, \jj} \vvvert
                    \le
                    \lambdasym \, \eta_\ell(u_\ell^{k, \jj}).
                \end{align}
        \end{itemize}
    \item[\rmfamily(II)]
        {\ttfamily MARK.}
        Determine a set \(\MM_\ell \subseteq \TT_\ell\)
        of minimal cardinality that satisfies
        the Dörfler marking criterion
        \[
            \theta \, \eta_\ell(u_\ell^{\kk,\jj})^2
            \le
            \eta_\ell(\MM_\ell, u_\ell^{\kk,\jj})^2.
        \]
        \item[\rmfamily(III)]
            {\ttfamily REFINE.}
            Generate the new mesh
            \(
                \TT_{\ell+1}
                \coloneqq
                \texttt{refine}(\TT_\ell, \MM_\ell)
            \)
            by NVB and define
            \(
                u_{\ell+1}^{0,0}
                \coloneqq
                u_{\ell+1}^{0,\jj}
                \coloneqq
                u_\ell^{\kk,\jj}
            \)
            (nested iteration).
\end{enumerate}

\noindent
{\bfseries Output:}
Sequences of successively
refined triangulations \(\TT_\ell\),
discrete approximations \(u_\ell^{k,j}\),
and corresponding error estimators
\(\eta_\ell(u_\ell^{k,j})\).
\medskip

\begin{remark}
    Let us comment on the stopping 
    criteria~\eqref{BMP:eq:aisfem:algstop}--\eqref{BMP:eq:aisfem:symstop} of 
    Algorithm~\ref{BMP:algorithm:aisfem}. 
    \begin{enumerate}[label={\normalfont (\roman*)}]
        \item Since the algebraic solver is contractive, the term 
        \(\vvvert u_\ell^{k, j} - u_\ell^{k, j-1} \vvvert\) 
        in~\eqref{BMP:eq:aisfem:algstop} provides computable a-posteriori 
        error control on the algebraic error 
        \(\vvvert u_\ell^{k, \star} - u_\ell^{k, j} \vvvert\). With the interpretation that 
        \(\eta_\ell(u_\ell^{k,j}) \approx \eta_\ell(u_\ell^{\star})\) 
        measures the discretization error and 
        \(\vvvert u_\ell^{k, j} - u_\ell^{k-1, \jj} \vvvert 
        \approx  \vvvert u_\ell^{k, \star} - u_\ell^{k, j} \vvvert \) 
        measures the symmetrization error, the stopping 
        criterion~\eqref{BMP:eq:aisfem:algstop} can be interpreted as numerical equilibration of the algebraic error with the sum of 
        discretization and symmetrization error.
        \item With the same understanding, the stopping criterion~\eqref{BMP:eq:aisfem:symstop} for the 
        symmetrization can be interpreted as numerical equilibration of  
        symmetrization and discretization error.
        \item These interpretations can be made rigorous with the help of 
        Lemma~\ref{BMP:lem:contraction_aisfem} below. Moreover, 
        Lemma~\ref{BMP:prop:aposteriori-jk-aisfem} below proves that 
        the full error \(\vvvert u^{\star} - u_\ell^{k, j} \vvvert\) 
        is indeed controlled by the sum of the terms 
        in~\eqref{BMP:eq:aisfem:algstop}. 
        An additional consequence 
        of~\eqref{BMP:eq:aisfem:algstop}--\eqref{BMP:eq:aisfem:symstop} is 
        that \(\eta_\ell(u_\ell^{\kk,\jj}) \) provides a computable upper 
        bound for the error \(\vvvert u^{\star} - u_\ell^{\kk,\jj} \vvvert\).
        \item The innermost \(j\)-loop of the algebraic solver indeed turns out 
        to be finite; see Lemma~\ref{BMP:lemma:unif-steps-aisfem} 
        below. The \(k\)-loop of the symmetrization can formally fail to 
        terminate. In this case, 
        Theorem~\ref{BMP:thm:full-linear-convergence-aisfem} below proves that 
        \(u^\star = u_\ell^\star\) and 
        \( \eta_\ell(u_\ell^{k,\jj}) \rightarrow  
        \eta_\ell(u_\ell^{\star}) = 0\) as \(k \rightarrow \infty\).
    \end{enumerate}
\end{remark}

Analogous notation to Subsection~\ref{BMP:sec:quasi-errors}
applies for the index set
\[
    \QQ
    \coloneqq
    \{
        (\ell, k, j) \in \N_0^3 
        \colon
        u_\ell^{k, j}
        \text{ is defined in Algorithm~\ref{BMP:algorithm:aisfem}}
    \},
\]
the lexicographic ordering on \(\QQ\)
\[
    (\ell', k', j') \le (\ell, k, j)
    :\Longleftrightarrow
    u_{\ell'}^{k', j'} \text{ is defined no later than }
    u_\ell^{k,j}
    \text{ in Algorithm~\ref{BMP:algorithm:aisfem}},
\]
the stopping indices
\begin{align*}
    \underline{\ell}
    &\coloneqq
    \sup \{ \ell \in \N_0 \colon (\ell,0,0) \in \QQ \}
    \in \N_0 \cup \{\infty\},
    \\
    \underline{k}[\ell]
    &\coloneqq
    \sup \{ k \in \N_0 \colon (\ell,k,0) \in \QQ \}
    \in \N_0 \cup \{\infty\}
    \text{ for } (\ell,0,0) \in \QQ,
    \\
    \underline{j}[\ell,k]
    &\coloneqq
    \sup \{ j \in \N_0 \colon (\ell,k,j) \in \QQ \}
    \in \N_0 \cup \{\infty\}
    \text{ for } (\ell,k,0) \in \QQ,
\end{align*}
the total step counter
\[
    \vert \ell, k,j \vert
    \coloneqq
    \#\{
        (\ell', k',j') \in \QQ
        \,:\,
        (\ell', k',j') \leq (\ell, k,j)
    \}
    \in
    \N_0
    \text{ for all }
    (\ell, k,j) \in \QQ,
\]
and the computational cost, for all \((\ell, k, j) \in \QQ\)
\[
    \cost(\ell, k, j)
    \coloneqq
    \sum_{\substack{
        (\ell^\prime, k^\prime, j^\prime) \in \QQ
        \\
        \vert \ell^\prime, k^\prime, j^\prime \vert
        \le \vert \ell, k, j\vert
    }}
    \# \TT_\ell
    \quad\text{and}\quad
    \cost_p(\ell, k, j)
    \coloneqq
    \sum_{\substack{
            (\ell^\prime, k^\prime, j^\prime) \in \QQ
            \\
            \vert \ell^\prime, k^\prime, j^\prime \vert
            \le \vert \ell, k, j\vert
    }}
    \dim(\XX_\ell).
\]
Moreover, the quasi-error is modified to include the additional 
error related to the symmetrization: Indeed, the first term is the 
error between exact discrete solution and its inexactly symmetrized 
approximation
\[
    \Eta_\ell^{k,j}
    \coloneqq
    \vvvert u_\ell^\star - u_\ell^{k,j} \vvvert
    +
    \vvvert u_\ell^{k,\star} - u_\ell^{k,j} \vvvert
    + \eta_\ell(u_\ell^{k,j}).
\]

Importantly, the following result states that the innermost \(j\)-loop is 
known to terminate. 
\begin{lemma}[bounded number of algebraic solvers steps {\cite[Lem.~3.2]{aisfem}}]
    \label{BMP:lemma:unif-steps-aisfem}
    Consider arbitrary parameters $0 < \theta \le 1$, $0 < \lambdalin $, 
    $0 < \lambdaalg$. Then, 
    Algorithm~\ref{BMP:algorithm:aisfem}
    guarantees that $\jj[\ell,k]  < \infty$
    for all $(\ell,k,0) \in \QQ$, i.e., the number of algebraic solver steps is finite.
    \qed
\end{lemma}

Note that there holds a-posteriori error control at each step 
of the adaptive algorithm, which can again be used to terminate 
Algorithm~\ref{BMP:algorithm:aisfem}.
\begin{lemma}[a-posteriori error control
    {\cite[Prop.~3.5]{aisfem}}]
	\label{BMP:prop:aposteriori-jk-aisfem}
    For
	any $(\ell,k,j) \in \QQ$ with $k \ge 1$ and $j \ge 1$, it holds that
    \begin{equation}
        \label{BMP:eq:aposteriori-jk-aisfem-1}
		\vvvert u^\star - u_\ell^{k,j}\vvvert
		\lesssim \big[ \eta_\ell(u_\ell^{k,j}) 
        + \vvvert u_\ell^{k,j} - u_\ell^{k-1,\jj}\vvvert
        + \vvvert u_\ell^{k,j} - u_\ell^{k,j-1}\vvvert \big],
    \end{equation}
    and for \( k = \kk[\ell]\), \( j = \jj[\ell,\kk]\)
    \begin{equation}
        \label{BMP:eq:aposteriori-jk-aisfem-2}
		\vvvert u^\star - u_\ell^{\kk,\jj}\vvvert
		\lesssim  \eta_\ell(u_\ell^{\kk,\jj}) .
    \end{equation}
\end{lemma}

\subsection{Convergence and complexity}

Algorithm~\ref{BMP:algorithm:aisfem} does not compute the
exact Zarantonello iterates with
contraction~\eqref{BMP:eq:Zarantonello_contraction}.
Nevertheless, the following lemma establishes that the
contraction property also holds for the inexact symmetrization
in order to replace Lemma~\ref{BMP:lem:contraction_quasi-error}.
\begin{lemma}[contraction of inexact Zarantonello iteration
    {\cite[Lem.~5.1]{bfmps2023}}]
    \label{BMP:lem:contraction_aisfem}
    There exists \(0<\lambdaalg^\star\) and \(0 < \qsym <1\) 
    such that, for all \(0<\lambdaalg<\lambdaalg^\star\), it holds
    \begin{align*}
        \vvvert u_\ell^\star - u_\ell^{k,\jj} \vvvert
        &\leq
        \qsym\,
        \vvvert u_\ell^\star - u_\ell^{k-1,\jj} \vvvert
        \quad\text{for all }
        (\ell, k, \jj) \in \QQ
        \text{ with } 1 \leq k < \kk[\ell],
        \\
        \vvvert u_\ell^\star - u_\ell^{\kk,\jj} \vvvert
        &\leq
        \qsym^\star\,
        \vvvert u_\ell^\star - u_\ell^{\kk-1,\jj} \vvvert
        +
        \frac{2 \qalg}{1 - \qalg}
        \lambdaalg
        \lambdasym
        \,
        \eta_\ell(u_\ell^{\kk,\jj})
        \quad\text{for all }
        (\ell, \kk, \jj) \in \QQ.
    \end{align*}
\end{lemma}

Since the non-symmetric problem~\eqref{BMP:eq:general-second-order}
ensures the Galerkin orthogonality~\eqref{BMP:eq:orthogonality} for 
\(b(\cdot,\cdot)\), but not for its symmetric part \(a(\cdot,\cdot)\),
the Pythagorean identity~\eqref{BMP:eq:Pythagoras} is
replaced by the quasi-orthogonality introduced as additional axiom
in~\cite{cfpp2014}.
The following generalization from~\cite{feischl2022}
has been applied in the optimal-complexity analysis in \cite{bfmps2023}.
\begin{lemma}[quasi-orthogonality {\cite[Prop.~2]{bfmps2023}}]
    \label{BMP:lem:quasi-orthogonality}
    There exists \(0 < \alpha \leq 1\) such that the
    following holds: For any sequence 
    \(\XX_\ell \subseteq \XX_{\ell+1} \subset \XX\)
    of nested subspaces, the corresponding exact Galerkin solutions
    \(u_\ell^\star \in \XX_\ell\)
    to~\eqref{BMP:eq:weak_nonsymmetric}
    satisfy
    \begin{align}\label{BMP:eq:quasi-orthogonality}
        \sum_{\ell' = \ell}^{\ell+N}
        \vvvert u_{\ell'+1}^\star - u_{\ell'}^\star \vvvert^2
        \lesssim
        (N+1)^{1 - \alpha}
        \,
        \vvvert u^\star - u_\ell^\star \vvvert^2
        \quad\text{for all }
        \ell, N \in \N_0.
    \end{align}
\end{lemma}

The two previous lemmas allow for the following
generalization of the convergence results from
Subsection~\ref{BMP:sec:convergence}.
Note that, in contrast to
Theorem~\ref{BMP:thm:full-linear-convergence} 
and Theorem~\ref{BMP:thm:full-linear-convergence-goafem}, 
full R-linear convergence is subject to some requirements on the 
stopping parameter 
\(\lambdaalg\)
stemming from Lemma~\ref{BMP:lem:contraction_aisfem}.
In this sense, the following 
theorem fails to state parameter-robust convergence.
Finally, we note that the present formulation 
is proved in~\cite{bfmps2023}, while an earlier 
result~\cite{aisfem} had even more severe restrictions 
on \(\lambdaalg\) and 
\(\lambdasym\).
\begin{theorem}[full R-linear convergence of Algorithm~\ref{BMP:algorithm:aisfem}, {\cite[Thm.~14]{bfmps2023}}]
    \label{BMP:thm:full-linear-convergence-aisfem}
    Let $0 < \theta \le 1$ and \(0 < \lambdasym \le 1\) be arbitrary.
    Let \(\delta > 0\) be sufficiently small to guarantee
    contraction~\eqref{BMP:eq:Zarantonello_contraction}
    of the exact Zarantonello iteration.
    Recall \(\lambdaalg^\star\) from Lemma~\ref{BMP:lem:contraction_aisfem}.
    Then, there exists an upper bound 
    \(\lambda^\star  > 0\)
    such that, for all 
    \(0< \lambdaalg < \lambdaalg^\star\) with 
    \( \lambdaalg \lambdasym < \lambda^\star\),
    the following holds:
    Algorithm~\ref{BMP:algorithm:aisfem}
    guarantees the existence of 
    $\Clin > 0$ and $0 < \qlin < 1$
    such that
    \begin{equation}
        \label{BMP:eq:full-linear-convergence-aisfem}
        \Eta_\ell^{k, j}
        \lesssim
        \qlin^{|\ell,k,j| - |\ell',k',j'|} \,
        \Eta_{\ell'}^{k',j'}
        \quad\text{for all }
        (\ell,k,j), (\ell',k',j') \in \QQ
        \text{ with }
        |\ell,k,j| > |\ell',k',j'|.
    \end{equation}
    In particular, this yields convergence
    \begin{align*}
        \vvvert u^\star - u_\ell^{k, j} \vvvert
        \lesssim \Eta_\ell^{k, j}
        \le
        \Clin \,
        \qlin^{
            \vert \ell, k, j \vert
        } \,
        \Eta_0^{0,0} \rightarrow 0 \quad \text{as }
        \vert \ell, k, j \vert \rightarrow \infty.
        \qed
    \end{align*}
\end{theorem}

With validity of linear 
convergence~\eqref{BMP:eq:full-linear-convergence-aisfem}, one can follow the 
lines of the proof of Theorem~\ref{BMP:thm:optimal-complexity-goafem} to 
conclude optimal complexity.
\begin{theorem}[optimal computational complexity {\cite{bfmps2023,aisfem}}]
    \label{BMP:thm:aisfem:complexity}
    Under the assumptions of
    Theorem~\ref{BMP:thm:full-linear-convergence-aisfem},
    there exist additional upper bounds 
    \(\theta^\star, \lambdasym^\star > 0\)
    such that, for all
    \(0 < \theta < \theta^\star\),
    \(0 < \lambdasym < \lambdasym^\star\),
    and
    \(0 < \lambdaalg < \lambdaalg^\star\),
    with \( \lambdaalg \lambdasym < \lambda^\star\),
    the following holds:
    Algorithm~\ref{BMP:algorithm:aisfem} guarantees,
    for all \(s > 0\), that
    \begin{align}
        \label{BMP:eq:aisfem:optimal:complexity}
        \| u^\star \|_{\A_s}
        \lesssim
        \sup_{(\ell,k,j) \in \QQ}
        \cost(\ell, k, j)^s
        \, \Eta_\ell^{k,j}
        \lesssim
        \max\{
            \| u^\star \|_{\A_s},
            \Eta_{0}^{0,0}
        \}.
    \end{align}
\end{theorem}

Based on the ideas from Section~\ref{BMP:sec:goafem}, one can extend the 
results in the Theorems~\ref{BMP:thm:full-linear-convergence-aisfem}--\ref{BMP:thm:aisfem:complexity}
for general second-order linear elliptic 
PDEs from standard AFEM to GOAFEM. Details are found in the 
recent work~\cite{goaisfem}.

\subsection{Numerical example}
\label{BMP:sec:aisfem:example}

The performance of Algorithm~\ref{BMP:algorithm:aisfem}
is tested using a benchmark example
on the L-shaped domain
\(\Omega \coloneqq (-1, 1)^2 \setminus [0, 1)^2\)
with dominating convection
\(\boldsymbol{b} \coloneqq (-10, -10)^\top\)
and \(c = 0\).
It seeks \(u^\star \in H^1(\Omega)\) satisfying
\[
    -\Delta u^\star +
    \begin{pmatrix}
        -10
        \\
        -10
    \end{pmatrix}
    \cdot
    \nabla u^\star
    =
    1
    \text{ in } \Omega
    \quad\text{with}\quad
    u^\star
    =
    0
    \text{ on } \partial\Omega.
\]

\begin{figure}
    \centering
    \subfloat[Adaptively refined mesh on level \(\ell = 11\) 
    with \(1\,336\) triangles.]{%
        \label{BMP:fig:aisfem:mesh}
        \begin{tikzpicture}
    \begin{axis}[%
        axis equal image,%
        width=4.8cm,%
        xmin=-1.15, xmax=1.15,%
        ymin=-1.15, ymax=1.15,%
        font=\footnotesize%
    ]
        \addplot graphics [xmin=-1, xmax=1, ymin=-1, ymax=1]
        {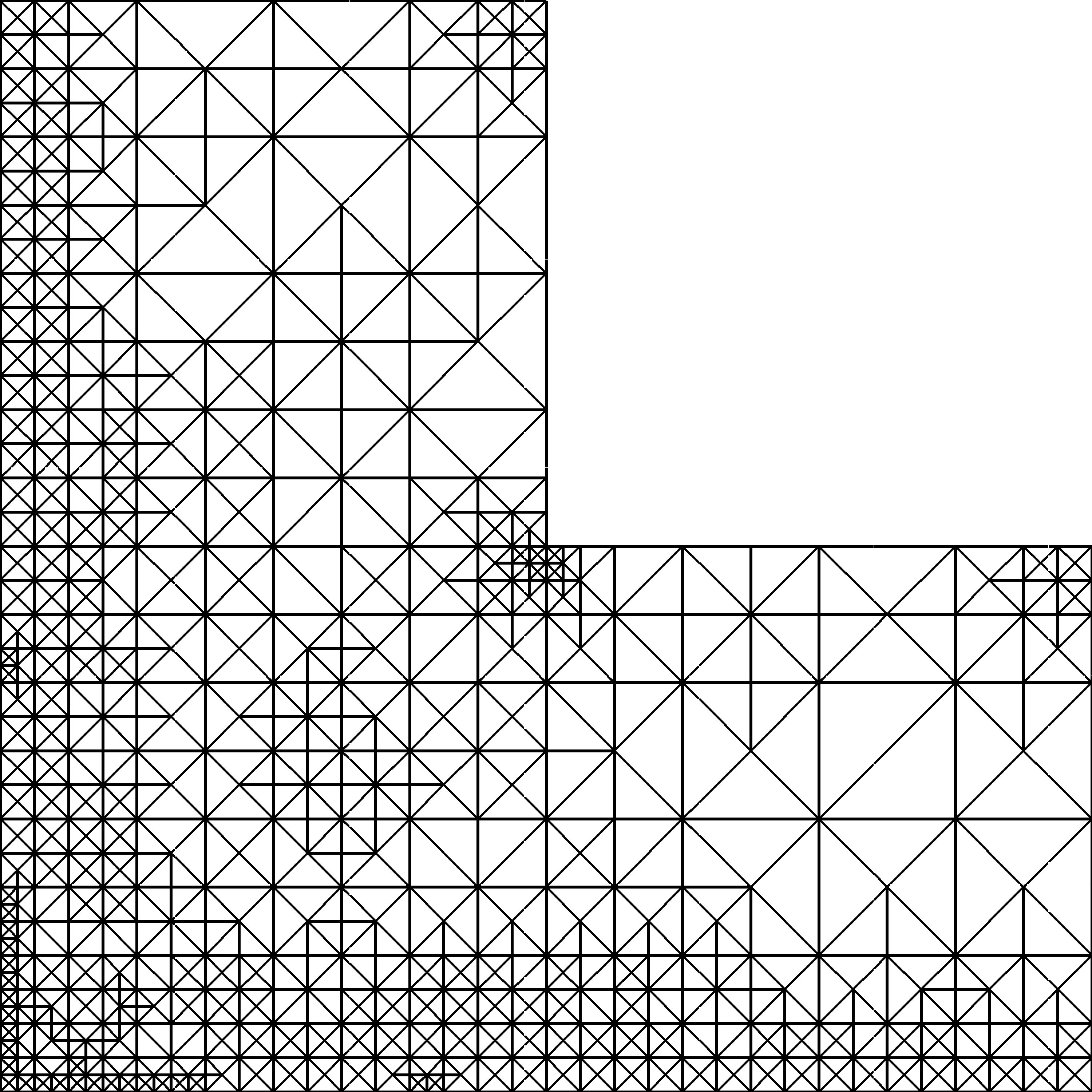};
    \end{axis}
\end{tikzpicture}
    }
    \hfil
    \subfloat[Discrete solution \(u_{11}^{\kk,\jj} \in \XX_{11}\).]{%
        \label{BMP:fig:aisfem:solution}
        \begin{tikzpicture}
    \pgfplotsset{/pgf/number format/fixed}
    \begin{axis}[%
        width=4.8cm,%
        xmin=-1.1, xmax=1.1,%
        ymin=-1.1, ymax=1.1,%
        zmin=-0.009, zmax=0.1,%
        font=\footnotesize,%
    ]
        \addplot3 graphics [%
            points={%
                (-1,-1,0) => (0.4,268.7-180.6)
                (1,-1,0) => (220.2,268.7-268.2)
                (0,1,0) => (375.5,268.7-150.8)
                (-0.75,-0.75,0.0823709) => (61.0,268.7-1.1)
            }%
            ]
            {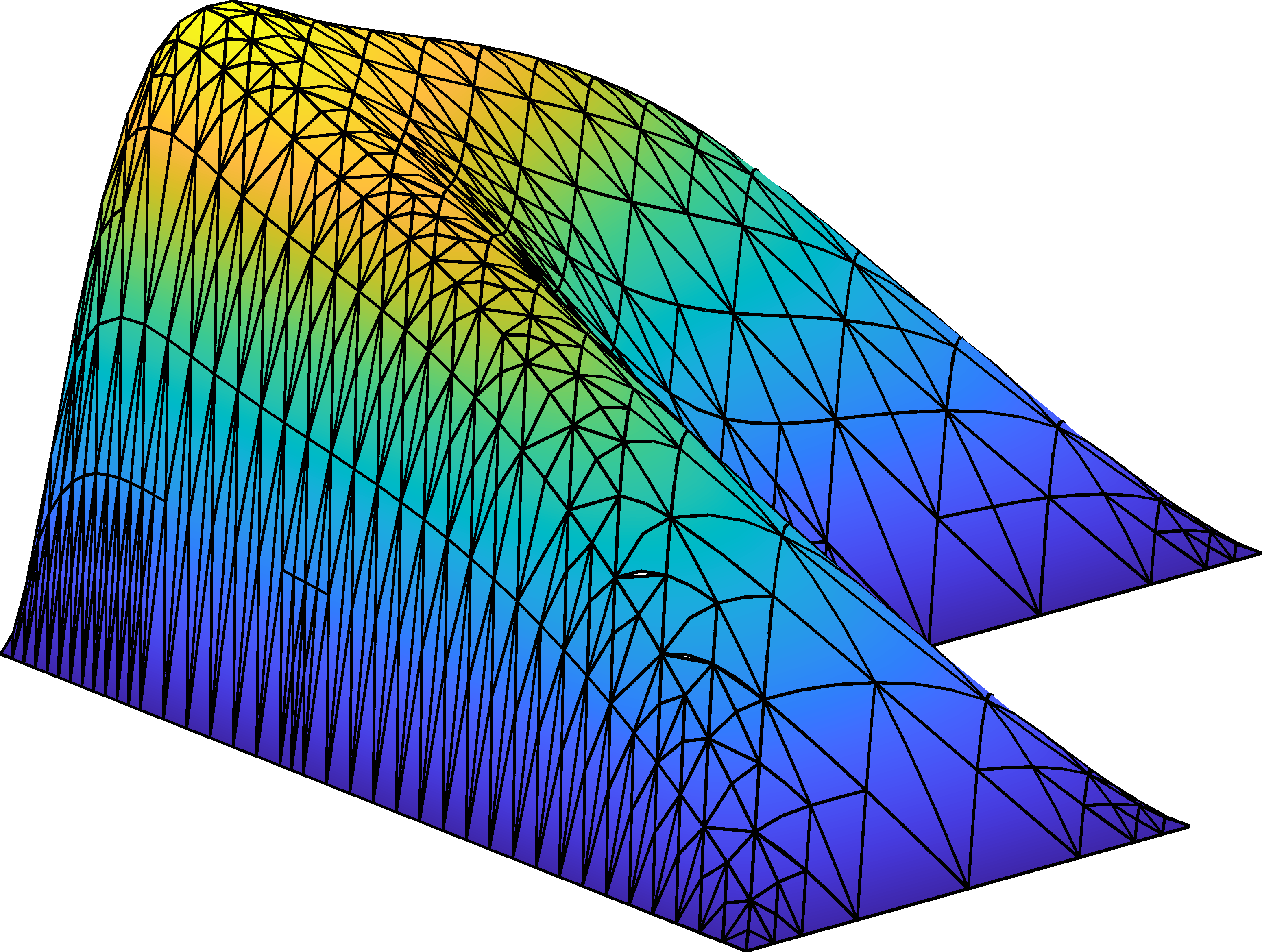};
    \end{axis}
\end{tikzpicture}
    }
    \caption{%
        Mesh and corresponding discrete solution
        for the benchmark problem from
        Subsection~\ref{BMP:sec:aisfem:example}.
        The results are generated by Algorithm~\ref{BMP:algorithm:aisfem} with
        polynomial degree \(p = 2\), damping parameter \(\delta = 0.1\),
        bulk parameter \(\theta = 0.3\),
        and stopping parameters \(\lambdasym = \lambdaalg = 0.1\).
    }
    \label{BMP:fig:aisfem:mesh_solution}
\end{figure}

\begin{figure}
    \centering
    \input{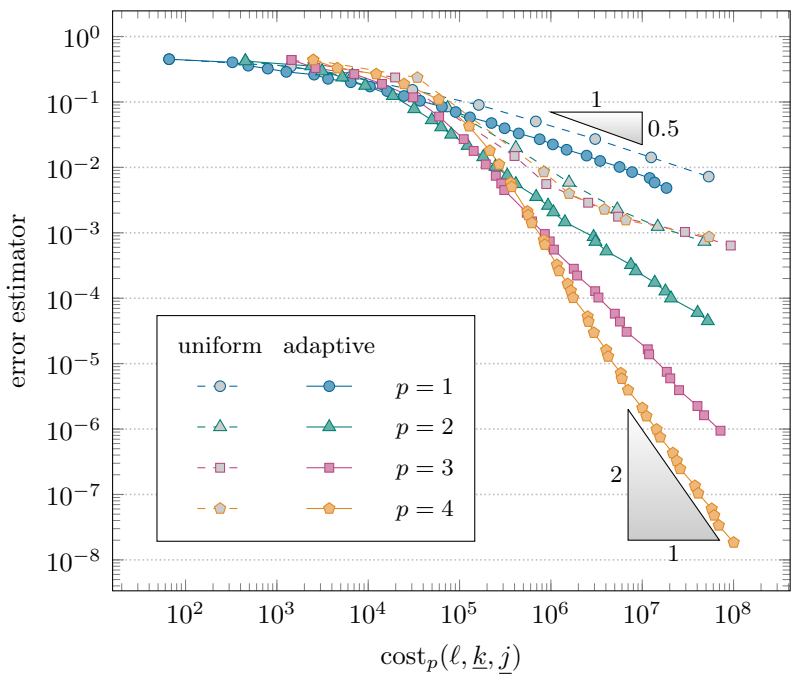}
    \caption{
        Convergence plot of the 
        adaptive Algorithm~\ref{BMP:algorithm:aisfem}
        to solve the benchmark problem from
        Subsection~\ref{BMP:sec:aisfem:example}
        for various polynomial degrees.
        The damping parameter reads \(\delta = 0.1\)
        and the adaptivity parameters \(\theta = 0.3\)
        and \(\lambdasym = \lambdaalg = 0.1\).
        All graphs display values of
        the residual-based error estimator
        \(\eta_\ell^{\kk,\jj} (u_\ell^{\kk,\jj})\) from~\eqref{BMP:eq:aisfem:estimator}
        for the final iterates on the level \(\ell \in
        \N_0\).
    }
    \label{BMP:fig:aisfem:convergence}
\end{figure}

\begin{figure}
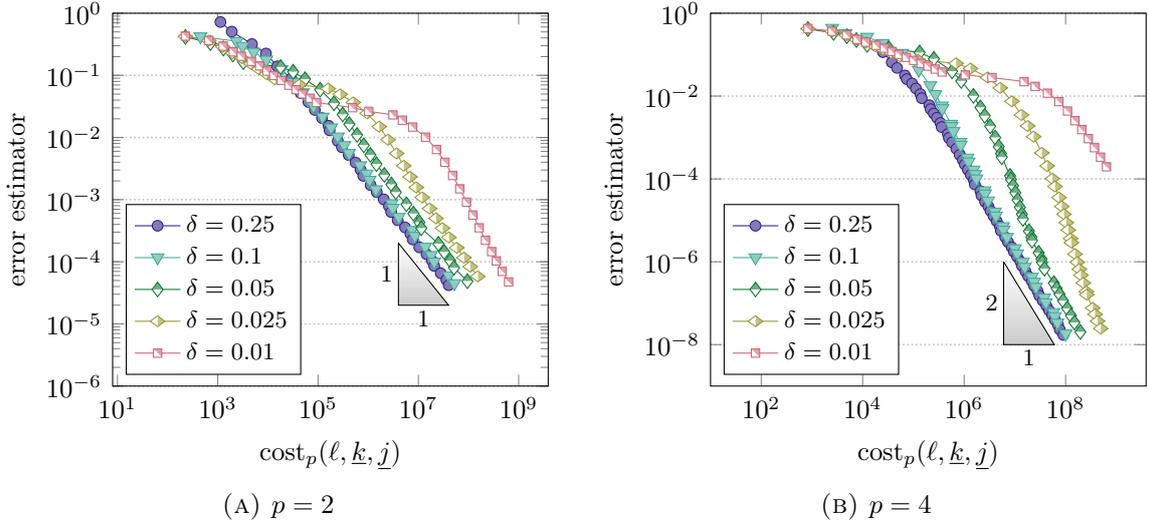

    \centering
    \subfloat[\(p = 2\)]{%
        \label{BMP:fig:aisfem:damping:p2}
        \input{chapter_afem_BMP/figures/Fig16a_AISFEM_damping_p2.tex}
    }
    \hfil
    \subfloat[\(p = 4\)]{%
        \label{BMP:fig:aisfem:damping:p4}
        \input{chapter_afem_BMP/figures/Fig16b_AISFEM_damping_p4.tex}
    }
    \caption{
        Convergence plot of the 
        adaptive Algorithm~\ref{BMP:algorithm:aisfem}
        to solve the benchmark problem from
        Subsection~\ref{BMP:sec:aisfem:example}
        for various damping parameters \(\delta\)
        and polynomial degrees \(p\).
        The adaptivity parameters read \(\theta = 0.3\),
        and \(\lambdasym = \lambdaalg = 0.1\).
        All graphs display values of
        the residual-based error estimator
        \(\eta_\ell^{\kk,\jj} (u_\ell^{\kk,\jj})\) from~\eqref{BMP:eq:aisfem:estimator}
        for the final iterates on the level \(\ell \in
        \N_0\).
    }
    \label{BMP:fig:aisfem:damping}
\end{figure}

In addition to increased refinement due to the singularity
at the reentrant corner,
the strong convection term causes additional refinement
at the lower and left boundary layer;
see Figure~\ref{BMP:fig:aisfem:mesh}.
This is due to the large gradients of the discrete solution
as shown in Figure~\ref{BMP:fig:aisfem:solution}.
Therewith, the adaptive algorithm recovers
optimal convergence rates of the residual-based error estimator
for several polynomial degrees
in Figure~\ref{BMP:fig:aisfem:convergence}.
A closer investigation of the involved stopping parameters
\(\lambdasym\) and \(\lambdaalg\) in Table~\ref{BMP:table:aisfem:parameters}
in the fashion of Table~\ref{BMP:tab:Kellogg:parameters} in
Section~\ref{BMP:sec:numerics}
supports the choice of a comparably small value
\(\lambdasym = 0.1\)
for an optimal performance independently of the parameters
\(\theta\) and \(\lambdaalg\).
In the present case of a dominating convection term,
the choice of the damping parameter \(\delta\) is crucial.
A large norm of \(\boldsymbol{b}\) results in a strong deviation
from the symmetric bilinear form \(a (\cdot, \cdot) \) and, thus,
the Zarantonello iteration requires a small \(\delta \ll 1\)
to ensure the contraction property~\eqref{BMP:eq:Zarantonello_contraction}.
However, a small \(\delta\) leads to a slow convergence.
Figure~\ref{BMP:fig:aisfem:damping} displays the
convergence of the error estimator for various damping
parameters.
For all larger values of \(\delta\), no significant convergence
could be achieved.
This justifies the choice of \(\delta = 0.1\) in the
remaining experiments of this section.

\begin{table}
    \subfloat[\(\theta = 0.3\)]{%
        \label{BMP:table:aisfem:theta30}
        %
%
\pgfplotstableset{%
    trim cells = true,%
    col sep = comma,%
    row sep = newline%
}
%
%
\pgfplotstableset{%
    fixed zerofill,%
    column type = {r},%
    every head row/.style = {before row = \toprule, after row = \midrule},%
    every last row/.style = {after row = \bottomrule},%
    columns/lambdaAlg/.style = {%
        fixed,%
        precision = 1,%
        column  name={{$\lambdaalg \setminus$}},%
    },%
    columns/lambdaSym0.9/.style = {zerofill, column name = {$0.9$}},%
    columns/lambdaSym0.7/.style = {zerofill, column name = {$0.7$}},%
    columns/lambdaSym0.5/.style = {zerofill, column name = {$0.5$}},%
    columns/lambdaSym0.3/.style = {zerofill, column name = {$0.3$}},%
    columns/lambdaSym0.1/.style = {zerofill, column name = {$\lambdasym$ \qquad $0.1$}},%
    rowmin/.style = {%
        postproc cell content/.append style={
            /pgfplots/table/@cell content/.add={\cellcolor{colorrowmin}}{}
        },
    },
    colmin/.style = {%
        postproc cell content/.append style={
            /pgfplots/table/@cell content/.add={\cellcolor{colorcolmin}}{}
        }
    },
    bothmin/.style = {%
        postproc cell content/.append style={
            /pgfplots/table/@cell content/.add={\cellcolor{colorbothmin}}{}
        }
    },
    highlight col min/.code 2 args = {%
        \pgfmathtruncatemacro\rowindex{#1-1}
        \edef\setstyles{%
            \noexpand\pgfplotstableset{%
                every row \rowindex\noexpand\space column #2/.style={colmin}
            }%
        }\setstyles
    },
    highlight row min/.code 2 args = {%
        \pgfmathtruncatemacro\rowindex{#1-1}
        \edef\setstyles{%
            \noexpand\pgfplotstableset{%
                every row \rowindex\noexpand\space column #2/.style={rowmin}
            }%
        }\setstyles
    },
    highlight both min/.code 2 args = {%
        \pgfmathtruncatemacro\rowindex{#1-1}
        \edef\setstyles{%
            \noexpand\pgfplotstableset{%
                every row \rowindex\noexpand\space column #2/.style={bothmin}
            }%
        }\setstyles
    },
}
%
%
\pgfplotstableread{
lambdaAlg,lambdaSym0.1,lambdaSym0.3,lambdaSym0.5,lambdaSym0.7,lambdaSym0.9
0.1,0.023721362097746,0.0247140670156256,0.0313330091067829,0.0367161460277687,0.0479096381110999
0.3,0.0190169068894587,0.0215034453768888,0.0287215050264754,0.0431922313818154,0.0517092724738407
0.5,0.0190092960032165,0.0250181399496875,0.0321007059058708,0.0442245776389518,0.0532836915128865
0.7,0.0190315370557475,0.0239195363059229,0.029929701540542,0.0400143429446262,0.0545333178994897
0.9,0.0203221607961539,0.0233720638214165,0.0291920854836623,0.0437863884056162,0.054571479729614
}\data
%
%
\pgfplotstabletypeset[%
    highlight col min={2}{2},
    highlight col min={2}{3},
    highlight col min={1}{4},
    highlight col min={1}{5},
    highlight row min={1}{1},
    highlight row min={2}{1},
    highlight row min={3}{1},
    highlight row min={4}{1},
    highlight row min={5}{1},
    highlight both min={2}{1},
    highlight both min={3}{1},
    highlight both min={4}{1},
]{\data}

    }

    \subfloat[\(\theta = 0.5\)]{%
        \label{BMP:table:aisfem:theta50}
        %
%
\pgfplotstableset{%
    trim cells = true,%
    col sep = comma,%
    row sep = newline%
}
%
%
\pgfplotstableset{%
    fixed zerofill,%
    column type = {r},%
    every head row/.style = {before row = \toprule, after row = \midrule},%
    every last row/.style = {after row = \bottomrule},%
    columns/lambdaAlg/.style = {%
        fixed,%
        precision = 1,%
        column  name={{$\lambdaalg \setminus$}},%
    },%
    columns/lambdaSym0.9/.style = {zerofill, column name = {$0.9$}},%
    columns/lambdaSym0.7/.style = {zerofill, column name = {$0.7$}},%
    columns/lambdaSym0.5/.style = {zerofill, column name = {$0.5$}},%
    columns/lambdaSym0.3/.style = {zerofill, column name = {$0.3$}},%
    columns/lambdaSym0.1/.style = {zerofill, column name = {$\lambdasym$ \qquad $0.1$}},%
    rowmin/.style = {%
        postproc cell content/.append style={
            /pgfplots/table/@cell content/.add={\cellcolor{colorrowmin}}{}
        },
    },
    colmin/.style = {%
        postproc cell content/.append style={
            /pgfplots/table/@cell content/.add={\cellcolor{colorcolmin}}{}
        }
    },
    bothmin/.style = {%
        postproc cell content/.append style={
            /pgfplots/table/@cell content/.add={\cellcolor{colorbothmin}}{}
        }
    },
    highlight col min/.code 2 args = {%
        \pgfmathtruncatemacro\rowindex{#1-1}
        \edef\setstyles{%
            \noexpand\pgfplotstableset{%
                every row \rowindex\noexpand\space column #2/.style={colmin}
            }%
        }\setstyles
    },
    highlight row min/.code 2 args = {%
        \pgfmathtruncatemacro\rowindex{#1-1}
        \edef\setstyles{%
            \noexpand\pgfplotstableset{%
                every row \rowindex\noexpand\space column #2/.style={rowmin}
            }%
        }\setstyles
    },
    highlight both min/.code 2 args = {%
        \pgfmathtruncatemacro\rowindex{#1-1}
        \edef\setstyles{%
            \noexpand\pgfplotstableset{%
                every row \rowindex\noexpand\space column #2/.style={bothmin}
            }%
        }\setstyles
    },
}
%
%
\pgfplotstableread{
lambdaAlg,lambdaSym0.1,lambdaSym0.3,lambdaSym0.5,lambdaSym0.7,lambdaSym0.9                               
0.1,0.0205125621480996,0.0224694329083864,0.0329314508806981,0.0335878926551897,0.0489343923370489
0.3,0.0178148458415776,0.0250980079256744,0.035169153507588,0.0404170354001837,0.0524715038138848 
0.5,0.0179157834995145,0.018929501025451,0.0343662319420812,0.0422165619752556,0.0626148225947492 
0.7,0.0177816609277936,0.0197356905219535,0.0348112985581539,0.0406511739768613,0.0580408012240886
0.9,0.0178010056598496,0.0192292843134526,0.0349462718181699,0.0426003490722103,0.0591567302259911
}\data
%
%
\pgfplotstabletypeset[%
    highlight col min={2}{1},
    highlight col min={4}{1},
    highlight col min={5}{1},
    highlight col min={3}{2},
    highlight col min={1}{3},
    highlight col min={1}{4},
    highlight col min={1}{5},
    highlight row min={1}{1},
    highlight row min={3}{1},
    highlight both min={2}{1},
    highlight both min={4}{1},
    highlight both min={5}{1},
]{\data}

    }

    \subfloat[\(\theta = 0.7\)]{%
        \label{BMP:table:aisfem:theta70}
        %
%
\pgfplotstableset{%
    trim cells = true,%
    col sep = comma,%
    row sep = newline%
}
%
%
\pgfplotstableset{%
    fixed zerofill,%
    column type = {r},%
    every head row/.style = {before row = \toprule, after row = \midrule},%
    every last row/.style = {after row = \bottomrule},%
    columns/lambdaAlg/.style = {%
        fixed,%
        precision = 1,%
        column  name={$\lambdaalg \setminus$},%
    },%
    columns/lambdaSym0.9/.style = {zerofill, column name = {$0.9$}},%
    columns/lambdaSym0.7/.style = {zerofill, column name = {$0.7$}},%
    columns/lambdaSym0.5/.style = {zerofill, column name = {$0.5$}},%
    columns/lambdaSym0.3/.style = {zerofill, column name = {$0.3$}},%
    columns/lambdaSym0.1/.style = {zerofill, column name = {$\lambdasym$ \qquad $0.1$}},%
    rowmin/.style = {%
        postproc cell content/.append style={
            /pgfplots/table/@cell content/.add={\cellcolor{colorrowmin}}{}
        },
    },
    colmin/.style = {%
        postproc cell content/.append style={
            /pgfplots/table/@cell content/.add={\cellcolor{colorcolmin}}{}
        }
    },
    bothmin/.style = {%
        postproc cell content/.append style={
            /pgfplots/table/@cell content/.add={\cellcolor{colorbothmin}}{}
        }
    },
    highlight col min/.code 2 args = {%
        \pgfmathtruncatemacro\rowindex{#1-1}
        \edef\setstyles{%
            \noexpand\pgfplotstableset{%
                every row \rowindex\noexpand\space column #2/.style={colmin}
            }%
        }\setstyles
    },
    highlight row min/.code 2 args = {%
        \pgfmathtruncatemacro\rowindex{#1-1}
        \edef\setstyles{%
            \noexpand\pgfplotstableset{%
                every row \rowindex\noexpand\space column #2/.style={rowmin}
            }%
        }\setstyles
    },
    highlight both min/.code 2 args = {%
        \pgfmathtruncatemacro\rowindex{#1-1}
        \edef\setstyles{%
            \noexpand\pgfplotstableset{%
                every row \rowindex\noexpand\space column #2/.style={bothmin}
            }%
        }\setstyles
    },
}
%
%
\pgfplotstableread{
lambdaAlg,lambdaSym0.1,lambdaSym0.3,lambdaSym0.5,lambdaSym0.7,lambdaSym0.9                                 
0.1,0.0231241702192104,0.0232317115657584,0.0327703906707659,0.0435460737165996,0.043934212603266
0.3,0.0167959274844345,0.0252755456129996,0.0313900143494311,0.0528232137249066,0.065282169440460
0.5,0.0166175041340151,0.0256779327859868,0.0335498589853174,0.0536555857350605,0.067120032204729
0.7,0.0165252195196711,0.0266164098345122,0.0337368883814511,0.0558195933296381,0.067023807810640
0.9,0.0170235249844349,0.025779544351121,0.0345356212146683,0.0554371301417296,0.0655973916169193
}\data
%
%
\pgfplotstabletypeset[%
    highlight col min={1}{2},
    highlight col min={2}{3},
    highlight col min={1}{4},
    highlight col min={1}{5},
    highlight row min={1}{1},
    highlight row min={2}{1},
    highlight row min={3}{1},
    highlight row min={4}{1},
    highlight row min={5}{1},
    highlight both min={4}{1},
]{\data}

    }

    \caption{
        \normalfont\footnotesize\quad
        Investigation of the influence of the adaptivity parameters
        \(\lambdasym\) and \(\lambdaalg\) on the performance
        of Algorithm~\ref{BMP:algorithm:aisfem}
        to solve the benchmark problem from
        Subsection~\ref{BMP:sec:aisfem:example}
        for various polynomial degrees.
        The damping parameter reads \(\delta = 0.1\)
        and the polynomial degree \(p = 2\).
        The best choice per column
        is marked in yellow, per row in blue, and for both in green.
        In either case \(\theta \in \{0.3, 0.5, 0.7 \}\), the best choice 
        is \(\lambdasym = 0.1\) and \(\lambdaalg = 0.7\).
    }
    \label{BMP:table:aisfem:parameters}
\end{table}

%
%

\section{Extension 3: Non-linear problems}
\label{BMP:sec:ailfem}

The Zarantonello iteration from
Subsection~\ref{BMP:sec:aisfem}
is one possible choice for the linearization
of non-linear problems.
This leads to an analogous adaptive algorithm
to Algorithm~\ref{BMP:algorithm:aisfem},
where the Zarantonello symmetrization~\eqref{BMP:eq:Zarantonello}
is replaced by a suitable linearization scheme
with non-linear residual.
The convergence results from
Theorems~\ref{BMP:thm:full-linear-convergence-aisfem}--\ref{BMP:thm:aisfem:complexity}
apply verbatim in this setting.
For further details, the reader is referred to~\cite{bfmps2023} and 
\cite{hpsv2021}.
Instead, this section presents a novel energy-based
approach from \cite{mps2024} that allows for a parameter-robust
alternative.

\subsection{Strongly monotone model problem}

Throughout this subsection,
consider a non-linear mapping
\(\boldsymbol{A}: \R^d \to \R^d\)
in the non-linear elliptic PDE
\begin{align}
    \label{BMP:eq:nonlinear}
    - \div(\boldsymbol{A}(\nabla u^\star))
    =
    f - \div \boldsymbol{f}
    \text{ in } \Omega
    \quad\text{and}\quad
    u^\star = 0
    \text{ on } \partial\Omega.
\end{align}
The analysis in \cite{mps2024} is based on an
energy-minimization setting and, thus, is restricted
to scalar non-linearities of the form
\(\boldsymbol{A}(\xi) \coloneqq \mu(\vert \xi \vert^2) \xi\)
for some continuous \(\mu: \R_{\geq 0} \to \R_{\geq 0}\).
The weak formulation of~\eqref{BMP:eq:nonlinear}
seeks \(u^\star \in \XX\) such that
\begin{equation}
    \label{BMP:eq:weak_nonlinear}
    \langle \mathcal{A}(u^\star), v \rangle
    =
    F(v)
    \quad\text{for all }
    v \in \XX, 
\end{equation}
where \(\mathcal{A}: \XX \to \XX^*\)
denotes the non-linear operator 
defined by
\[
    \mathcal{A}(u)
    \coloneqq
    \langle \boldsymbol{A}(\nabla u),
    \nabla(\cdot) \rangle_{L^2(\Omega)}
    \coloneqq
    \int_\Omega \boldsymbol{A}(\nabla u) \cdot \nabla (\cdot) \, \d{x}.
\]
The non-linear coefficient \(\mu\) is assumed to satisfy a linear
growth condition, 
i.e.,
there exist \(0 < \alpha \leq L\) such that
\begin{align}
    \label{BMP:eq:mu_growth}
	\alpha (t-s)
	\le
	\mu(t^2) t - \mu(s^2) s
	\le
	\frac{L}{3} \, (t-s)
	\quad \text{for all }
	t \ge s \ge 0,
\end{align}
which is, for instance, satisfied for typical non-linearities in the 
context of magnetostatics.
This condition ensures that \(\mathcal{A}\) is strongly monotone
and Lipschitz continuous, i.e.,
for all \(u,v,w \in \XX\), there holds
\[
    \alpha\,
    \vvvert u - v \vvvert^2
    \leq
    \langle \mathcal{A}(u) - \mathcal{A}(v), u - v \rangle
    \text{ and }
    \langle \mathcal{A}(u) - \mathcal{A}(v), w \rangle
    \leq
    L\,
    \vvvert u - v \vvvert\,
    \vvvert w \vvvert.
\]
This provides existence and uniqueness of the solution \(u^\star\in\XX\) 
of the weak formulation~\eqref{BMP:eq:weak_nonlinear}; see~\cite[Sect.~25.4]{Zeidler1990}.

The discretization uses piecewise affine functions
since stability~\eqref{BMP:axiom:stability}
is only known for this case \cite[Rem.~6.2]{gmz2012}, while empirically also larger 
\(p > 1\) works in practice.
The discrete formulation seeks
\(u_H^\star \in \XX_H 
\)
with
\begin{equation}
    \label{BMP:eq:discrete_nonlinear}
    \langle \mathcal{A}(u_H^\star), v_h \rangle
    =
    F(v_H)
    \quad\text{for all } v_H \in \XX_H.
\end{equation}
Existence and uniqueness of \(u_H^\star\in\XX_H\) follow from the abstract 
arguments that also apply to the weak 
formulation~\eqref{BMP:eq:weak_nonlinear}.
The corresponding residual-based a-posteriori
error estimator reads
\(
    \eta_H(u_H)^2
    \coloneqq
    \sum_{T \in \TT_H}
    \eta_H(T, u_H)^2
\)
with the local contributions
\begin{equation}
    \label{BMP:eq:ailfem:estimator}
    \begin{split}
        \eta_H(T, u_H)^2
        &\coloneqq
        \vert T \vert^{2/d} \,
        \Vert f + \div(\mu(\vert \nabla u_H \vert^2) \nabla u_H - \boldsymbol{f})
        \Vert_{L^2(T)}^2
        \\
        &\phantom{{}\coloneqq{}}
        +
        \vert T \vert^{1/d}
        \Vert
            \lbracket
            (\mu(\vert \nabla u_H \vert^2) \nabla u_H - \boldsymbol{f}) \cdot n
            \rbracket
        \Vert_{L^2(\partial T \cap \Omega)}^2.
    \end{split}
\end{equation}
It is well-known from~\cite[Section~10.1]{cfpp2014} that \(\eta_H\)
satisfies~\eqref{BMP:axiom:stability}--\eqref{BMP:axiom:discrete_reliability}.

The (non-linear) discrete 
problem~\eqref{BMP:eq:discrete_nonlinear},
could be treated with the Zarantonello iteration or alternatively 
a damped Newton method. However, in the 
present setting we employ the Ka\v{c}anov 
linearization~\cite{kacanov}, see also~\cite{hpw2021}, which does not 
rely on the choice of appropriate additional 
damping parameters, thus guaranteeing parameter-robust 
full R-linear convergence. Instead, the Ka\v{c}anov linearization 
exploits the multiplicative
structure of the non-linearity
\(
	\boldsymbol{A}(\nabla u_H)
	=
	\mu(| \nabla u_H |^2) \, \nabla u_H
\).
The corresponding iteration mapping $\Phi_H \colon \XX_H \to \XX_H$ 
is determined via
\begin{equation}\label{BMP:eq:kacanov_iter}
	\langle \mu(| \nabla u_H |^2) \, \nabla \Phi_H(u_H),
    \nabla v_H \rangle_{L^2(\Omega)}
    =
    F(v_H)
	\quad
	\text{for all}
	\quad
	u_H, v_H \in \XX_H.
\end{equation}
Since \( \alpha \le \mu(| \nabla u_H |^2) \le L/3 \)
from~\eqref{BMP:eq:mu_growth} with \(s=0\), this discrete linear system is 
well-posed and admits a unique solution \(\Phi_H(u_H) \in \XX_H\) so that 
\(\Phi_H\) is indeed well-defined. Moreover, \eqref{BMP:eq:kacanov_iter} 
takes the form of the linear system from 
Section~\ref{BMP:sec:model_problem}
and can thus be solved using the appropriate
iterative algebraic solver
\(\Psi_H: \XX^* \times \XX_H \to \XX_H\)
from Subsection~\ref{BMP:sec:solve}.
Further details are omitted here and the reader is referred to 
\cite[Sect.~2.3]{mps2024}.
As in Section~\ref{BMP:sec:aisfem}, the adaptive algorithm will thus 
require a triple index \((\ell, k, j)\): The index \(\ell \in \N_0\) denotes 
the mesh-level, the index \(k \in \N_0\) denotes the Ka\v{c}anov step, and 
the index \(j \in \N_0\) denotes the step of the algebraic solver, i.e., 
\[ 
    u_\ell^{k,j} 
    =
    \Psi_\ell(u_\ell^{k,\star}, u_\ell^{k,j-1}) 
    \approx
    u_\ell^{k,\star}
    \coloneqq
    \Phi_\ell (u_\ell^{k-1,\jj})
    \approx
    u_\ell^\star, 
\] 
where neither \( u_\ell^{k,\star}\) nor 
\(u_\ell^{\star}\) are computed, but only \(u_\ell^{k,j}\).

Solving the weak problem~\eqref{BMP:eq:weak_nonlinear} (resp.~the discrete 
problem~\eqref{BMP:eq:discrete_nonlinear}) can indeed be 
seen as minimizing the energy \(\mathcal{E}\)
\begin{equation}
    \label{BMP:eq:energy}
    \mathcal{E}(v)
    \coloneqq
    \frac{1}{2}
    \int_\Omega
    \int_{0}^{\vert \nabla v \vert^2}
    \mu(t)\, \d t\, \d x
    -
    F(v)
    \quad\text{for all }
    v \in \XX
\end{equation}
over \(v \in \XX\) with unique minimizer \(u^\star \in \XX \) (resp.~over 
\(v_H \in \XX_H\) with unique minimizer \(u^\star_H \in \XX_H \)); see, e.g., 
\cite{mps2024}.
The corresponding energy difference
\[
    \dist^2(v, w)
    \coloneqq
    \mathcal{E}(w) - \mathcal{E}(v)
    \quad\text{for all }
    v, w \in \XX
\]
satisfies 
\[
    \frac{\alpha}{2} \, \vvvert v_{H} - u_{H}^\star \vvvert^2
    \le \dist^2(u_{H}^\star, v_{H})
    \le \frac{L}{2} \, \vvvert v_{H} - u_{H}^\star \vvvert^2
    \quad \text{for all } v_{H} \in \XX_{H};
\]
see, e.g.,~\cite[Lemma~5.1]{ghps2018}. Since the Ka\v{c}anov linearization can 
be proved to be contractive in the energy difference, we employ 
\(\dist^2(\cdot, \cdot)\) 
in a novel stopping criterion
for the adaptive algorithm in~\cite{mps2024}.

\begin{remark} 
    The reader is also referred to \cite{bbimp2022,bbimp2023,bps2024}
    for results addressing locally Lipschitz semi-linear problems, 
    where linearization is based on the Zarantonello iteration from 
    Section~\ref{BMP:sec:aisfem}.
\end{remark}

\subsection{Algorithm 4: Adaptive iteratively linearized FEM}
\labeltext{4}{BMP:algorithm:ailfem}

\noindent
{\bfseries Input:}
Initial mesh \(\TT_0\),
polynomial degree \(p \in \N\),
initial iterate $u_0^{0,0} \coloneqq u_0^{0,\underline{j}} \coloneqq 0$,
marking parameter \(0 < \theta \le 1\),
stopping parameter $\lambdalin > 0$,
and algebraic solver parameter $0 < \rho < 1$.
Moreover, initialize $\alpha_{\textup{min}} > 0$, 
$J_{\textup{max}} \in \N$ with
arbitrary values.
\smallskip

\noindent
{\bfseries Adaptive loop:}
For all \(\ell = 0, 1, 2, \dots\),
repeat the steps~(I)--(III):
\begin{enumerate}[leftmargin=2.5em]
    \item[\rmfamily(I)]
        {\ttfamily SOLVE \& ESTIMATE.}
        For all \(k = 1, 2, 3, \dots\),
        repeat the steps (a)--(c):
        \begin{itemize}[leftmargin=1.7em]
            \item[\rmfamily(a)]
                \textbf{Algebraic solver loop:}
                For all \(j = 1, 2, 3, \dots\),
                repeat the steps~(i)--(ii):
                \begin{itemize}[leftmargin=2em]
                    \item[\rmfamily(i)]
                        Compute
                        \(
                            u_\ell^{k,j}
                            \coloneqq
                            \Psi_\ell(u_\ell^{k,\star}, u_\ell^{k,j-1})
                        \)
                        (for the approximation of the theoretical quantity
                        $u_\ell^{k, \star}
                        \coloneqq
                        \Phi_\ell(u_\ell^{k-1,\underline j})
                        \in \XX_\ell$ solving~\eqref{BMP:eq:kacanov_iter} 
                        for \(u_H = u_\ell^{k-1,\underline j}\)), 
                        refinement indicators
                        \(\eta_\ell(T, u_\ell^{k,j})\) for
                        all \(T \in \TT_\ell\),
                        and the energy difference
                        \(\dist^2(u_\ell^{k,j}, u_\ell^{k-1,\jj})\).
                    \item[\rmfamily(ii)]
                        Compute
                        \(
                            \alpha_\ell^{k,j}
                            \coloneqq
                            \dist^2(u_\ell^{k,j}, u_\ell^{k-1,\jj})
                            /
                            \vvvert u_\ell^{k,j} - u_\ell^{k-1,\jj} \vvvert^2
                        \).
                    \item[\rmfamily(iii)]
                        Terminate \(j\)-loop with \(\jj[\ell,k] \coloneqq j\)
                        and employ nested iteration \(u_\ell^{k, 0} \coloneqq u_\ell^{k-1,\jj}\)
                        provided that either
                        \[
                            \alpha_\ell^{k,j} \ge \alpha_{\textup{min}}
                            \quad\text{or}\quad
                            u_\ell^{k,j} = u_\ell^{k-1,\jj}
                            \quad\text{or}\quad
                            \big[
                                \, \alpha_\ell^{k,j} > 0
                                \text{ \& } 
                                j > J_{\textup{max}} \,
                            \big]
                        \]
                \end{itemize}
            \item[\rmfamily(b)]
                If $\jj[\ell,k] > J_{\textup{max}}$,
                then update $J_{\textup{max}} \coloneqq \jj[\ell,k]$ and
                $\alpha_{\textup{min}} \coloneqq \rho \, \alpha_{\textup{min}}$.
            \item[\rmfamily(c)]
                Terminate \(k\)-loop and
                define \(\kk[\ell] \coloneqq k\)
                provided that
                \begin{align}
                    \label{BMP:eq:ailfem:linstopping}
                    \dist^2(u_\ell^{k,\jj},u_\ell^{k-1,\jj})
                    \le
                    \lambdalin\eta_\ell(u_\ell^{k,\jj})^2.
                \end{align}
        \end{itemize}
    \item[\rmfamily(II)]
        {\ttfamily MARK.}
        Determine a set \(\MM_\ell \subseteq \TT_\ell\)
        of minimal cardinality that satisfies
        the D\"orfler marking criterion
        \[
            \theta \, \eta_\ell(u_\ell^{\kk,\jj})^2
            \le
            \eta_\ell(\MM_\ell, u_\ell^{\kk,\jj})^2.
        \]
        \item[\rmfamily(III)]
            {\ttfamily REFINE.}
            Generate the new mesh
            \(
                \TT_{\ell+1}
                \coloneqq
                \texttt{refine}(\TT_\ell, \MM_\ell)
            \)
            by NVB and define
            \(
                u_{\ell+1}^{0,0}
                \coloneqq
                u_{\ell+1}^{0,\jj}
                \coloneqq
                u_\ell^{\kk,\jj}
            \)
            (nested iteration).
\end{enumerate}

In practice, Algorithm~\ref{BMP:algorithm:ailfem} 
will be stopped after finite time if the 
number of degrees of freedom exceeds a given threshold or if the computational 
time has reached a maximum or if the approximation is sufficiently accurate.

This algorithm coincides with Algorithm~\ref{BMP:algorithm:aisfem}
except for the novel stopping condition in step~(I.a.iii)
for the algebraic solver loop
and in step~(I.c) for the linearization loop that exploits the inherent 
energy structure of the PDE~\eqref{BMP:eq:nonlinear}.

\begin{remark}
    Let us comment on the stopping criteria in 
    Algorithm~\ref{BMP:algorithm:ailfem} for the Ka\v{c}anov 
    linearization and the algebraic solver.
    \begin{enumerate}[label={\normalfont (\roman*)}]
    \item Note the novel and parameter-free stopping 
    criterion for the algebraic solver loop in step~(I.a.iii) of 
    Algorithm~\ref{BMP:algorithm:ailfem} when compared to 
    Algorithm~\ref{BMP:algorithm:aisfem}. It 
    is motivated by the fact that the Ka\v{c}anov 
    linearization guarantees existence of a constant 
    \( C_{\text{nrg}}^\star>0\) such that
    \begin{align}
        \label{BMP:eq:energy_kacanov}
        \alpha 
        \coloneqq
            \dist^2(u_\ell^{k,\star}, u_\ell^{k-1,\jj})
            /
            \vvvert u_\ell^{k,\star} - u_\ell^{k-1,\jj} \vvvert^2
            \ge 
            C_{\text{nrg}}^\star>0;
    \end{align}
    see, e.g.,~\cite{hpw2021}.
    Since this holds analytically for \(u_\ell^{k,\star}\), the algebraic 
    solver should iterate until this condition is satisfied for the 
    computed approximation \(u_\ell^{k,j}\). Indeed, the 
    estimate~\eqref{BMP:eq:energy_kacanov} ensures that the 
    nested linearization with algebraic solver is contractive in energy, i.e., 
    \[
        \dist^2(u_\ell^{\star}, u_\ell^{k,\jj})
        \le
        \qlin \, \dist^2(u_\ell^{\star}, u_\ell^{k-1,\jj}),
        \]
    where \(0<\qlin< 1\) is independent of \(\ell\) and \(k\). 
    The latter contraction 
    is the crucial ingredient in proving full R-linear convergence for the 
    adaptive algorithm, where the proof follows the above lines 
    essentially replacing 
    \( \vvvert \cdot \vvvert^2 \) by the energy difference \(\dist^2 (\cdot,\cdot)\).
    \item With this understanding, the stopping 
    criterion~\eqref{BMP:eq:ailfem:linstopping} for the linearization 
    has the same interpretation as for Algorithm~\ref{BMP:algorithm:afem}. We numerically equilibrate the 
    linearization error 
    \(\dist^2(u_\ell^{\star}, u_\ell^{k,\jj})^{1/2} 
    \lesssim \dist^2(u_\ell^{k,\jj}, u_\ell^{k-1,\jj})^{1/2}\)
    with the discretization error 
    \(\vvvert u^{\star} - u_\ell^\star \vvvert 
    \lesssim \eta_\ell(u_\ell^{\star}) \) measured in terms of 
    \(\eta_\ell(u_\ell^{k,\jj})\).
    \item While the innermost \(j\)-loop of the algebraic solver 
    is known to terminate (see Lemma~\ref{BMP:lemma:unif-steps-ailfem} below), the \(k\)-loop of the linearization can fail. However, in 
    this case Theorem~\ref{theorem:linearconv} predicts 
    \(u^\star =u_\ell^\star \) and 
    \(\eta_\ell(u_\ell^{k,\jj}) 
    \rightarrow \eta_\ell(u_\ell^{\star}) = 0\) as \(k \rightarrow \infty\).
    \end{enumerate}
\end{remark}

An analogous notation applies for the totally ordered index
set \(\QQ \, \subset \N_0^3\) as in 
Subsection~\ref{BMP:sec:aisfem:algorithm}.
The corresponding quasi-error reads
\begin{equation}
    \label{BMP:eq:energy_difference}
    \Eta_\ell^{k,j}
    \coloneqq
    \vvvert u_\ell^\star - u_\ell^{k,j} \vvvert
    + \vvvert u_\ell^{k,\star} - u_\ell^{k,j} \vvvert
    + \eta_\ell(u_\ell^{k,j}).
\end{equation}

A first important observation is that the innermost \(j\)-loop is known to terminate. 
Unlike Algorithm~\ref{BMP:algorithm:aisfem} (see 
Lemma~\ref{BMP:lemma:unif-steps-aisfem}), the energy structure exploited 
by the novel stopping criterion of Algorithm~\ref{BMP:algorithm:ailfem} even 
allows to conclude that the number of algebraic solver steps per linearization 
step is uniformly bounded.

\begin{lemma}[uniform bound on number of algebraic solvers steps {\cite[Prop.~4]{mps2024}}]
    \label{BMP:lemma:unif-steps-ailfem}
    Consider arbitrary parameters $0 < \theta \le 1$, $0 < \lambdalin $, 
    $0 < \rho < 1$,
    \(\alpha_{\min} > 0\), \(J_{\max} \in \N\). Then, 
    Algorithm~\ref{BMP:algorithm:ailfem}
    guarantees the existence of $j_0 \in \N$ such that $\jj[\ell,k] \le j_0$
    for all $(\ell,k,0) \in \QQ$, i.e., the number of algebraic solver steps is finite and even uniformly bounded.
    \qed
\end{lemma}

The following proposition provides a-posteriori error control at each step 
of the adaptive algorithm, which can again be used to terminate 
Algorithm~\ref{BMP:algorithm:ailfem}.
\begin{lemma}[a-posteriori error control
    {\cite[Prop.~5]{mps2024}}]
	\label{BMP:prop:aposteriori-jk-ailfem}
    For
	any $(\ell,k,j) \in \QQ$ with $k \ge 1$ and $j \ge 1$, it holds that
	\begin{equation}\label{BMP:eq:aposteriori-jk-ailfem}
		\vvvert u^\star - u_\ell^{k,j}\vvvert
		\lesssim \big[ \eta_\ell(u_\ell^{k,j}) 
        + \vvvert u_\ell^{k,j} - u_\ell^{k-1,\jj}\vvvert
        + \vvvert u_\ell^{k,j} - u_\ell^{k,j-1}\vvvert \big]
	\end{equation}
    and for \( k = \kk[\ell]\), \( j = \jj[\ell,\kk]\)
    \[
		\vvvert u^\star - u_\ell^{\kk,\jj}\vvvert
		\lesssim  \eta_\ell(u_\ell^{\kk,\jj}).
        \qed
    \]
\end{lemma}

\subsection{Convergence and complexity}

By definition, the energy difference from~\eqref{BMP:eq:energy_difference}
satisfies
\[
    \dist^2(u_\ell^{k,\jj}, u_\ell^{k-1,\jj})
    =
    \dist^2(u_\ell^\star, u_\ell^{k-1,\jj})
    -
    \dist^2(u_\ell^\star, u_\ell^{k,\jj}),
\]
where all terms are guaranteed to be non-negative.
This replaces the role of the Pythago\-rean
identity~\eqref{BMP:eq:Pythagoras} in the energy-based proof
in~\cite{mps2024} to establish the following convergence results
for Algorithm~\ref{BMP:algorithm:ailfem}.
Importantly, the theorems state parameter robust 
full R-linear convergence and, for sufficiently small adaptivity parameters, 
optimal complexity of the adaptive algorithm.

\begin{theorem}[parameter-robust full R-linear convergence {\cite[Thm.~6]{mps2024}}]
	\label{theorem:linearconv}
	Consider arbitrary parameters $0 < \theta \le 1$, 
	$0 < \lambdalin$, and $0 < \rho < 1$ 
    as well as initialization 
    values $\alpha_{\textup{min}} > 0$, 
    $J_{\textup{max}} \in \N$.
	Then, Algorithm~\ref{BMP:algorithm:ailfem} guarantees
	the existence of constants $0 <\qlin< 1$ and $\Clin > 0$ such that
    \[
        \Eta_\ell^{k,j}
        \le
        \Clin\,
        \qlin^{|\ell,k,j| - |\ell',k',j'|} \,
        \Eta_{\ell'}^{k',j'}
        \quad\text{for all } 
        (\ell,k,j), (\ell',k',j') \in \QQ
        \text{ with } 
        |\ell,k,j| > |\ell',k',j'|.
    \]
    In particular, this yields parameter-robust convergence 
    \begin{align*}
        \vvvert u^\star - u_\ell^{k, j} \vvvert
        \lesssim \Eta_\ell^{k, j}
        \le
        \Clin \,
        \qlin^{
            \vert \ell, k, j \vert
        } \,
        \Eta_0^{0,0} \rightarrow 0 \quad \text{as }
        \vert \ell, k, j \vert \rightarrow \infty.
        \qed
    \end{align*}
\end{theorem}

\begin{theorem}[optimal complexity {\cite[Thm.~10]{mps2024}}]
    \label{BMP:thm:ailfem:optimal_complexity}
    There exist upper bounds 
    \(0 < \theta^\star \leq 1\) and \(\lambdalin^\star > 0\)
    such that, for all
    \(0 < \theta \leq \theta^\star\) and
    \(0 < \lambdalin \leq \lambdalin^\star\),
    the following holds:
	Algorithm~\ref{BMP:algorithm:ailfem} guarantees that
	\[
		\Vert u^\star \Vert_{\A_s }
		\lesssim
		\sup_{(\ell,k, j) \in \QQ}
        \cost(\ell, k, j)^s
		\Eta_\ell^{k,j}
        \lesssim
		\max\{ \Vert u^\star \Vert_{\A_s}, \, \Eta_{0}^{0,0}\}
        \text{ for all } s > 0.
        \qed
	\]
\end{theorem}

\subsection{Numerical example}
\label{BMP:sec:ailfem:example}

\begin{figure}
    \centering
    \subfloat[Adaptively refined mesh on level \(\ell = 9\) with \(5\,246\) triangles.]{%
        \label{BMP:fig:ailfem:Mesh}
        \begin{tikzpicture}
    \begin{axis}[%
        axis equal image,%
        width=4.6cm,%
        xmin=-1.15, xmax=1.15,%
        ymin=-1.15, ymax=1.15,%
        font=\footnotesize%
    ]
        \addplot graphics [xmin=-1, xmax=1, ymin=-1, ymax=1]
        {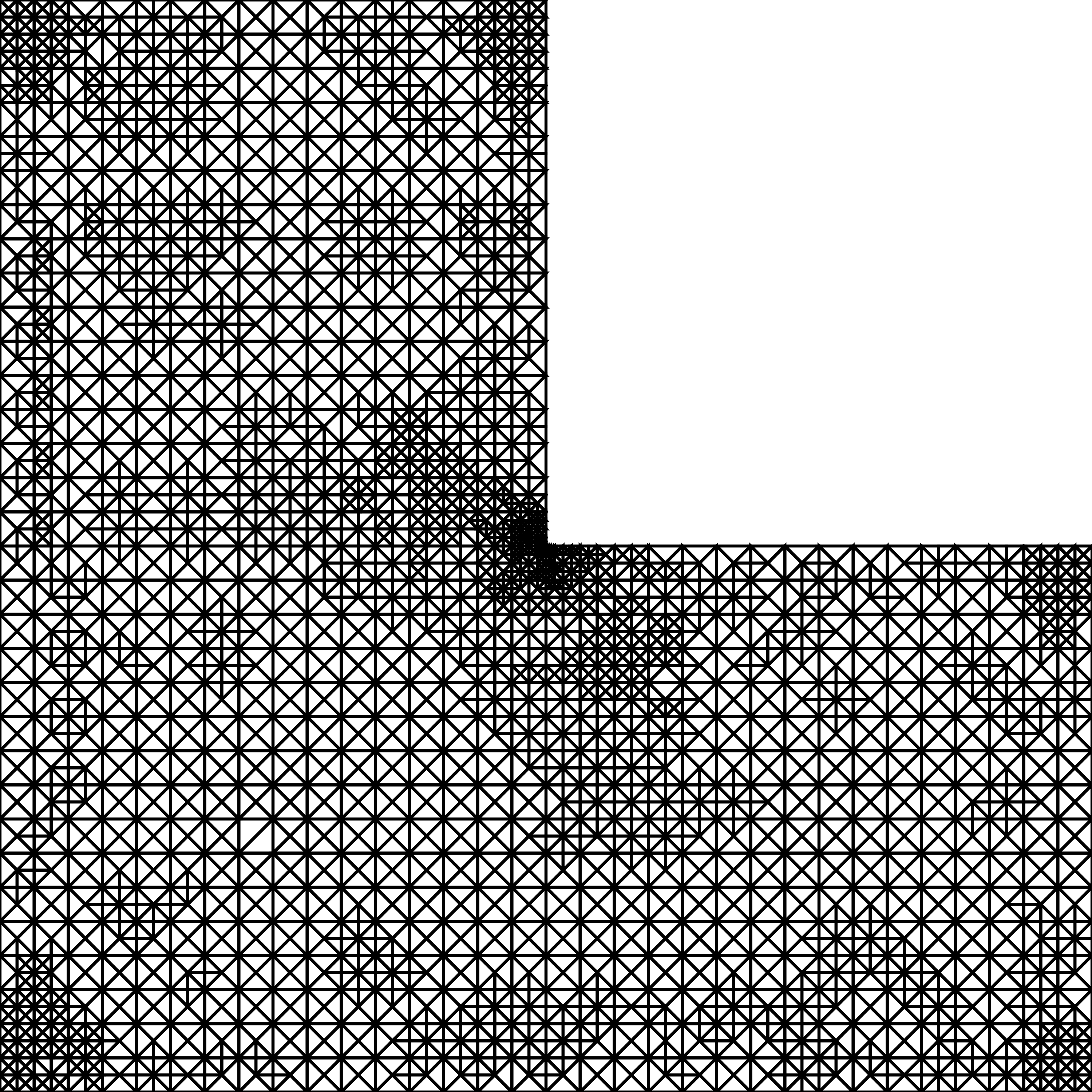};
    \end{axis}
\end{tikzpicture}
    }
    \hfil
    \subfloat[Discrete solution \(u_{9}^{\kk,\jj} \in \XX_{9}\).]{%
        \label{BMP:fig:ailfem:Solution}
        \begin{tikzpicture}
    \pgfplotsset{/pgf/number format/fixed}
    \begin{axis}[%
        width=4.6cm,%
        xmin=-1.1, xmax=1.1,%
        ymin=-1.1, ymax=1.1,%
        zmin=-0.005, zmax=0.05,%
        font=\footnotesize,%
    ]
        \addplot3 graphics [%
                points={%
                    (-1,-1,0) => (0.4,279.7-191.8)
                    (1,-1,0) => (220.1,279.7-279.2)
                    (0,1,0) => (371.6,279.7-161.9)
                    (-0.3125,-0.3125,0.0431357) => (165.8,279.7-6.5)
                }%
            ]
            {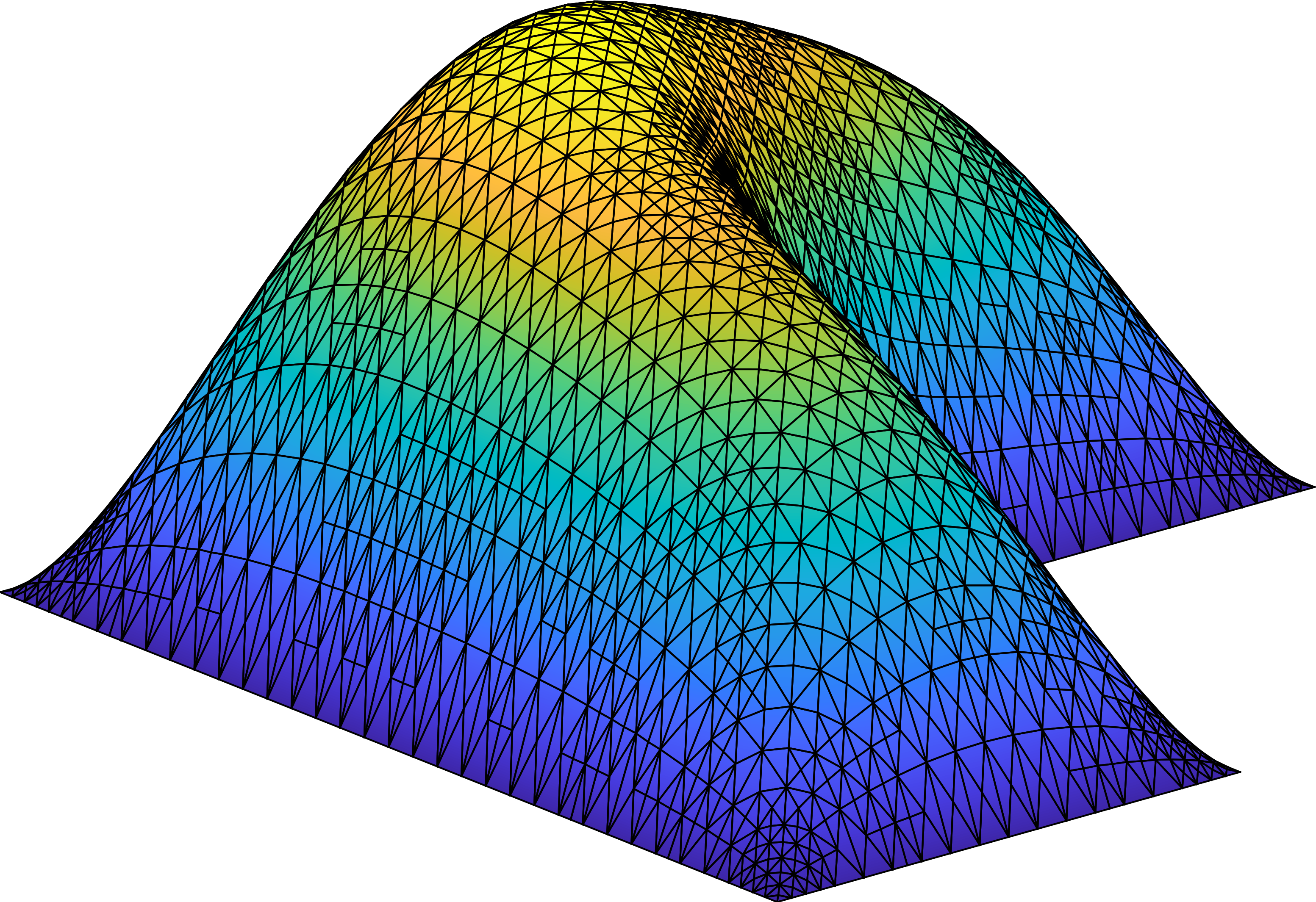};
    \end{axis}
\end{tikzpicture}
    }
    \caption{%
        Mesh and corresponding discrete solution
        for the benchmark problem from
        Subsection~\ref{BMP:sec:ailfem:example}.
        The results are generated by Algorithm~\ref{BMP:algorithm:ailfem} with polynomial degree \(p=1\),
        bulk parameter \(\theta = 0.3\), and stopping parameter \(\lambdalin = 0.7\).
    }
    \label{BMP:fig:ailfem:mesh_solution}
\end{figure}

\begin{figure}
    \centering
    \input{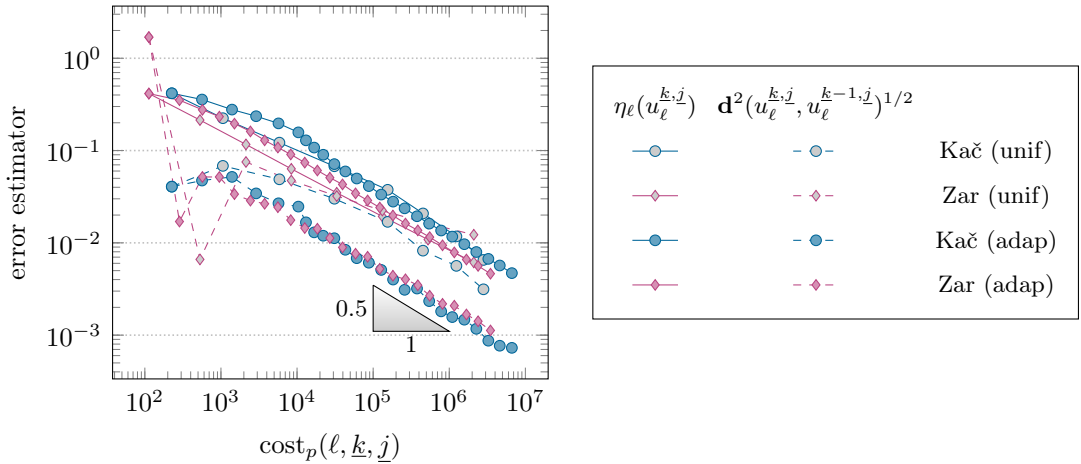}
    \caption{
        Convergence plot of the
        adaptive Algorithm~\ref{BMP:algorithm:aisfem} with Zarantonello 
        linearization and the adaptive Algorithm~\ref{BMP:algorithm:ailfem} 
        with Ka\v{c}anov linearization, with polynomial degree \(p=1\), 
        to solve the benchmark problem from
        Subsection~\ref{BMP:sec:ailfem:example}
        for two different linearizations.
        The damping parameter of the Zarantonello iteration
        reads \(\delta = 1/3\) and the adaptivity parameters
        \(\theta = 0.3\) and \(\lambdalin = \lambdaalg = 0.7\), whereas 
        Ka\v{c}anov linearization uses initial parameters 
        \(\alpha_{\textup{min}} = 100\) and \( J_{\textup{max}} =1\).
        The graphs display values of
        the residual-based error estimator
        \(\eta_\ell (u_\ell^{\kk,\jj})\) from~\eqref{BMP:eq:ailfem:estimator}
        for the final iterates on the level \(\ell \in \N_0\)
        and the energy differences
        \( \dist^2(u_{\ell}^{\kk,\jj},u_{\ell}^{\kk-1,\jj})^{1/2}\).
    }
    \label{BMP:fig:ailfem:linearization}
\end{figure}

\begin{figure}
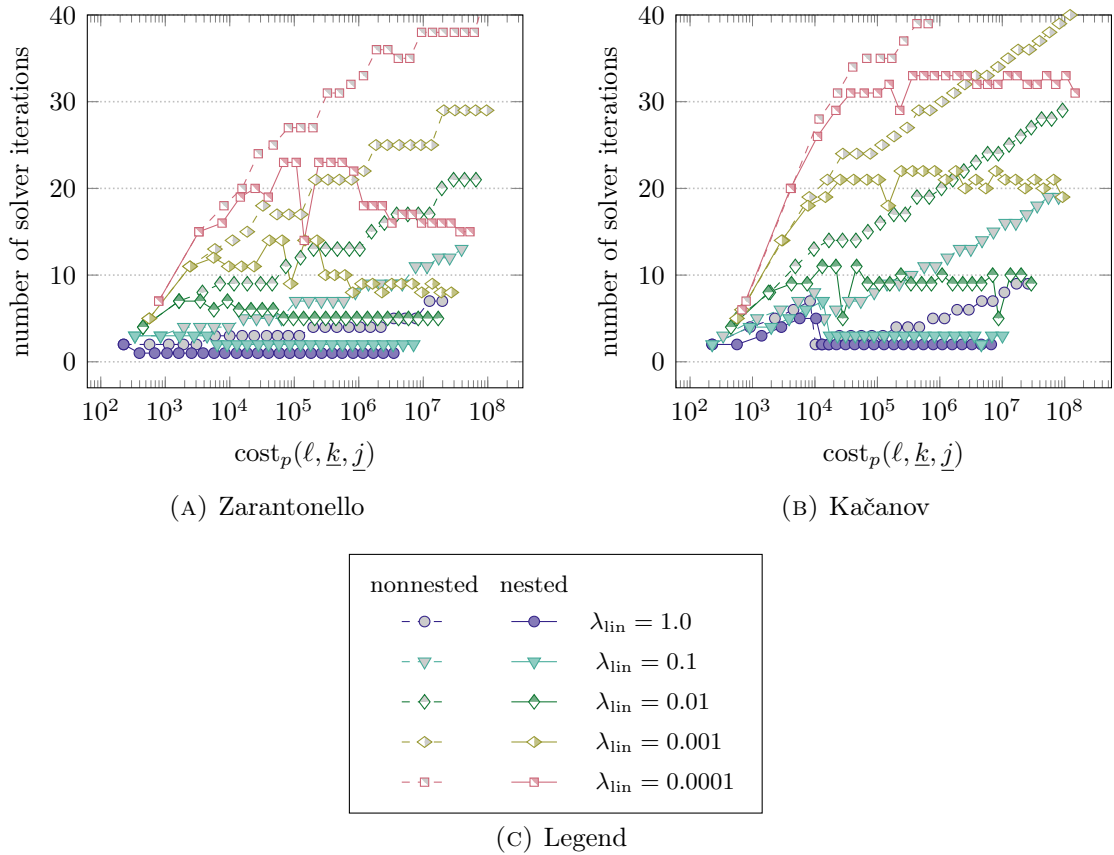

    \centering
    \subfloat[Zarantonello]{%
        \label{BMP:fig:ailfem:iterations:Zarantonello}
        \input{chapter_afem_BMP/figures/Fig19a_AILFEM_Zarantonello_iterations.tex}
    }
    \hfil
    \subfloat[Ka\v{c}anov]{%
        \label{BMP:fig:ailfem:iterations:Kacanov}
        \input{chapter_afem_BMP/figures/Fig19b_AILFEM_Kacanov_iterations.tex}
    }

    \subfloat[Legend]{%
        \label{BMP:fig:ailfem:iterations:legend}
        \begin{tikzpicture}
    \definecolor{col1}{HTML}{332288}
    \definecolor{col2}{HTML}{88CCEE}
    \definecolor{col3}{HTML}{44AA99}
    \definecolor{col4}{HTML}{117733}
    \definecolor{col5}{HTML}{999933}
    \definecolor{col6}{HTML}{DDCC77}
    \definecolor{col7}{HTML}{CC6677}
    \definecolor{col8}{HTML}{882255}
    \definecolor{col9}{HTML}{AA4499}
    \pgfplotsset{%
        degdefault/.style = {%
            mark = *,%
            mark size = 2pt,%
            every mark/.append style = {solid},%
            gray,%
            every mark/.append style = {fill = gray!60!white}%
        },%
        marker1/.style = {%
            degdefault,%
            col1,%
            every mark/.append style = {fill = col1!60!white}%
        },%
        marker2/.style = {%
            degdefault,%
            mark = triangle*,%
            mark size = 2.75pt,%
            every mark/.append style = {rotate = 180},
            col3,%
            every mark/.append style = {fill = col3!60!white}%
        },%
        marker3/.style = {%
            degdefault,%
            mark = halfdiamond*,%
            every mark/.append style = {rotate = 180},%
            mark size = 2.75pt,%
            col4,%
            every mark/.append style = {fill = col4!60!white}%
        },%
        marker4/.style = {%
            degdefault,%
            mark = halfdiamond*,%
            every mark/.append style = {rotate = 90},%
            mark size = 2.75pt,%
            col5,%
            every mark/.append style = {fill = col5!60!white}%
        },%
        marker5/.style = {%
            degdefault,%
            mark = halfsquare*,%
            every mark/.append style = {rotate = 135},%
            mark size = 2.2pt,%
            col7,%
            every mark/.append style = {fill = col7!60!white}%
        },%
        marker6/.style = {%
            degdefault,%
            mark = halfsquare*,%
            every mark/.append style = {rotate = 315},%
            mark size = 2.2pt,%
            col8,%
            every mark/.append style = {fill = col8!60!white}%
        },%
        nonnested/.style = {%
            dashed,%
            every mark/.append style = {%
                fill = black!20!white
            }%
        },%
        nested/.style = {%
            solid%
        },%
    }
    \matrix [%
        matrix of nodes,%
        anchor = center,%
        font = \scriptsize,%
        draw = black,%
        fill = white,%
    ] at (legend) {
        nonnested & nested &  \\
        \ref*{leg:Zar:iter:nonnested:lamlin1} & \ref*{leg:Zar:iter:nest:lamlin1} & $\lambda_{\text{lin}} = 1.0\phantom{0000}$ \\
        \ref*{leg:Zar:iter:nonnested:lamlin0.1} & \ref*{leg:Zar:iter:nest:lamlin0.1} & $\lambda_{\text{lin}} = 0.1\phantom{000}$ \\
        \ref*{leg:Zar:iter:nonnested:lamlin0.01} & \ref*{leg:Zar:iter:nest:lamlin0.01} & $\lambda_{\text{lin}} = 0.01\phantom{00}$ \\
        \ref*{leg:Zar:iter:nonnested:lamlin0.001} & \ref*{leg:Zar:iter:nest:lamlin0.001} & $\lambda_{\text{lin}} = 0.001\phantom{0}$ \\
        \ref*{leg:Zar:iter:nonnested:lamlin0.0001} & \ref*{leg:Zar:iter:nest:lamlin0.0001} & $\lambda_{\text{lin}} = 0.0001$ \\
    };
\end{tikzpicture}
    }

    \caption{%
        Number of iteration steps
        of algebraic iteration
        for the benchmark problem from
        Subsection~\ref{BMP:sec:ailfem:example}
        with nested and non-nested iteration
        for various stopping parameters \(\lambdalin\).
        The damping parameter of the Zarantonello iteration reads \(\delta = 1/3\) and
        the adaptivity parameters \(\theta = 0.3\) and \( \lambdaalg = 0.7\),
        whereas Ka\v{c}anov linearization additionally uses initial parameters 
        \(\alpha_{\textup{min}} = 100\) and \( J_{\textup{max}} =1\).
    }
    \label{BMP:fig:ailfem:iterations}
\end{figure}

Consider the benchmark problem~\eqref{BMP:eq:weak_nonlinear}
on the L-shaped domain
\(\Omega \coloneqq (-1,1)^2 \setminus [0,1)^2\)
with scalar non-linear coefficient
\(\mu(t) \coloneqq 2 + (1 + t)^{-2}\).
The adaptive algorithm creates meshes with increased
resolution of the reentrant corner
where the exact solution exhibits a singularity.
This is illustrated in
Figure~\ref{BMP:fig:ailfem:mesh_solution}.
Consequently, Figure~\ref{BMP:fig:ailfem:linearization}
demonstrates that the adaptive algorithm converges
with the optimal rate for the error estimator
and the energy difference.
The computation with the Zarantonello linearization
employs an algorithm analogous to
Algorithm~\ref{BMP:algorithm:aisfem} 
whereas the Ka\v{c}anov linearization is utilized
in Algorithm~\ref{BMP:algorithm:ailfem}.
Due to the nested structure of the algorithm
and the \(h\)-robustness of the iterative
algebraic solver, the number of algebraic solver steps
(accumulated over each level \(\ell\))
is uniformly bounded as depicted in
Figures~\ref{BMP:fig:ailfem:iterations}.
It turns out that the novel stopping criterion in Algorithm~\ref{BMP:algorithm:ailfem}
appears more robust with respect to a disadvantageous choice
of the initial iterate for the algebraic solver.

%
%

\section{Further generalizations and comments}

We comment that the results on adaptive algorithms treated in this chapter 
are not restricted to a specific linear or non-linear 
model problem but are rather formulated 
in an abstract setting that covers a wide range of problems. 
The only requirements are the properties of the a-posteriori error
estimator as established in~\cite{cfpp2014} and presented 
by~\eqref{BMP:axiom:stability}--\eqref{BMP:axiom:discrete_reliability} on 
page~\pageref{BMP:axiom:stability}, 
a contraction property as~\eqref{BMP:eq:algebra} on 
page~\pageref{BMP:eq:algebra} for the integrated 
iterative solver, as well as the 
quasi-orthogonality~\eqref{BMP:eq:quasi-orthogonality} on 
page~\pageref{BMP:eq:quasi-orthogonality}. 
The latter is satisfied for any well-posed
saddle-point problem \cite{feischl2022}, but as 
exemplified in Section~\ref{BMP:sec:ailfem} can be proved for certain 
energy minimization problems~\cite{ghps2021,hpw2021,bfmps2023}.
In the following, we briefly comment on potential 
generalizations.

\textbf{More general problems.} 
Note that adaptivity has been successfully applied to a variety of problems, 
in particular, to variational problems, such as the obstacle problem; see, e.g.,
\cite{PagePraetorius2013} for linear convergence
as well as
\cite{CarstensenHu2015} for optimal convergence rates with respect to 
the dimension of the underlying discrete space, provided that the 
variational inequality is solved exactly. 
In practice, the iterative solution of variational inequalities 
typically employs a semi-smooth Newton method as formulated in
\cite{HintermuellerItoKunisch2002}.
However, the super-linear convergence result in~\cite{HintermuellerItoKunisch2002}
holds only for initial iterates sufficiently close to the
exact solution and thus the required contraction property is
only guaranteed locally.
Consequently, the application within our framework of 
optimal complexity is subject to future research.

\textbf{Boundary conditions.} 
While the treatment of Robin or Neumann boundary conditions is rather 
immediate, see, e.g.,~\cite{cfpp2014}, the treatment of 
inhomogeneous Dirichlet boundary conditions 
needs more careful consideration.
Available results rely on the discretization of the Dirichlet data followed by 
the solution of the discrete weak formulation; 
see, e.g.,~\cite{BartelsCarstensenDolzmann2004,afkpp2013,cfpp2014}.
In two dimensions, the nodal interpolation of the boundary
data suffices to ensure optimal convergence rates; see, e.g.,
\cite{fpp2014} for the lowest-order case \(p=1\) and 
\cite{BringmannCarstensen2017} for a generalization to arbitrary order \(p \ge 1\).
It turned out that, in any dimension, each \(H^{1/2}\)-stable
operator such as the Scott-Zhang quasi-interpolation
is an admissible choice for the approximation of
the Dirichlet boundary data \cite{afkpp2013};
see also \cite{cfpp2014}.
Alternatively, the \(L^2\) projection can be used 
\cite{BartelsCarstensenDolzmann2004}. However, beyond the 2D case 
\cite{fpp2014,BringmannCarstensen2017} or the Scott-Zhang 
quasi-interpolation for \(d \ge 2\) 
\cite{cfpp2014}, the state-of-the-art analysis requires an extended marking 
strategy to equibalance the approximation error of the Dirichlet data and the 
discretization error of the weak form \cite{afkpp2013}.

\textbf{Equivalent estimators.}
Finally let us comment on other a-posteriori error estimators 
which were applied in the context of adaptive mesh refinement. 
The work~\cite{KreuzerSiebert2011} considers a variety of a-posteriori 
error estimators (including equilibrated fluxes
\cite{Vejchodsky2004,BraessSchoeberl2008}) for the Poisson model problem and 
proves linear convergence with optimal convergence rates for an adaptive 
algorithm with exact solver.
The technical constraint that not only the marked elements but also all their 
neighbors must be refined, has been removed in \cite[Sect.~8]{cfpp2014}, where 
the standard adaptive algorithm with exact solver (exemplified for 
Zienkiewicz-Zhu-type estimators \cite{ZienkiewiczZhu1987}) is considered and 
analyzed: 
The key idea (for linear elliptic PDEs) is the observation that the usual 
estimators are locally equivalent, i.e., the use of D\"orfler marking for one 
of the usual a-posteriori error estimators also yields implicitly 
the D\"orfler marking for the standard residual a-posteriori error 
estimator (with a different marking parameter though).
Then the usual analysis for the residual a-posteriori error estimator applies 
and proves linear convergence with optimal 
rates~\cite{KreuzerSiebert2011,cfpp2014}.
For linear PDEs, this idea extends to the adaptive algorithm with contractive 
solver, while the extension to non-linear PDEs remains for future work.


\sloppy
\printbibliography

\end{document}